\documentclass[openright,a4paper,american,11pt]{article}

\usepackage[latin1]{inputenc}

\usepackage{amsmath}
\usepackage{amssymb}

\usepackage{babel}

\usepackage{amsfonts}
\input xy
\xyoption{all}
\usepackage[dvips]{graphicx}
\usepackage{boxedminipage}
\usepackage{color}

\newtheorem{thm}{Theorem}[section]
\newtheorem{defi}[thm]{Definition}
\newtheorem{prop}[thm]{Proposition}
\newtheorem{lemma}[thm]{Lemma}
\newtheorem{cor}[thm]{Corollary}
\newcommand{\findemo}{\hfill
                     $\Box$ \vspace{1.5 ex}}

\begin{document}
\setlength{\parindent}{0mm}
\title{SYMMETRIC PLANE CURVES OF DEGREE  $7$~: PSEUDOHOLOMORPHIC AND ALGEBRAIC CLASSIFICATIONS}
\author{Erwan Brugallé}
\date{}
\maketitle
\begin{abstract}
This paper is motivated by the real symplectic isotopy problem : does
there exist a nonsingular real pseudoholomorphic curve not 
isotopic in the projective plane to any real algebraic curve of the same
degree? 
Here, we focus our study on symmetric real curves on the
projective plane. We give a classification of real schemes (resp.
complex schemes) realizable by symmetric real curves  of degree $7$
with respect to the type of the curve (resp. $M$-symmetric real curves
of degree $7$). 
In particular, we exhibit two real schemes which are realizable by
real symmetric dividing pseudoholomorphic curves of degree $7$ on the
projective plane but not by algebraic ones. 
\end{abstract} 
\tableofcontents

\section{Introduction and statements of results}

\subsection{Brief History}
The origin of the topology of real algebraic curves can be traced back to
the paper of A. Harnack \cite{Har} which was published in 1876.
D. Hilbert formulated the following question in the 16-th problem of his famous list (see \cite{Hil}):
up to ambient isotopy of
$\mathbb RP^2$, how can be arranged the connected components of the
real point set of a nonsingular real algebraic curve of a given degree
$m$ on $\mathbb RP^2$?

In 1900, the answer was known only for $m\le5$ and nowadays the problem is also solved
for $m=6$ (D. A. Gudkov 1969, \cite{Gud}) and $m=7$ (O. Ya. Viro
1979).
\begin{thm}[Viro, \cite{V1}, \cite{V3}]\label{class7}
 Any nonsingular real algebraic curve of degree $7$ on $\mathbb RP^2$ has one of the following real schemes~:
\begin{itemize}
\item $\langle J\amalg \alpha\amalg 1\langle\beta\rangle\rangle$
 with $\alpha+\beta\le 14$, $0\le\alpha\le 13$, $1\le\beta\le 13$,
\item $\langle J\amalg \alpha\rangle$ with $0\le\alpha\le 15$,
\item $\langle J\amalg 1\langle 1 \langle 1\rangle\rangle\rangle.$
\end{itemize}
Moreover, any of these 121 real schemes is realizable by nonsingular real algebraic  curves of degree $7$ on $\mathbb RP^2$.
\end{thm}
\textbf{Remark. }The notations used to encode real schemes are the
usual ones introduced in \cite{V3}. For example,
$\langle J\amalg \alpha\amalg 1\langle\beta\rangle\rangle$ means a non-contractible (in $\mathbb RP^2$) component,
an oval with $\beta$ ovals in its interior and $\alpha$ ovals in its exterior. All the ovals $\alpha$ and $\beta$ lie outside each other. An oval is a contractible (in $\mathbb RP^2$) component of the real curve.

\vspace{2ex}
In the study of the topology of real algebraic varieties, one can distinguish two directions~: \textit{constructions} and \textit{prohibitions}.
Historically, the first prohibition result, the Harnack
theorem,
asserts that the real part of a real algebraic curve of genus $g$ cannot have more than
$g+1$ connected components (curves with the maximal number of connected components are called \textit{$M$-curves}).
In \cite{Arn}, V. I. Arnold paved the way to the use of powerful
topological methods in the study of this problem.
Almost all known prohibitions
can be obtained via a topological study of the double covering of $\mathbb CP^2$ branched along the complex point set of  a curve of even degree
or looking at the braid defined by the intersection of a curve and a certain sphere $S^3$ in $\mathbb CP^2$ (the interested reader can see the surveys \cite{DK}, \cite{V3} and \cite{Wil}).
In particular, almost all known restrictions on the topology of real algebraic curves are also valid for a
category of  more flexible objects~: \textit{pseudoholomorphic curves}. These objects, introduced by M. Gromov to study symplectic $4-$manifolds
in \cite{Gro}, share a lot of properties with the algebraic curves (for example,
the Harnack theorem is still valid in the real
pseudoholomorphic case) and are much easier to deal with. Indeed, there are methods to construct pseudoholomorphic curves
which are not necessarily algebraic~: S. Yu. Orevkov (see \cite{O1}) proved that it is sufficient to show the quasipositivity of some braid,
and I. Itenberg and E. Shustin (see \cite{IS}) proved that any T-construction gives a pseudoholomorphic curve
 even if
the convexity condition on the triangulation is not fulfilled.
Up to now, the classification up to isotopy of real pseudoholomorphic $M$-curves
of degree $8$ on $\mathbb RP^2$ is achieved by Orevkov \cite{O2}. However,
it remains $6$ real schemes for which it is unknown whether they are
realizable by algebraic $M$-curves of degree $8$ on $\mathbb RP^2$.

Then a natural question arises~: does there exist a real scheme on
$\mathbb RP^2$ which would be realizable pseudoholomorphically but not
algebraically? This problem is the real counterpart of the symplectic
isotopy problem.
Strictly speaking (i.e., dealing with nonsingular
curves on the projective plane), this  problem is still open.
However,
dealing with
singular curves or curves in ruled surfaces up to fiberwise isotopy, the answer is \textit{yes}~: S. Fiedler Le Touz\'e and S. Orevkov
(see \cite{O2}, \cite{FTO}) exhibit mutual arrangements of two
nonsingular curves on $\mathbb RP^2$ which are realizable pseudoholomorphically but not
algebraically (examples with many irreducible components are easy to construct, see the ``pseudoholomorphic Pappus 
theorem'' in
\cite{FTO}), and
Orevkov and Shustin (see \cite{OS1} and \cite{OS2}) exhibit
nonsingular \textit{$\mathcal L$-schemes} on the second rational
geometrically ruled surface $\Sigma_2$ which are realizable 
pseudoholomorphically but not algebraically. One can note that
in our definition of pseudoholomorphic curves on ruled surfaces,
(which is the same than 
Orevkov and Shustin's one), we consider only almost complex structures
on $\Sigma_n$
for which the exceptional divisor (if any) is pseudoholomorphic (see
section \ref{prelim}). Forgetting this condition, 
 J-Y. Welschinger constructed 
in \cite{W} examples of real pseudoholomorphic curves on $\Sigma_n$ for $n\ge2$
which are not isotopic to 
any real algebraic curve realizing the same homology class.

However in the case of
nonsingular curves on $\mathbb RP^2$, for each known algebraic classification, the pseudoholomorphic classification is the same
(and even the proofs for both classifications are alike!).
Thus, Theorem \ref{class7} is still true replacing "algebraic" by "pseudoholomorphic".
The real symplectic isotopy problem
 turns out to be difficult. So, one can tackle a simpler question, looking for example at complex curves
which admit more symmetries
than the action of $\mathbb Z/2\mathbb Z$ given by the complex conjugation.
The first natural action is an additional holomorphic
action of $\mathbb Z/2\mathbb Z$, which can be given by a symmetry of
the projective plane (see section \ref{general}).
Such a real plane curve, invariant
under
a symmetry, is called a \textit{symmetric curve}. The systematic study of symmetric curves was started by T. Fiedler (\cite{Fie}) and continued
by S. Trille (\cite{Tr}, \cite{Tr2}).
The \textit{rigid isotopy classes} (two nonsingular curves of degree $m$ on $\mathbb RP^2$ are said to be rigidly isotopic if
they belong to the same connected component
of the complement of the discriminant hypersurface in the space of curves of degree $m$) of nonsingular sextics
in $\mathbb RP^2$ which contain a symmetric curve can be obtained from \cite{It1}.
Recently (see \cite{It}), using auxiliary conics, I. Itenberg and V. Itenberg found an elementary proof of this classification.
Once again, algebraic and pseudoholomorphic classifications coincide.
On the other hand, in \cite{OS2}, Orevkov and Shustin
showed that there exists a  real scheme which is realizable by nonsingular symmetric pseudoholomorphic $M$-curves of
degree $8$, but which is not realizable by real symmetric algebraic curves of degree $8$.

Hence, it is natural to wonder about the degree $7$ and this is the subject of this paper. It turns out that  the classification of real schemes which are realizable by
 nonsingular symmetric curves of degree $7$ on $\mathbb RP^2$
are the same in both algebraic and pseudoholomorphic cases, as well as  the classification of \textit{complex schemes} which are realizable by
 nonsingular symmetric $M$-curves of degree $7$ on $\mathbb RP^2$  (Theorem
\ref{Main thm1} and Corollary \ref{Main thm2}). However, if we look at real schemes which are realizable by \textit{nonsingular dividing
symmetric curves} of degree $7$ on $\mathbb RP^2$, the answers are different. In Theorems \ref{Main thm3} and \ref{Main thm4}, 
we exhibit two real schemes which are realizable by real symmetric dividing pseudoholomorphic curves of degree $7$ on $\mathbb RP^2$ but not
by algebraic ones.

We state our classification results in the next subsection. In section
\ref{prelim}, we give some definitions and properties of the objects
used in this paper (\textit{rational geometrically ruled surfaces,
  braid  associated to an $\mathcal L$-scheme}). 
The algebraic prohibitions are obtained by means of \textit{real
  trigonal graphs} and \textit{combs} associated to an $\mathcal
L$-scheme. We present these objects and their link with real algebraic
trigonal curves in section \ref{comb theo}.
We  give there  an efficient algorithm to deal with combs.
In section \ref{pseudo state}, we give results related to the pseudoholomorphic category. Algebraic results are given in section \ref{algebraic state}.

\textbf{Acknowledgments : }I am very grateful to I. Itenberg for
advices and encouragements, to S. Orevkov for useful discussions and
to T. Fiedler for encouragements. 

\subsection{Classification results}\label{statement}
First we state two simple lemmas.
Surprisingly, we did not succeed to find these statements in the literature.
The proof of these prohibitions is a straightforward application of the Bezout
theorem, the Fiedler orientation alternating rule (see \cite{V3}) and the Rokhlin-Mischachev orientation formula (see \cite{Rok} or section \ref{prelim}).
All the constructions are performed in \cite{Soum} and \cite{FLT2}. A
real algebraic (or pseudoholomorphic) curve is said to be a \textit{dividing} curve if the
set of real points of the curve disconnects the set of complex points.
\begin{lemma}\label{class7 div}
 Any nonsingular dividing real pseudoholomorphic curve of degree $7$ on $\mathbb RP^2$ has one of the following real schemes~:
\begin{itemize}
\item $\langle J\amalg \alpha\amalg 1\langle\beta\rangle\rangle$
 with $\alpha+\beta\le 14$, $\alpha+\beta=0$ mod $2$, $0\le\alpha\le 13$, $1\le\beta\le 13$, if $\alpha=0$ then $\beta\ne 2, 6, 8$ and if $\alpha=1$ then $\beta\ge 5$,
\item $\langle J\amalg \alpha\rangle$ with $7\le\alpha\le 15$, $\alpha=1$ mod $2$,
\item $\langle J\amalg 1\amalg 1\langle 1 \langle 1\rangle\rangle\rangle.$
\end{itemize}
Moreover, any of these real schemes is realizable by nonsingular dividing real algebraic  curves of degree $7$ on $\mathbb RP^2$.
\end{lemma}
\begin{lemma}\label{class7 nondiv}
 Any nonsingular non-dividing real pseudoholomorphic curve of degree $7$ on $\mathbb RP^2$ has one of the following real schemes~:
\begin{itemize}
\item $\langle J\amalg \alpha\amalg 1\langle\beta\rangle\rangle$
 with $\alpha+\beta\le 13$, $0\le\alpha\le 12$, $1\le\beta\le 13$,
\item $\langle J\amalg \alpha\rangle$ with $0\le\alpha\le 14.$
\end{itemize}
Moreover, any of these real schemes is realizable by nonsingular non-dividing real algebraic  curves of degree $7$ on $\mathbb RP^2$.
\end{lemma}
We now introduce \textit{symmetric curves}. Denote by $s$ the holomorphic involution of $\mathbb CP^2$ given by
$[x:y:z]\mapsto [x:-y:z]$.
A real curve on $\mathbb RP^2$ is called \textit{symmetric} if $s(A)=A$.

Here we state the five main classifications of this article.
\begin{thm}\label{Main thm1}
The following real schemes are not realizable by nonsingular symmetric real pseudoholomorphic curves of degree $7$  on $\mathbb RP^2$~:
\begin{itemize}
\item $\langle J\amalg (14-\alpha)\amalg 1\langle \alpha\rangle \rangle\textrm{  with } \alpha=6,7,8,9$,
\item $\langle J\amalg (13-\alpha)\amalg 1\langle \alpha\rangle \rangle \textrm{  with } \alpha=6,7,9$.
\end{itemize}
Moreover, any other real scheme realizable by nonsingular real algebraic curves of degree $7$ on $\mathbb RP^2$
is realizable by nonsingular symmetric real algebraic curves of degree $7$ on $\mathbb RP^2$.
\end{thm}
\textit{Proof. }
The pseudoholomorphic prohibitions are proved in Propositions \ref{first M} and \ref{prohib thm 1}. All the other curves are constructed
algebraically in Propositions \ref{main constr},
\ref{change1}, \ref{constr thm1 3},
and
Corollary \ref{constr thm1 1}.~\findemo

If $C$ is a symmetric curve of degree $7$, the quotient curve $C/s$  is a curve of
bidegree $(3,1)$ on $\Sigma_2$ which has a special position with
respect to the base $\{y=0\}$ (see section \ref{general}). 
In the case
of $M$-curves of degree $7$, we give a classification of the possible mutual
arrangements of those two curves.

\begin{figure}[h]
\centering
\begin{tabular}{cccc}
\includegraphics[height=2cm, angle=0]{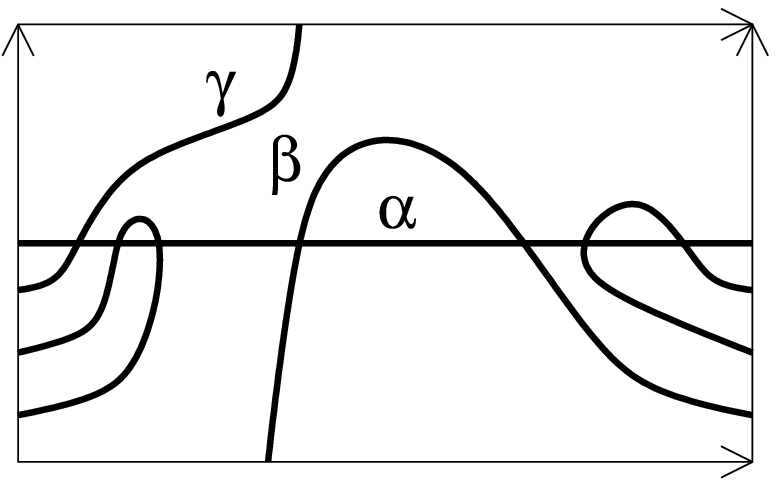}&
\includegraphics[height=2cm, angle=0]{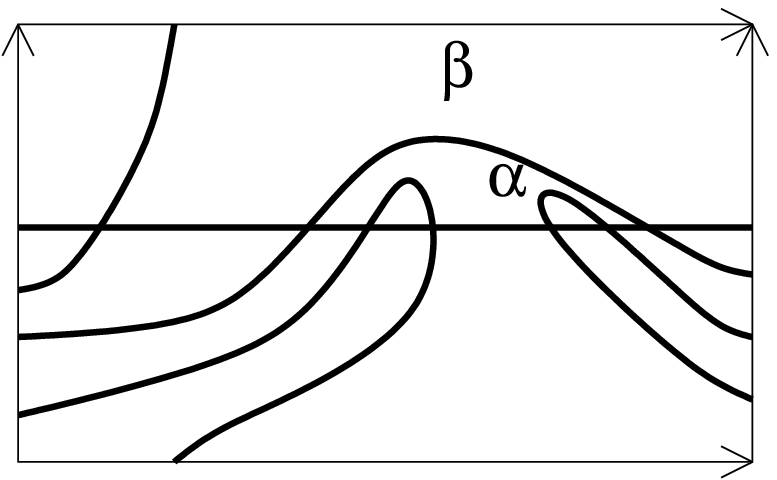}&
\includegraphics[height=2cm, angle=0]{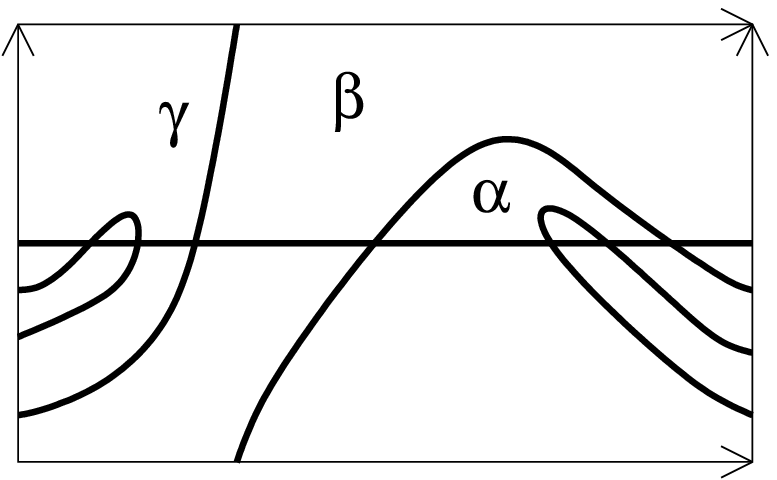}&
\includegraphics[height=2cm, angle=0]{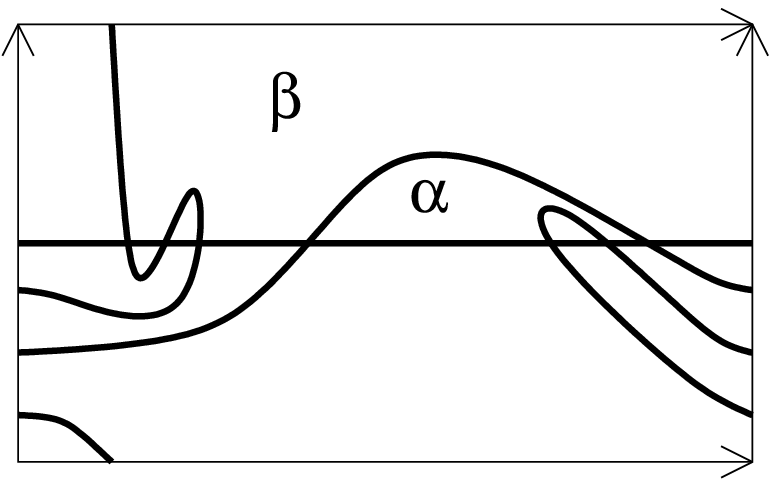}
\\ a)&b)&
c)&d)
\end{tabular}
\caption{}
\label{pr even}
\end{figure}

As explained in section \ref{ruled}, the real part of $\Sigma_2$ is a torus. In Figure \ref{pr even}, the rectangles with parallel edges identified according to the arrows represent $\mathbb R\Sigma_2$. 
The two horizontal edges represent the real part of the exceptional divisor $E$, and the two vertical edges represent the real part of a fiber. An $\mathcal L$-scheme is a collection of embedded circles in $\mathbb R\Sigma_n$ up to isotopy which respects the pencil of lines in $\mathbb R\Sigma_n$(see section \ref{ruled}). In the $\mathcal L$-schemes depicted in Figure \ref{pr even}, the numbers $\alpha, \beta$ and $\gamma$ represent as many ovals lying all outside each other.

\begin{thm}\label{M quotient}
Let $C$ be a nonsingular symmetric pseudoholomorphic $M$-curve of degree $7$ on $\mathbb
RP^2$. Then the $\mathcal L$-scheme realized by
the union of its quotient curve $C/s$ and the base $\{y=0\}$ in
$\Sigma_2$ is one of those
depicted in 
\begin{itemize}
\item Figure \ref{pr even}a) with $(\alpha,\beta,\gamma)=(1,5,0)$,
  $(5,1,0)$, $(0,1,5)$ and $(0,5,1)$,
\item Figures \ref{pr even}b) and d) with $(\alpha,\beta)=(1,5)$ and 
  $(5,1)$,
\item Figure \ref{pr even}c) with $(\alpha,\beta,\gamma)=(2,4,0)$,
  $(6,0,0)$, $(0,0,6)$ and $(0,4,2)$.
\end{itemize}

Moreover, all these arrangements are realizable by the quotient curves
of nonsingular symmetric algebraic
$M$-curves of degree $7$ on $\mathbb 
RP^2$.
\end{thm}
\textit{Proof. }Let $C$ be a symmetric $M$-curve of degree $7$ on
$\mathbb RP^2$. According to Theorem \ref{max sym}, Proposition \ref{inter max},
the Bezout theorem, Lemmas \ref{lem bezout1} and \ref{convexe},
 the only possibilities for the $\mathcal L$-schemes realized by the
 union of the quotient curve of $C$ and the base $\{y=0\}$
  are depicted in Figures \ref{prohib even}a) and c) with
  $\alpha+\beta+\gamma=6$ and in Figures \ref{prohib even}b) and d) with
  $\alpha+\beta=6$. Now the prohibitions follow from Lemmas \ref{pr1},
  \ref{pr2}, \ref{pr3} and \ref{pr4}. All the constructions are
  performed in Proposition \ref{main constr}.~\findemo

\textbf{Remark. }Theorem \ref{M quotient} has an immediate consequence
: if a symmetric $M$-curve of degree $7$ has 
a nest (see section \ref{real curves}) and at least 2 inner ovals, then the $\mathcal L$-scheme realized by
the union of its quotient curve and the base $\{y=0\}$ in $\Sigma_2$
is uniquely 
determined by the real scheme of the symmetric curve. If a symmetric
$M$-curve has a nest and only one inner oval, or has no nest, then
there are two possibilities for the $\mathcal L$-scheme realized by
the union of its quotient curve and the base $\{y=0\}$.

\vspace{2ex}
This latter theorem allows us to classify the complex schemes realized
by symmetric $M$-curves of degree $7$. The notations used to encode complex schemes are the
usual ones proposed in \cite{V3}.
\begin{cor}\label{Main thm2}
A complex scheme is realizable by  nonsingular symmetric real algebraic (or pseudoholomorphic) $M$-curves of degree $7$
on $\mathbb RP^2$  if
and only if it is contained in the following list~:
\begin{center}
\begin{tabular}{lllll}
$\langle J\amalg 9_+\amalg 6_-\rangle_I$
&&
$\langle J\amalg 4_+ \amalg 6_-
\amalg 1_-\langle 3_+ \amalg 1_-\rangle \rangle_I$ 
&&
$\langle J\amalg  2_-
\amalg 1_-\langle 7_+ \amalg 5_-\rangle \rangle_I$
\\
 $\langle J\amalg 7_+ \amalg 6_-
\amalg 1_-\langle 1_+\rangle \rangle_I$
&&
 $\langle J\amalg 5_+\amalg  4_-
\amalg 1_-\langle 3_+\amalg 2_- \rangle \rangle_I$
&&
 $\langle J\amalg 1_+
\amalg 1_-\langle 7_+\amalg 6_- \rangle \rangle_I$
\\
$\langle J\amalg 5_+ \amalg 7_-
\amalg 1_-\langle 2_+\rangle \rangle_I$ 
&&
$\langle J\amalg 1_+ \amalg 3_-
\amalg 1_-\langle 6_+ \amalg 4_-\rangle \rangle_I$
&
&

\\
 $\langle J\amalg 6_+ \amalg 5_-
\amalg 1_+\langle 1_+\amalg 2_-\rangle \rangle_I$
& &
 $\langle J\amalg 2_+ \amalg 1_-
\amalg 1_+\langle 5_+\amalg 6_-\rangle \rangle_I$
&&
\end{tabular}
\end{center}
\end{cor}
\textit{Proof. }If $C$ is a nonsingular dividing symmetric curve of
degree $7$ in 
$\mathbb RP^2$, the $\mathcal L$-scheme realized by the union of its
quotient curve and the base $\{y=0\}$ determines uniquely the complex
orientations 
of the initial  symmetric curve (see section \ref{general}).
Now the corollary follows from Theorem \ref{M quotient}.~\findemo 

Looking at dividing symmetric curves of degree 7, the
pseudoholomorphic and the algebraic classifications differ.

\begin{thm}[Pseudoholomorphic classification]\label{Main thm3}
The following real schemes are not realizable by nonsingular dividing symmetric real pseudoholomorphic curves of degree $7$
on $\mathbb RP^2$~:
 \begin{center}
$\langle J\amalg 2\amalg 1\langle 10\rangle \rangle$, 
$\langle J\amalg 6\amalg 1\langle 6\rangle \rangle$ and 
$\langle J\amalg 4\amalg 1\langle 4\rangle \rangle$.
\end{center}
Moreover, any other real scheme  mentioned in  Lemma \ref{class7 div} and   not forbidden by Theorem \ref{Main thm1}
 is realizable by nonsingular dividing symmetric real pseudoholomorphic curves of degree $7$ on $\mathbb RP^2$; any
real scheme  mentioned in Lemma \ref{class7 nondiv} which is  not forbidden by Theorem \ref{Main thm1}
 is realizable by non-dividing nonsingular symmetric real pseudoholomorphic curves of degree $7$ on $\mathbb RP^2$.
\end{thm}
\textit{Proof. }
All the pseudoholomorphic prohibitions are proved in Propositions \ref{restr thm 4 3} and \ref{non 4 4 I}. All the constructions are done in Propositions \ref{restr thm 4 5},
\ref{restr thm 4 6}, \ref{restr thm 4 7}, \ref{restr thm 4 8} and in Corollary \ref{restr thm 4 4}, and \ref{pseudo non alg 2}.~\findemo

\begin{thm}[Algebraic classification]\label{Main thm4}
The real schemes
$$\langle J\amalg 8\amalg 1\langle 4\rangle \rangle\textrm{ and }\langle J\amalg 4\amalg 1\langle 8\rangle \rangle$$
are not realizable by
a nonsingular dividing symmetric real algebraic curve of degree $7$ on $\mathbb RP^2$.
Any other real scheme which is realizable by nonsingular dividing symmetric real pseudoholomorphic curves of degree $7$ on $\mathbb RP^2$ is realizable by
nonsingular dividing symmetric real algebraic curves of degree $7$ on $\mathbb RP^2$.

Any real scheme which is realizable by  non-dividing nonsingular symmetric real pseudoholomorphic curves of degree $7$ on $\mathbb RP^2$  is  realizable
by non-dividing nonsingular symmetric real algebraic curves of degree $7$ on $\mathbb RP^2$.
\end{thm}
\textit{Proof. }
The two algebraic prohibitions are proved in Propositions \ref{only possibility}, \ref{prop 84} and \ref{prop 48}. All the constructions are done in Propositions \ref{restr thm 4 5},
\ref{restr thm 4 6}, \ref{restr thm 4 7}, \ref{restr thm 4 8}, and in Corollary \ref{restr thm 4 4}.
\findemo

\section{Preliminaries}\label{prelim}

\subsection{Rational geometrically ruled surfaces}\label{ruled}
Let us define \textit{the $n^{th}$ rational geometrically ruled surface}, denoted by $\Sigma_n$, the surface obtained by taking four copies of $\mathbb C^2$ with coordinates $(x,y)$, $(x_2,y_2)$, $(x_3,y_3)$
and $(x_4,y_4)$, and by gluing them along $(\mathbb C^*)^2$ with the identifications $(x_2,y_2)=(1/x,y/x^n)$, $(x_3,y_3)=(x,1/y)$
and $(x_4,y_4)=(1/x,x^n/y)$. Define on $\Sigma_n$ the algebraic curve $E$ (resp. $B$ and $F$) by the equation $\{y_3=0\}$ (resp. $\{y=0\}$ and $\{x=0\}$). The coordinate system $(x,y)$
is called \textit{standard}.
The projection $\pi$ : $(x,y)\mapsto x$ on $\Sigma_n$ defines a $\mathbb CP^1$-bundle over $\mathbb CP^1$.
The intersection numbers of $B$ and $F$ are $B\circ B=n$, $F\circ F=0$ and $B\circ F=1$. The surface $\Sigma_n$ has a natural real structure induced by the complex conjugation on $\mathbb C^2$, and
the real part of $\Sigma_n$ is a torus if $n$ is even and a Klein bottle if $n$ is odd. The restriction of $\pi$ on $\mathbb R\Sigma_n$ defines a pencil of lines denoted by $\mathcal L$.

 The group $H_2(\Sigma_n;\mathbb Z)$ is isomorphic to $\mathbb Z\times\mathbb Z$ and is generated by the classes of $B$ and $F$. Moreover, one has $E=B-nF$. An algebraic or pseudoholomorphic curve on $\Sigma_n$ is
said to be of \textit{bidegree} $(k,l)$ if it realizes the homology class $kB+lF$ in $H_2(\Sigma_n;\mathbb Z)$. 
\begin{defi}
A curve of bidegree $(1,0)$ is called a base.

A curve of bidegree $(3,0)$ is called a trigonal curve.
\end{defi}

In the rational geometrically ruled surfaces, we study real curves up
to \textit{isotopy with respect to $\mathcal L$}. Two curves are said
to be 
isotopic with respect to the fibration $\mathcal L$ if there exists an
isotopy of $\mathbb R\Sigma_n$ which brings the first curve to the second one, and
which maps each line of $\mathcal L$ to another line of $\mathcal L$.

The work of Gromov (see \cite{Gro}) shows that the  following proposition, trivial in the algebraic case, it is also true in the pseudoholomorphic case.
\begin{prop}
There exists a unique base of $\Sigma_2$ which passes through $3$ generic points.
\end{prop}

An arrangement of circles $A$, which may be nodal, in $\mathbb R\Sigma_n$ up to
isotopy of $\mathbb R\Sigma_n\setminus E$ 
which respects the pencil of lines $\mathcal L$ is called
an \textit{$\mathcal L$-scheme}. 

\vspace{1ex}
In this paper, we represent $\mathbb R\Sigma_n$ as a rectangle whose parallel edges are identified according to the arrows. The two horizontal edges represent the real part of the exceptional divisor $E$, and the two vertical edges represent the real part of a fiber. 
In the $\mathcal L$-schemes depicted in this paper, the numbers $\alpha, \beta$ and $\gamma$ represent as many ovals lying all outside each other.

\subsection{Real curves}\label{real curves}
A \textit{real pseudoholomorphic curve} $C$ on $\mathbb CP^2$ or
$\Sigma_n$ is an immersed Riemann surface which is a $J$-holomorphic
curve in some tame almost complex structure $J$ such that the
exceptional section (in $\Sigma_n$ with $n\ge 1$) is $J$-holomorphic
(see \cite{Gro}),
$conj(C)=C$, and $conj_*\circ J_p=-J_p\circ conj_*$, where $conj$ is
the complex conjugation and $p$ is any point of $C$. If the curve $C$ is immersed in $\mathbb CP^2$ (resp. $\Sigma_n$) and realizes the homology class $d[\mathbb CP^1]$ in $H_2(\mathbb CP^2;\mathbb Z)$ (resp. $kB+lF$ in $H_2(\Sigma_n;\mathbb Z)$), $C$ is said to be of degree $d$ (resp. of bidegree $(k,l)$).
For example, a real algebraic curve  of degree $d$ on $\mathbb CP^2$ is a real pseudoholomorphic curve of degree $d$ for the standard complex structure. All intersections of two $J$-holomorphic curves are positive, so the Bezout theorem is still true
for two $J$-holomorphic curves.

Gromov pointed out that given $2$ (resp. 5) given generic points in $\mathbb CP^2$,
there exists a unique $J$-line (resp. $J$-conic) passing through these points.
As soon as the degree is greater than $3$, such a statement is no more true. This is a direction to find
some difference between algebraic and pseudoholomorphic curves (see \cite{FTO} where the authors use pencils of cubics).

From now on, we state results about "curves" if they are valid in both
cases, algebraic and pseudoholomorphic. A nonsingular real curve of
genus $g$ is called an \textit{$(M-i)$-curve} 
if its real part has $g+1-i$ connected components. If $i=0$, we simply
speak about an \textit{$M$-curve} or \textit{maximal curve}.
A connected component of the real part of a nonsingular real curve on
a real surface is called an \textit{oval} if it is contractible in the
surface, and is called a \textit{pseudo-line} otherwise. 
In $\mathbb RP^2$ the complement of an oval is formed by two connected
components, one of which is homeomorphic to a disk  
(called the \textit{interior} of the oval) and the other to a M\"obius strip (called the \textit{exterior} of the oval).
Two ovals in $\mathbb RP^2$ are said to constitute an \textit{injective pair} if one of them is enclosed by the other. A set of ovals such that any two ovals of this set form an injective pair is called a nest.
For a nonsingular real curve of degree $7$ with a nest of depth $2$, the oval which contains some other ovals is called the
 \textit{non-empty oval} of the curve. The ovals lying
inside this oval are called the \textit{inner ovals} while those lying outside are called the
\textit{outer ovals}.

A nonsingular
real curve $C$ is a $2$-dimensional manifold and $C\setminus\mathbb RC$ is either connected or it has two connected components. In the former
case, we say that $\mathbb RC$ is a \textit{non dividing curve},
or of type $II$, and in the latter case, we say that $\mathbb RC$ is a \textit{dividing curve}, or of type $I$.

Let $C$ be a dividing curve.
The two halves of $C\setminus\mathbb RC$ induce two opposite
orientations on $\mathbb RC$ which are called \textit{complex
  orientations} 
of the curve. Fix such a complex orientation of $\mathbb RC$.
An injective pair of ovals in $\mathbb RP^2$ is called \textit{positive} if the orientations of the two ovals are induced by one of the orientations
of the annulus in $\mathbb RP^2$
bounded by the two ovals, and \textit{negative} otherwise. Now suppose that $C$ is of odd degree. If the integral homology classes realized by the odd component of the curve and an oval in the
M\"obius strip defined by the exterior of
this oval have the same sign, we say that this oval is \textit{negative}, and
\textit{positive} otherwise.
Denote by
$\Pi_+$ (resp. $\Pi_-$) the number of positive (resp. negative)  injective pairs of ovals and by
$\Lambda_+$ (resp. $\Lambda_-$) the number of positive (resp. negative) ovals.
\begin{prop}[Rokhlin-Mischachev's orientation formula, \cite{Rok}, \cite{Mis}]\label{orientation formula}
If $C$ is a dividing nonsingular curve of degree $2k+1$ on $\mathbb RP^2$ with $l$ ovals, then
$$\Lambda_+-\Lambda_-+2(\Pi_+-\Pi_-)=l-k(k+1)$$
\end{prop}
In this paper, we will use the the Fiedler orientations alternating rule  (see \cite{V3}, \cite{Fie2}) to determine complex orientations of dividing curves in $\Sigma_n$.

\vspace{1ex}
The following fact about curves on toric surfaces is well known (see,
for example, \cite{Ful}).
\begin{prop}\label{lattice}
Let $C$ be a nonsingular real curve with Newton polygon $\Delta$ on a toric surface. Then the genus of $C$ is equal to the number of
integer points in the interior of $\Delta$.
\end{prop}

\begin{figure}[h]
\centering
\begin{tabular}{cc}
\includegraphics[height=2.5cm, angle=0]{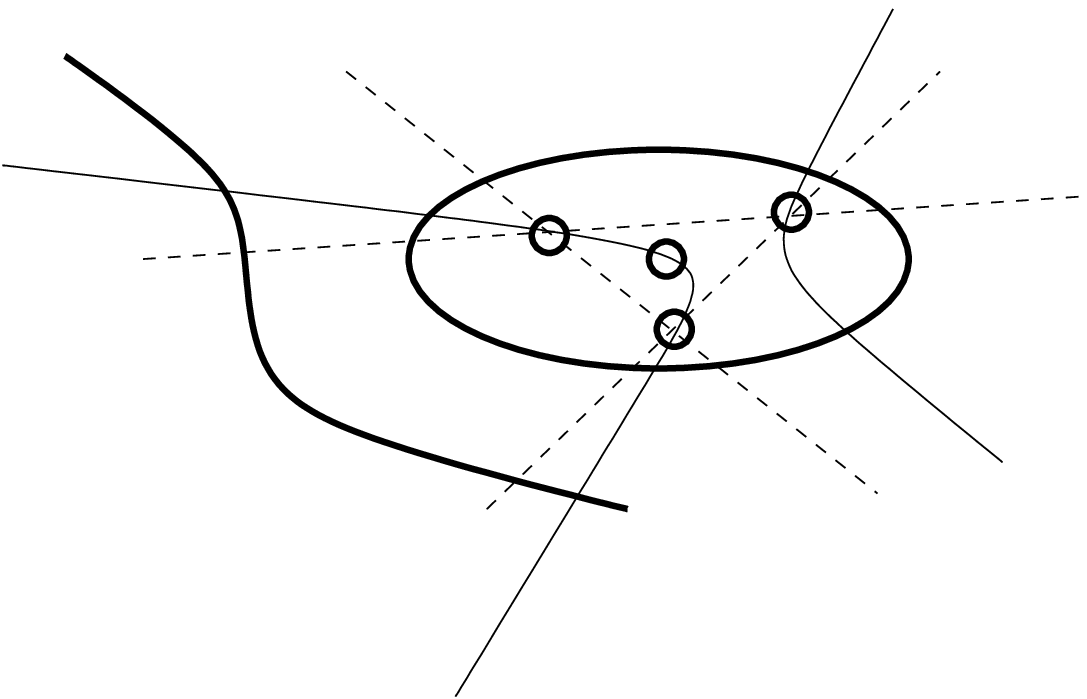}&
\includegraphics[height=2.5cm, angle=0]{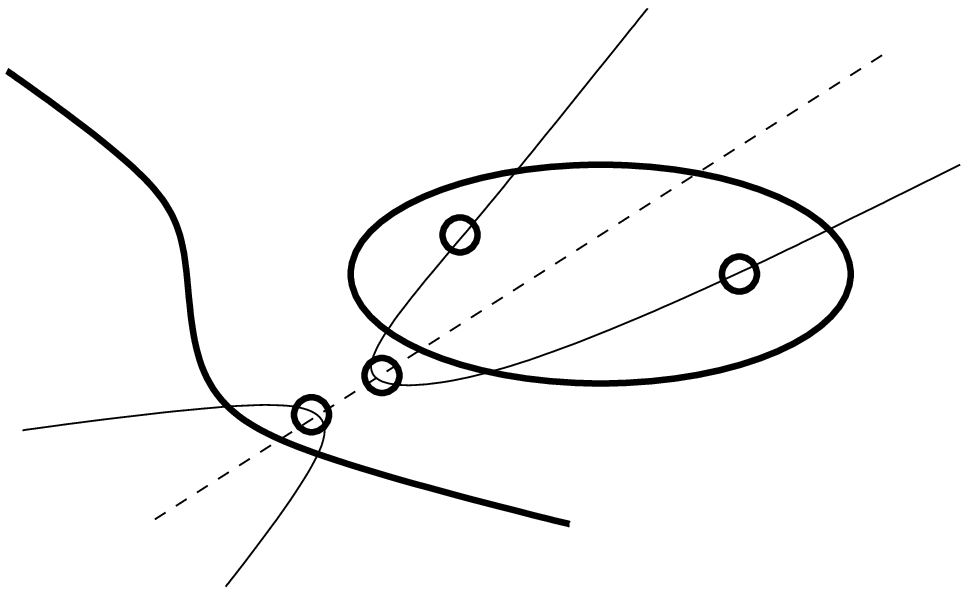}
\\ a)&b)
\end{tabular}
\caption{}
\label{bezout1}
\end{figure}
Given a curve of odd degree on $\mathbb RP^2$, one can speak about convexity~:
the segment defined by two points $a$ and $b$ is the
connected component of the line $(a,b)\setminus\{a,b\}$ which has an even number of intersection points
with the odd component of the curve.
\begin{lemma}\label{convexe}
 Let $C$ be  a real curve of degree $7$ with at least $6$
 ovals and a nest.

Pick a point in each inner oval. Then  these points are the vertices
of a convex polygon in $\mathbb RP^2$. 

Moreover, if a line $L$ passes through two outer ovals $O_1$ and $O_2$ and separates the inner ovals, then $O_1$ and $O_2$ are not
on the same connected component of $L\setminus (Int(O)\cup\mathcal J)$ where $O$ is the non-empty oval.
\end{lemma}
\textit{Proof. }Suppose there exist four empty ovals contradicting the
lemma as depicted in Figure \ref{bezout1}a) and b). Then the conic
passing through these ovals and another  one intersects the curve in
at least $16$ points which contradicts the Bezout theorem.~\findemo

\subsection{Braid theoretical methods}\label{braid}
Here we recall the method exposed by Orevkov in \cite{O1}. Our setting
is a little bit different from \cite{O1} since we
 consider curves which intersect the exceptional
divisor. However, all the
results stated in \cite{O1} remain valid for such curves. The proofs
of \cite{O1} should be adapted using the local model $z=constant$ for
the lines of the pencil, and $z=\pm\frac{1}{w}$
(resp.$z=\pm\frac{1}{w^2}$) for the curve at each
point of the curve which is not (resp. which is) a tangency point with
the pencil and which lie on the exceptional divisor.

The \textit{group of braids with $m$-strings} is defined as
$$B_m=\langle \sigma_1, ..., \sigma_{m-1}| \sigma_i\sigma_{i+1}\sigma_i=\sigma_{i+1}\sigma_{i}\sigma_{i+1},
 \sigma_i\sigma_j=\sigma_j\sigma_i \textrm{ if }|i-j|>1\rangle.$$

We call the \textit{exponent sum} of a braid $b=\prod_{j=1}^n \sigma_{i_j}^{k_j}$ the integer
$e(b)=\sum_{j=1}^n k_j. $
A braid $b$ is called quasipositive if it is a product of conjugated elements of $\sigma_1$ in $B_m$.

From now on, we fix a base $B$ of  $\Sigma_n$,
a fiber $l_{\infty}$ of $\mathcal L$ and an $\mathcal L$-scheme $A$ on
$\Sigma_n$. Suppose moreover that all the intersection points of $A$
and $E$ are nonsingular points of $A$, 
and that $A$ intersects each fiber in $m$ or $m-2$ real points
(counted with multiplicities) and $l_{\infty}$  in $m$ distinct real
points. Choose a standard coordinate system on $\Sigma_n$ such that
$l_\infty$ has equation $\{x_2=0\}$ (see section \ref{ruled}) and 
a trivialization of the $\mathbb CP^1$-bundle over $B\setminus(B\cap
l_\infty)$. 

\begin{figure}[h]
\begin{center}
\begin{tabular}{c}
\includegraphics[height=2cm, angle=0]{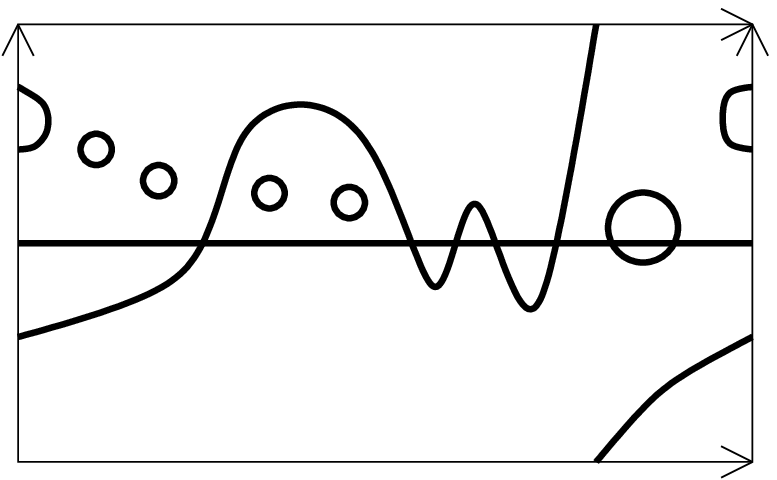}
\end{tabular}
\end{center}
\caption{}
\label{ex braid}
\end{figure}

Now, we describe how to encode such an $\mathcal L$-scheme
$A$. Examine the real part of the fibration from $x=-\infty$ to 
$x=+\infty$. Each time the pencil of real lines is not transversal to
$A$ or meets an intersection point 
of $A$ with $E$,  do the following~: 
\begin{itemize}
\item  if the pencil of lines has a tangency point $p$ with $A$ which is not
  on $E$,
  write $\supset_k$ if $A$ intersects a fiber in $m$ 
real points before $p$, and $\subset_k$ otherwise,
\item if the pencil of lines meets a double point $p$ of $A$, write
   $\subset_k\supset_k $ if the tangents are non-real and $\times_k$
   otherwise,  
\item if the pencil of lines has no tangency point with $A$ but meets
  an intersection point $p$ of $A$ and $E$,
 write $/$ if the branch of $A$
  passing through $p$ 
lies in the region $\{y>0\}$ before $p$, and
$\setminus$ otherwise. 
\item  if the pencil of lines has a tangency point $p$ with $A$ which
  is on  $E$, 
  write $\backslash\supset_{m-1}$ if $A$ intersects a fiber in $m$ 
real points before $p$, and $\subset_{m-1}/$ otherwise,
\end{itemize}
In the first two cases $k$ is the number of real intersection points
of the fiber and $A$ strictly below the point $p$, incremented by 1. 
 Thus, we have now
an  encoding $s_1\ldots s_r$ of the $\mathcal L-scheme$ $A$. In order
to obtain a braid from this encoding, perform the 
 following substitutions~:
\begin{itemize}
\item replace each $\times_j$ which appears between  $\subset_s$ and $\supset_t$
 by $\sigma_j^{-1}$,
\item replace each $\setminus$  which appears between  $\subset_s$ and $\supset_t$ by $\sigma_1\sigma_2\ldots\sigma_{m-1}$,
\item replace each $/$  which appears between  $\subset_s$ and $\supset_t$ by $\sigma_{m-1}\sigma_{m-2}\ldots\sigma_{1}$,
\item replace each subword $\supset_s\times_{i_{1,1}}\ldots\times_{i_{r_1,1}} ?_1 \times_{i_{1,2}}\ldots\times_{i_{r_2,2}} ?_2
\ldots\times_{i_{1,p}}\ldots\times_{i_{r_p,p}} \subset_t$
with $?_i=\setminus $ or $/$
by $\sigma_s^{-1}\tau_{s,m-1}\sigma_{i_{1,1}}^{-1}\ldots \sigma_{i_{r_1,1}}^{-1} v_1 \sigma_{i_{1,2}}^{-1}\ldots \sigma_{i_{r_2,2}}^{-1} v_2\ldots  v_p \sigma_{i_{1,p}}^{-1}\ldots \sigma_{i_{r_p,p}} ^{-1}\tau_{m-1,t}$,
\\where
\begin{center}
\begin{tabular}{l}

$v_l=\left\{\begin{array}{ll}
\sigma_1\sigma_2\ldots\sigma_{m-3}\sigma_{m-2}^2 & \textrm{ if
}?_l=\setminus,\\
\sigma_{m-2}^{-1}\sigma_{m-1}^2\sigma_{m-2}\ldots\sigma_1 & \textrm{
if }?_l=/,\\ 
\end{array}\right.$
\\ \\

$\tau_{s,t}=\left\{\begin{array}{ll}
(\sigma_{s+1}^{-1}\sigma_s)(\sigma_{s+2}^{-1}\sigma_{s+1})
\ldots(\sigma_t^{-1}\sigma_{t-1}) & \textrm{if }t>s,\\
(\sigma_{s-1}^{-1}\sigma_s)(\sigma_{s-2}^{-1}\sigma_{s-1})
\ldots(\sigma_t^{-1}\sigma_{t+1}) & \textrm{if }t<s,\\
1 & \textrm{if }t=s.
\end{array}\right.$
\end{tabular}

\end{center}
\end{itemize}
Then we obtain a braid $b_\mathbb R$. We define the braid associated
 to the $\mathcal L-scheme$ $A$, denoted $b_A$, 
 as the braid
$b_\mathbb R\Delta_m^n$, where $\Delta_m$ is the Garside element of
 $B_m$ and which is given by 
$$\Delta_m=(\sigma_1\ldots\sigma_{m-1})
 (\sigma_1\ldots\sigma_{m-2})\ldots(\sigma_1\sigma_2)\sigma_1.$$  
\textbf{Example. } The encoding and the braid corresponding to the real
$\mathcal L-scheme$ on $\Sigma_2$ depicted in Figure \ref{ex braid}
are (we have abbreviated the pattern $\subset_k\supset_k $ by $o_k$)~:
\begin{center}
\begin{tabular}{lll}
$\supset_3 o_3^2\times_1 o_2^2\times_1^4 / \subset_3\times_2^2\supset_3\subset_3$& and
$\sigma_3^{-3}\sigma_1^{-1}\sigma_2^{-1}\sigma_3\sigma_2^{-2}\sigma_3^{-1}\sigma_2\sigma_1^{-4}
\sigma_2^{-1}\sigma_3^{2}\sigma_2\sigma_1\sigma_2^{-2}\sigma_3^{-1}\Delta_4^2$
\end{tabular}
\end{center}

\vspace{2ex}
Orevkov has proved the following results, where $r$ is the number of (real) intersection points of
$A$ with $E$.
\begin{thm}[Orevkov, \cite{O1}]
The  $\mathcal L-scheme$ $A$ is realizable by a pseudoholomorphic curve of bidegree $(m,r)$ if and only if the braid $b_A$ is
quasipositive.
\end{thm}
\begin{prop}[Orevkov, \cite{O1}]\label{operations}
Let $A$ be an $\mathcal L-scheme$, and $A'$ be the $\mathcal L-scheme$ obtained
from $A$ by one of the following elementary operations~:
$$\begin{array}{c c c c c c cc}
\times_j\supset_{j\pm 1}\to\times_{j\pm 1}\supset_j &, & &
\subset_{j\pm 1}\times_j\to\subset_j\times_{j\pm 1} &, & &
\times_j u_k\to u_k\times_j &,
\end{array}$$
$$\begin{array}{c c c c c c cc c c}
\backslash\supset_{m-1}\leftrightarrow /\supset_1 &, & &
\subset_{m- 1}/\leftrightarrow \subset_1\backslash &, & &
\backslash u_k\leftrightarrow u_k\backslash &, & & / u_k\leftrightarrow u_k/ 
\end{array}$$
$$\begin{array}{c c c c c c}
o_k\subset_k\supset_{k- 1}&\longleftrightarrow&\subset_k\times_{k-1}\supset_{k}&\longleftrightarrow&
\subset_{k-1}\supset_{k}o_k&,
\end{array}$$
$$\begin{array}{c c c c c c cc}
\subset_j\supset_{j\pm 1}\to\emptyset &, & &
\subset_j\supset_{k}\to\supset_k\subset_j &, & &
o_j\to \emptyset &,
\end{array}$$
where $|k-j|>1$ and $u$ stands for  $\times$, $\subset$, or $\supset$.

Then if $A$ is realizable by a nonsingular real pseudoholomorphic
curve of bidegree $(m,r)$ in $\Sigma_n$, so is  $A'$.

Moreover, if $A$ is realizable by a dividing pseudoholomorphic curve of bidegree $(m,r)$ and if $A'$ is obtained from $A$ by one
of the previous elementary operations but the last one, so is $A'$.
\end{prop}

In this paper, we use the following tests to show that a given braid is not quasipositive.
Propositions \ref{alex} and \ref{double alex} are a
corollary of the Murasugi-Tristram inequality (see \cite{O1}) and Proposition \ref{square} is a corollary of the generalized Fox-Milnor theorem (see \cite{O2}).
\begin{prop}[Orevkov, \cite{O1}]\label{alex}
If a braid $b$ in $B_m$ is quasipositive and $e(b)<m-1$, then the Alexander polynomial of $b$ is identically zero.
\end{prop}
\begin{prop}[Orevkov, \cite{O1}]\label{double alex}
If a braid $b$ in $B_m$ is quasipositive and $e(b)=m-1$, then all the roots of the Alexander polynomial of $b$ situated on the unit circle are of order at least two.
\end{prop}
\begin{prop}[Orevkov, \cite{O2}]\label{square}
If a braid $b$ in $B_m$ is quasipositive and $e(b)=m-1$, then $det(b)$ is a square in $\mathbb Z$.
\end{prop}
One can refer to \cite{Lic} or \cite{O1} for the definitions of the Alexander polynomial and the determinant of a braid.

\subsection{General facts about symmetric curves on the real
  plane}\label{general}
In this article, we will use two results by Fiedler and one result by
Trille. As these results have never been published before, we give here an
outline of their proof.

 We denote by $B_0$ the line $\{y=0\}$. 
\begin{prop}[Fiedler, \cite{Fie}]\label{inter max}
If $C$ is a dividing symmetric curve of degree $d$ on $\mathbb RP^2$, then
\begin{displaymath}
Card(\mathbb RC\cap B_0)=d\textrm{ or }Card(\mathbb RC\cap B_0)=0
\end{displaymath}
\end{prop}
\textit{Proof. }The involution $s$ acts
locally on $C$ at all its fixed points like a rotation.
Suppose that $p$ is a point of
$\mathbb RC\cap B_0$. Then, $p$ is a fixed point of $s$ and $s$
exchanges the two halves of $C\setminus \mathbb RC$ locally at
$p$. Hence, $s$
exchanges globally the two halves of $C\setminus \mathbb RC$. So $s$
cannot have any non real fixed point and the proposition is proved.~\findemo

Thus, in the case of a dividing symmetric curve of odd degree, all the common points of the curve and $B_0$ are real.

\vspace{1ex}

The involutions $s$ and $conj$ commute, so $s\circ conj$ is an anti-holomorphic involution on $\mathbb CP^2$. The real part of this real structure
 is a real projective plane $\widetilde{\mathbb RP^2}=\{[x_0:ix_1:x_2]\in  \mathbb CP^2|
(x_0,x_1,x_2)\in \mathbb R^3\setminus\{0\}\}$, and it is clear that
$\widetilde{\mathbb RP^2}\cap \mathbb RP^2=\mathbb RB_0\cup\{[0:1:0]\}$.
A symmetric curve $C$ is real for the structures defined by $conj$ and
 $s\circ conj$. Denote $\widetilde{\mathbb RC}$ its real part for
 $s\circ conj$. This is called the \textit{mirror curve} of $\mathbb RC$.

For a maximal symmetric curve, the real scheme realized by the mirror curve is uniquely determined.
\begin{thm}[Fiedler, \cite{Fie}]\label{max sym}
The mirror curve of a maximal symmetric curve of degree $2k+1$ is a nest of depth $k$ with a pseudo-line.
\end{thm}
\textit{Proof. }Let $C$ be a maximal symmetric curve of degree $2k+1$
and put $C\setminus \mathbb RC=C^+\cup C^-$. According to Proposition
\ref{inter max}, all the fixed points of $s$ are real, and $s\circ conj$
defines an involution of $C^+$. Glue some disks to $C^+$ along
$\mathbb RC$ in order to obtain a sphere $S$ and extend $s\circ conj$
to the whole $S$. Then, one sees that the map $s\circ conj$ is isotopic to a
reflexion on $S$. Now, the lift of the fixed point set of $s\circ
conj$ on $C^+$ to $C$ shows that $\widetilde{\mathbb RC}$ contains at least $k$ ovals and a
pseudo-line and that these $k+1$ components divide $C$. So, $\widetilde{\mathbb RC}$  cannot
have other components, and according to the Rokhlin-Mischachev
orientation formula, it  is a nest of depth $k$ with a pseudo-line.~\findemo

Denote by $\mathcal L_p$ the pencil of lines through the point $[0:1:0]$ in $\mathbb CP^2$.
If $C$ is a real symmetric curve of degree $2k+1$ on $\mathbb RP^2$,
the curve $X=C/s$ is a real curve of bidegree $(k,1)$ on $\Sigma_2$ and is called the \textit{quotient curve} of $C$. The $\mathcal L$-scheme realized by $\mathbb RX$
is obtained by gluing the $\mathcal L_p$-schemes realized by $\mathbb
RC/s$ and $\widetilde{\mathbb RC}/s$ along $B_0$. Since there is no
ambiguity, we will denote by $B_0$ the (symmetric) line $\{y=0\}$ in $\mathbb
CP^2$ as well as its quotient curve in $\Sigma_2$.

Conversely, a symmetric curve of degree $2k+1$ is naturally associated  to any arrangement of a curve $X$ of bidegree $(k,1)$ and a base on $\Sigma_2$.

\begin{prop}\label{cplx orient}
Let $C$ be a nonsingular real pseudoholomorphic symmetric curve on
$\mathbb RP^2$. Then $C$ is smoothly and equivariantly isotopic to a
nonsingular real pseudoholomorphic symmetric curve $C'$ on $\mathbb
RP^2$ such that all the tangency points of the invariant components
of $\widetilde{\mathbb RC'}$ with the pencil of lines $\mathcal L_p$
lie on $B_0$.
\end{prop}
\textit{Proof. }Suppose that $\widetilde{\mathbb RC'}$ has tangency points with a line of $\mathcal L_p$ not on $B_0$. Then push all the corresponding tangent points of the quotient curve above $B_0$
 applying the first elementary operation of Proposition
 \ref{operations} to the curve $\widetilde{\mathbb RC'}\cup B_0$.
The resulting symmetric curve satisfies the conditions of the proposition.~\findemo

\begin{figure}[h]
      \centering
 \begin{tabular}{cccc}
\includegraphics[height=2cm, angle=0]{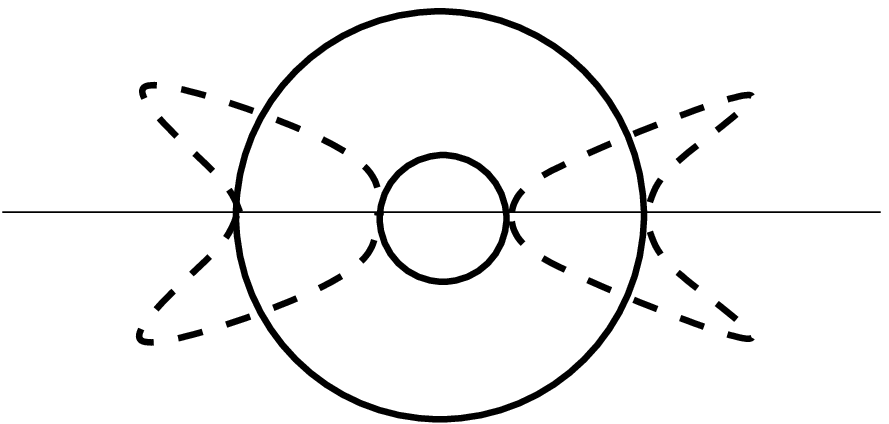}    &
\includegraphics[height=2cm, angle=0]{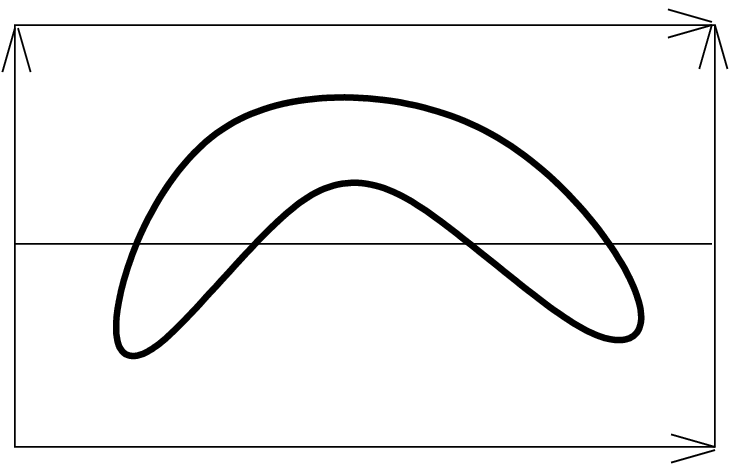}    &
\includegraphics[height=2cm, angle=0]{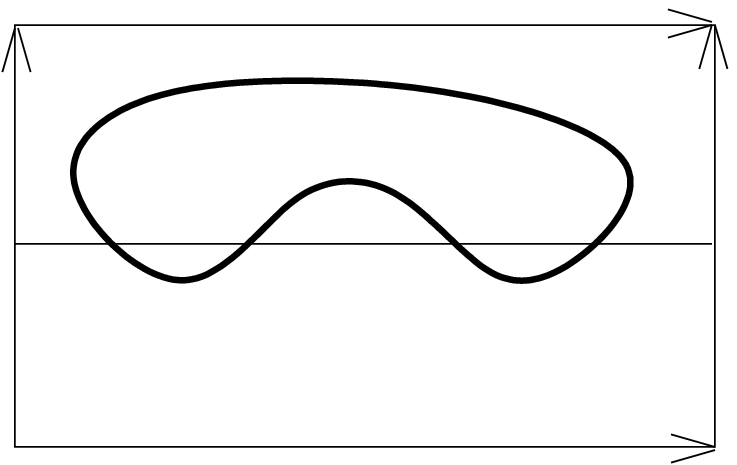}    &
 \includegraphics[height=2cm, angle=0]{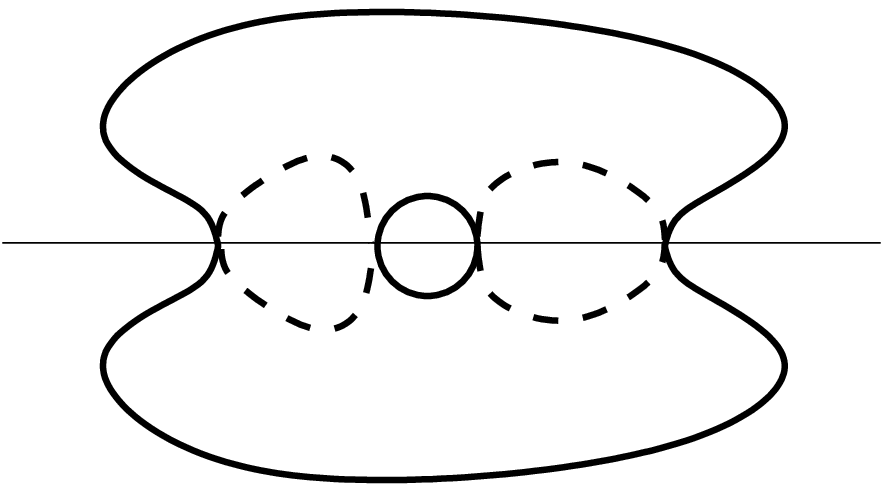}
\\ a)&b) &c)&d)
 \end{tabular}
\caption{}
 \label{sym isotopy}
\end{figure}

\textbf{Example. }The symmetric curves of degree $4$ depicted in Figures
\ref{sym isotopy}a) and d) are equivariantly isotopic in $\mathbb CP^2$.
The dashed curve represent their mirror curves. The corresponding quotient curves
 are depicted in Figures \ref{sym isotopy}b) and c).

The curves $\mathbb RC$ and $\mathbb RC'$ have the same complex
orientations. Using  the invariant components of $\widetilde{\mathbb
  RC'}$, one can apply the Fiedler orientations alternating rule  (see \cite{V3}) to 
$\mathbb RC'$. The following proposition is an application of this observation.

\begin{prop}\label{oval quotient}
If $C$ is a dividing symmetric curve of odd degree on $\mathbb RP^2$, then an oval of $\mathbb RC$ and an oval of $\widetilde{\mathbb RC}$
cannot intersect in more than $1$ point.
\end{prop}
\textit{Proof. }
Suppose that a curve $C$ contradicts the lemma and denote by $O$
(resp. $\widetilde O$) an invariant oval of $\mathbb RC$
(resp. $\widetilde{\mathbb RC}$) such that $O$ and $\widetilde O$ have
two points in common. By Proposition
\ref{cplx orient}, 
one can suppose that all the tangency points of the invariant components
of $\widetilde{\mathbb RC}$ with the pencil of lines $\mathcal L_p$
lie on $B_0$. 

\begin{figure}[h]
\begin{center}
\begin{tabular}{c}
\includegraphics[height=3cm, angle=0]{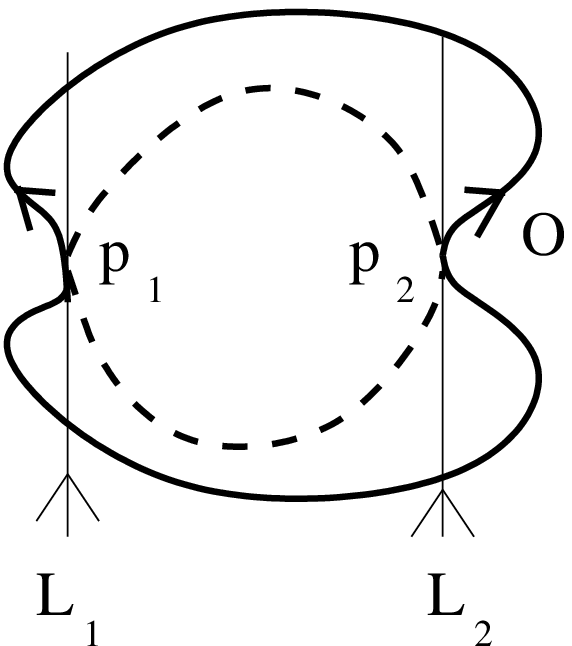}
\end{tabular}
\end{center}
\caption{}
\label{contr fied}
\end{figure}

The homology class realized by $\left( O\cup\widetilde
O\right)/s$  in $H_1(\mathbb R\Sigma_2,\mathbb Z)$ is 
equal to $a\mathbb RB_0$. We have $a=0$ because otherwise,  
 the quotient
curve $X$ would be singular. Indeed, $\mathbb  RX$ is equal to $b\mathbb
RB_0+\mathbb RF$ in this group, $\mathbb
RB_0\circ\mathbb
RB_0=\mathbb RF\circ\mathbb RF=0$ and $\mathbb
RB_0\circ\mathbb RF=1$. Hence,  $\left( O\cup\widetilde
O\right)/s$ is an oval of $\mathbb RX$ and
$O$ and $\widetilde O$ are arranged in $\mathbb CP^2$ as shown in Figure \ref{contr fied}. Then,
transport the orientation of $O$ at $p_1$ to $p_2$ using 
 the Fiedler orientations alternating rule. The two orientations are
 not consistent, so there is a contradiction.~\findemo

In general, there is no link between the type of the symmetric curve and the type of its quotient curve. However, if
both the symmetric curve and its mirror curve are of type I, there is no ambiguity.
\begin{prop}[Trille, \cite{Tr2}]\label{quotient dividing}
If a symmetric curve on $\mathbb RP^2$ and its mirror curve are of type $I$, so is the quotient curve.
\end{prop}
\textit{Proof. }The set $C\setminus \left ( \mathbb RC\cup
\widetilde{\mathbb RC}\right )$ has four connected components, so its
quotient by $s$ has two complex conjugated connected components.~\findemo

We finish this section by a simple observation. Denote by $\pi$ the
projection $\Sigma_2\to B_0$ and by $\pi_p$ the projection $\mathbb
CP^2\setminus\{p\}\to B_0$. Then the set $\pi_p^{-1}(\mathbb RB_0)\cap
C$
can be deduced from the set $\pi^{-1}(\mathbb RB_0)\cap X$. Hence  some
information on the complex orientations of $C$ (if any) can be deduced from
$\pi^{-1}(\mathbb RB_0)\cap X$. This observation will be very useful
in this paper since if $C$ is of degree $7$, the curve $X$ is of bidegree $(3,1)$. In particular, the set $\pi^{-1}(\mathbb RB_0)\cap X$ can
be extracted only from the $\mathcal L$-scheme realized by $X$.

\section{Real trigonal graphs on $\mathbb CP^1$ and real trigonal
  algebraic curves}\label{comb theo}
In \cite{O3}, Orevkov reformulates the existence of trigonal real algebraic
curves realizing a given $\mathcal L$-scheme on
$\Sigma_n$ in terms of some colored graphs on $\mathbb CP^1$. Guided
by this article, we give in this section an efficient algorithm to
check whether an $\mathcal L$-scheme is realizable by  a trigonal real
algebraic curve  on $\Sigma_n$. 

Using these so-called \textit{real trigonal graphs on $\mathbb CP^1$}, Orevkov
(\cite{O5}) obtained a classification of trigonal real algebraic curves on
$\Sigma_n$ up to isotopy which respects the pencil of lines, in terms
of gluing of cubics. 

Note that real trigonal graphs on $\mathbb CP^1$ are a particular
case of \textit{real rational graphs} defined in \cite{Br2} and used
to deal with \textit{root schemes}.

\subsection{Root scheme associated to a trigonal $\mathcal
  L$-scheme}\label{ro sc}
\begin{defi}
A root scheme is a k-uplet
$((l_1,m_1),\ldots,(l_k,m_k))\in(\{p,q,r\}\times\mathbb N)^k$ with $k$
a natural number (here, $p,q$ and $r$ are symbols and do not stand for
natural numbers).

\noindent A root scheme $((l_1,m_1),\ldots,(l_k,m_k))$ is realizable
by polynomials of degree $n$ if there exist two real polynomials in
one variable of degree $n$,
with no common roots, $P(x)$ and $Q(x)$ such that if
\mbox{$x_1< x_2< \ldots < x_k$}
are the real roots of $P,Q$ and $P+Q$, then
$l_i=p$ (resp., $q, r$) if $x_i$ is a root of $P$ (resp., $Q, P+Q$) and
$m_i$ is the multiplicity of $x_i$.

The polynomials $P$,  $Q$ and $P+Q$ are said to realize the root scheme $((l_1,m_1),\ldots,(l_k,m_k))$.
\end{defi}

In what follows, $n$ is a positive integer and
$C(x,y)=y^3+b_2(x)y+b_3(x)$ is a real polynomial, where $b_j(x)$ is a
real polynomial of degree $jn$ in $x$. By a suitable change of
coordinates in $\Sigma_n$, any real algebraic trigonal curve on $\Sigma_n$
can be put into this form. 

Denote by $D=-4b_2^3-27b_3^2$ the discriminant of $C$ with
respect to the variable $y$. The knowledge of the
root scheme realized by $-4b_2^3$, $27b_3^2$ and $D$  allows one to
recover the $\mathcal L$-scheme realized by  
$C$, up to the transformation $y\mapsto -y$. Indeed, the position of
$C$ with respect to the pencil of lines is given by the sign of the
double root of $C(x_0,Y)$ at each root $x_0$ of $D$, which is the sign of
$b_3(x_0)$.

Consider
a
trigonal $\mathcal L$-scheme
$A$ on $\Sigma_n$
such that
  $A$  intersects some fiber  $l_{\infty}$  in $3$
distinct real points. Consider also the encoding $s_1
\ldots s_r$ of $A$ defined in section \ref{braid}, using the symbols
$\subset$, $\supset $ and $\times$. In this encoding, replace all
the occurrences $\times_k$  by $\supset_k$ $\subset_k$.  This encoding is denoted by $r_1\ldots
r_q$. Define a root scheme $RS_A=(S_1,\ldots,S_q)$ where  the $S_i$'s are  sequences of pairs $(l_1,m_1),\ldots,(l_{k_i},m_{k_i})$, as
follows~:
\begin{itemize}
\item $S_1=\left\{\begin{array}{ll}
(r,1)& \textrm{ if }n\textrm{ is even and }r_1=\supset_k
\textrm{ and }r_q=\subset_k, \\ & \textrm{ or }n\textrm{ is odd and
}r_1=\supset_k\textrm{ and }r_q=\subset_{k\pm 1},\\ 

(q,2),(r,1)
& \textrm{ otherwise },\end{array}\right.$
\item for $i>1$,
\\$S_i=\left\{\begin{array}{ll}
(r,1)& \textrm{ if
}r_i=\subset_k\textrm{ and }r_{i-1}=\supset_{k},\\
&
 \textrm{or
}r_i=\supset_k\textrm{ and }r_{i-1}=\subset_{k},\\

(p,3),(q,2),(p,3),(r,1)&
\textrm{ if }r_i=\supset_k\textrm{ and }r_{i-1}=\subset_{k\pm 1},\\

(q,2),(r,1) & \textrm{ if
}r_i=\subset_k\textrm{ and
}r_{i-1}=\supset_{k\pm1}.\end{array}\right.$
\end{itemize}

\begin{defi}
The root scheme  $RS_A$ is called the root scheme associated
to $A$.
\end{defi}

This root scheme encode the mutual cyclic order in $\mathbb RP^1$ realized by the roots of the polynomials $b_2^3$ (encoded by the letter $p$), $b_3^2$ (encoded by the letter $q$)
and $D$ (encoded by the letter $r$).
 
\vspace{1ex}
\textbf{Example. }The encoding of the $\mathcal
L$-scheme depicted in Figure \ref{ex comb}a) is $\supset_2\subset_1\supset_2 o_2\subset_2$ and
the associated root scheme is 
$$\Big( (q,2),(r,1),(p,3), (q,2), (p,3), (r,1),(q,2),(r,1),(r,1),(r,1),(r,1) \Big).$$
\vspace{2ex}

The realizability of a root scheme associated to a trigonal $\mathcal
L$-scheme on 
$\Sigma_n$ can be studied via what we call real trigonal graphs.
 \begin{figure}[h]
      \centering
 \begin{tabular}{cccc}
\includegraphics[height=0.7cm, angle=0]{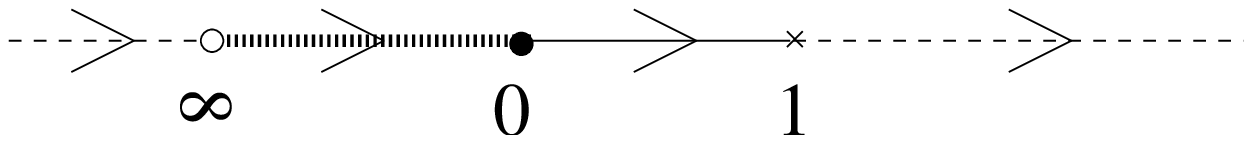}&
\includegraphics[height=1cm, angle=0]{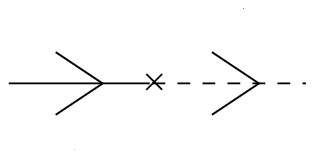}    &
\includegraphics[height=1cm, angle=0]{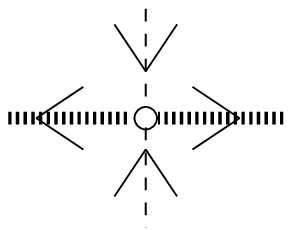}    &
\includegraphics[height=1.3cm, angle=0]{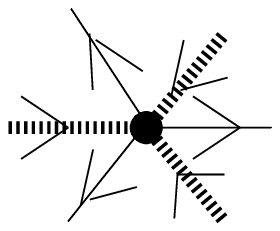}
\\ a)&b) &c)&d)
\end{tabular}
\caption{}
 \label{pregraph}
\end{figure}
\begin{defi}
Let $\Gamma$ be a graph
on $\mathbb CP^1$ invariant under the action
of the complex conjugation and $\pi:\Gamma\to\mathbb RP^1$ a
continuous map.
 Then the coloring and
orientation of $\mathbb RP^1$ shown in Figure~\ref{pregraph}a) defines
a coloring and an orientation of $\Gamma$ via $\pi$.

\noindent The graph $\Gamma$ equipped with this coloring and this
orientation is called a real trigonal graph of degree $n$ if 
\begin{itemize}
\item any vertex of $\Gamma$ has an even valence,
\item for any connected
component $D$ of $\mathbb CP^1\setminus\Gamma$, then
$\pi_{|\partial D}$ is a covering of $\mathbb RP^1$,
\item $\Gamma$ has exactly $6n$ vertices of the kind depicted in Figure
\ref{pregraph}b), $3n$ vertices  of the kind depicted on \ref{pregraph}c) and $2n$
vertices of the kind depicted on \ref{pregraph}d), and no other non-real multiple
points,
\item The set $\pi^{-1}([\infty;0])$ is connected.
\end{itemize}
\end{defi}
Since all the graphs on $\mathbb CP^1$ considered in this article are
invariant under the complex conjugation, we  draw only one
half of them.

\begin{figure}[h]
      \centering
 \begin{tabular}{cccc}
\includegraphics[height=2.5cm, angle=0]{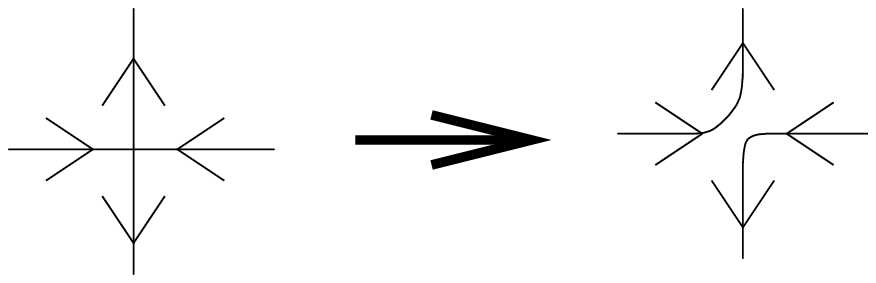}&
\includegraphics[height=1.5cm, angle=0]{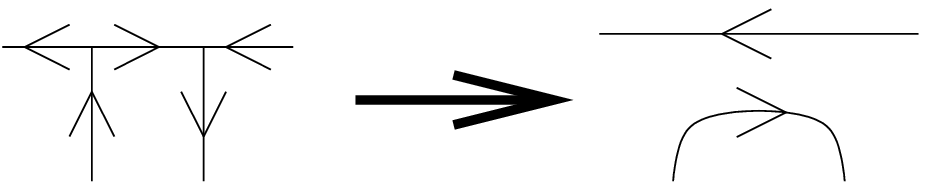}&
\includegraphics[height=1.5cm, angle=0]{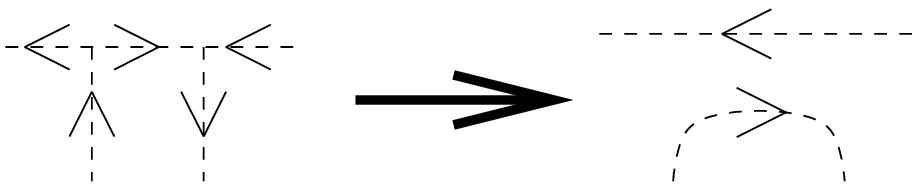}&
\includegraphics[height=1.5cm, angle=0]{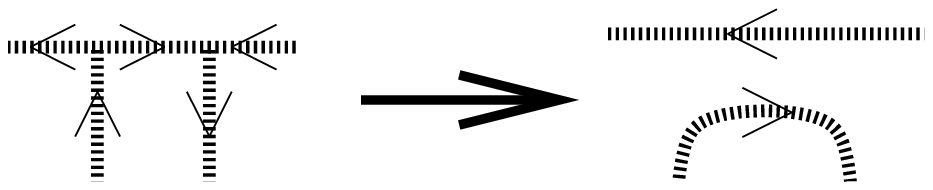}
\\ a)&b) &c)&d)
\end{tabular}
\caption{}
 \label{oper}
\end{figure}

Now, suppose that the trigonal curve $C$ on $\Sigma_n$
realizes the $\mathcal L$-scheme $A$ and that  $-4b_2^3$, $27b_3^2$ and
$D$ realize the root scheme $RS_A$. Color and orient $\mathbb
RP^1$ as depicted in Figure \ref{pregraph}a).
Consider the rational function
$f=\frac{-4b_2^3}{27b_3^2}$ defined on $\mathbb CP^1$ and let
$\Gamma$ be $f^{-1}(\mathbb{R}P^1)\subset\mathbb CP^1$ with the
coloring and the 
orientation induced by  those chosen on $\mathbb RP^1$. 
Let $\widetilde \Gamma$ be the colored
and oriented graph on $\mathbb CP^1$ obtained out of $\Gamma$ by
smoothing equivariantly its non-real double points as depicted in
Figure \ref{oper}a), and by performing operations depicted in Figures
\ref{oper}b), c)and d) in order to minimize the number
of its real double points.
The colored
and oriented graph on $\mathbb RP^1$ obtained as the intersection of
$\widetilde \Gamma$ and $\mathbb RP^1$ can clearly be 
extracted from $RS_A$. 
\begin{defi}
The colored and oriented graph on $\mathbb RP^1$ constructed above is
called the real graph associated
to $A$.
\end{defi}

\begin{figure}[h]
      \centering
 \begin{tabular}{cccc}
\includegraphics[height=2cm, angle=0]{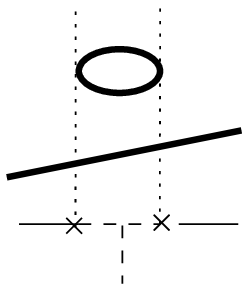}&
\includegraphics[height=2cm, angle=0]{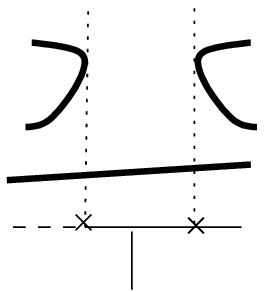}    &
\includegraphics[height=2cm, angle=0]{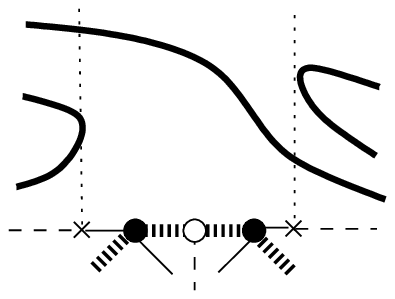}    &
\includegraphics[height=2cm, angle=0]{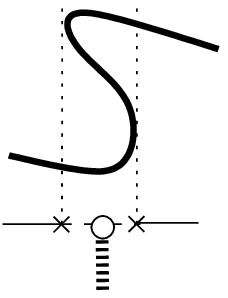}
\\ a)&b) &c)&d)
\end{tabular}
\caption{}
 \label{real graph}
\end{figure}

The real graph associated to $A$ is obtained from $A$ as depicted in
Figure \ref{real graph} (we omitted the arrows).

\begin{figure}[h]
      \centering
 \begin{tabular}{cc}
 \includegraphics[height=2cm, angle=0]{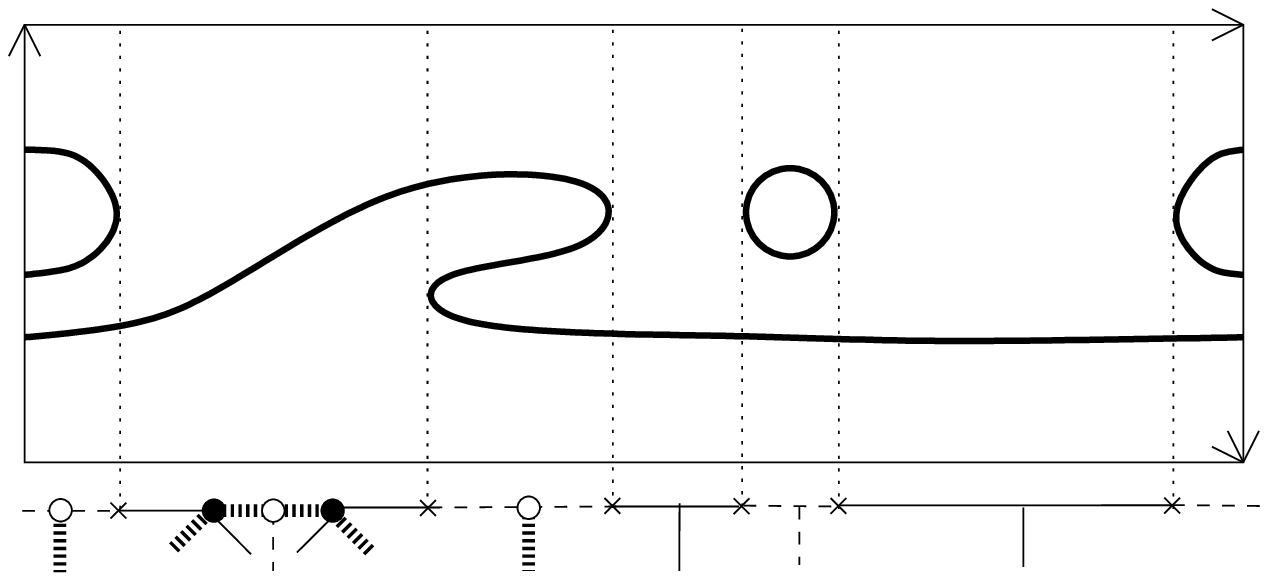} &
\includegraphics[height=2cm, angle=0]{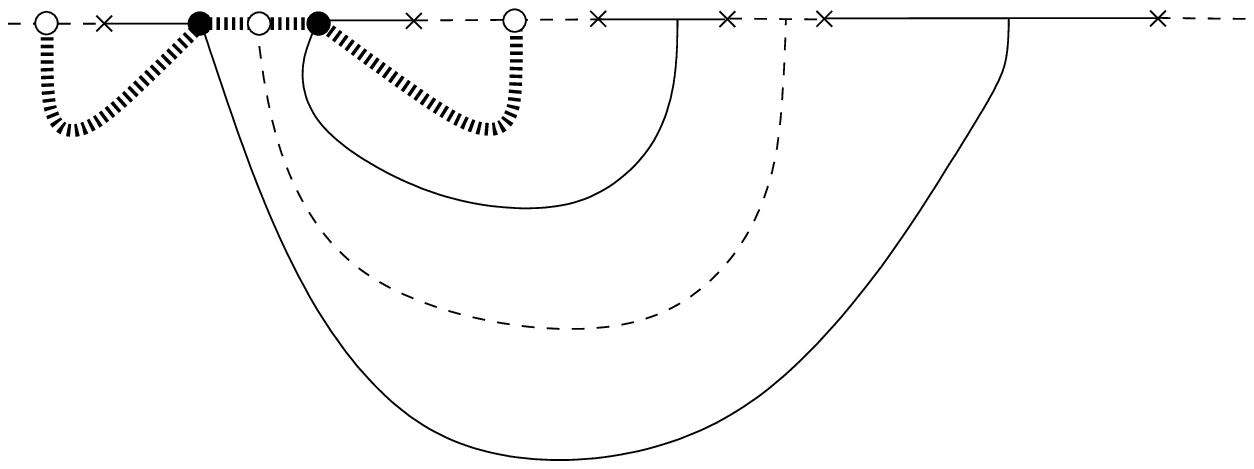} 
\\ a)&b)
 \end{tabular}
\caption{}
 \label{ex realgraph}
\end{figure}

\vspace{2ex}
\textbf{Example. }In Figure \ref{ex realgraph}a) we have depicted a trigonal $\mathcal
L$-scheme on $\Sigma_1$ and its real graph.

\vspace{2ex}

The importance of the real graph associated to an $\mathcal L$-scheme
is given by the following theorem.

\begin{thm}\label{constr trig}
Let $A$ be a trigonal  $\mathcal L$-scheme on $\Sigma_n$ and $G$ its
real graph. Then $A$ is realizable by  a nonsingular trigonal real algebraic curve  on $\Sigma_n$ if and only if there exists a trigonal graph $\Gamma$ of degree $n$ such that $\Gamma\cap \mathbb RP^1=G$.
\end{thm}
\textit{Proof. }If there exists such a trigonal graph of degree $n$,
the existence of the desired real algebraic curve is proved in
\cite{O3} (see also \cite{Shap1}).

Suppose there exists a real algebraic curve $C$ realizing $A$. Define
 $f$ and $\widetilde \Gamma$ as above.We will perform some
operations on one of the halves of $\mathbb CP^1\setminus\mathbb
RP^1$. The final picture will be obtained by gluing the obtained graph
with its image under the complex conjugation.

 \begin{figure}[h]
      \centering
 \begin{tabular}{ccccc}
\includegraphics[height=2.5cm, angle=0]{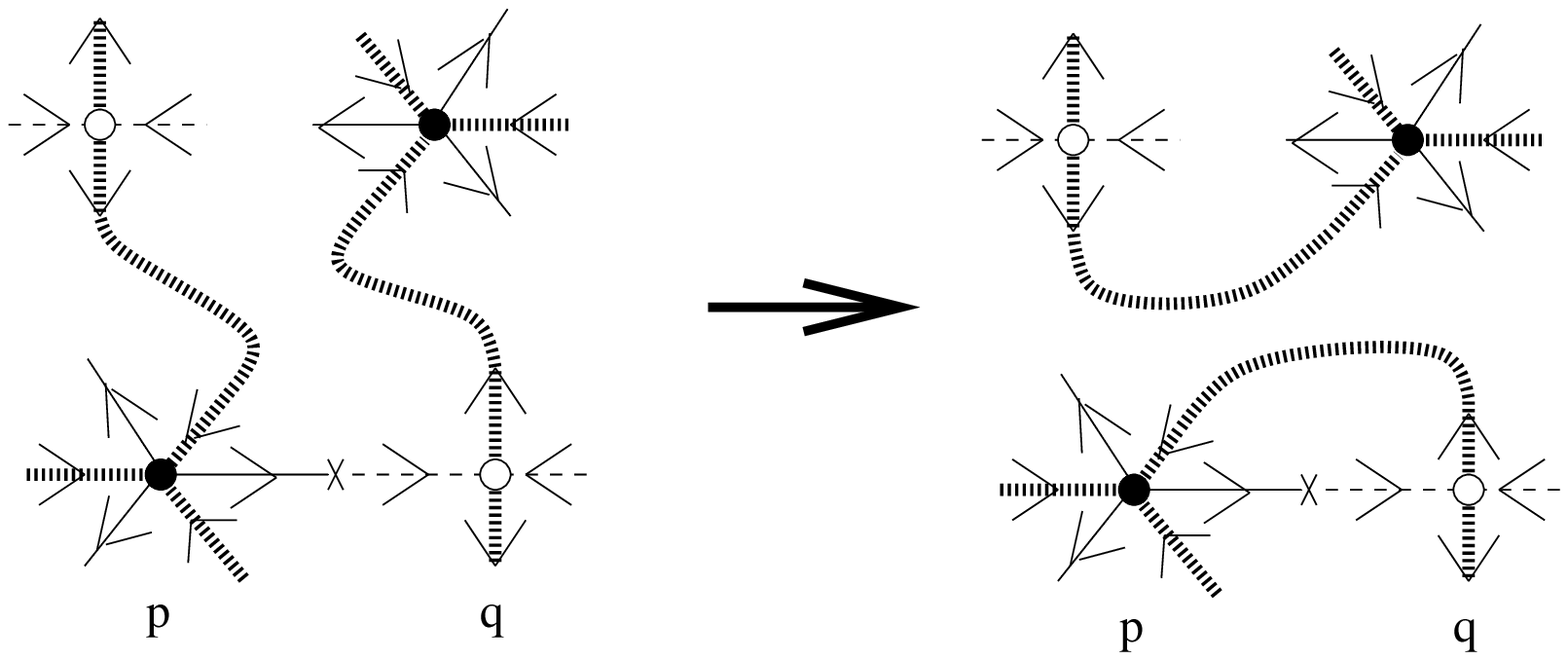}    &&&&
\includegraphics[height=2.5cm, angle=0]{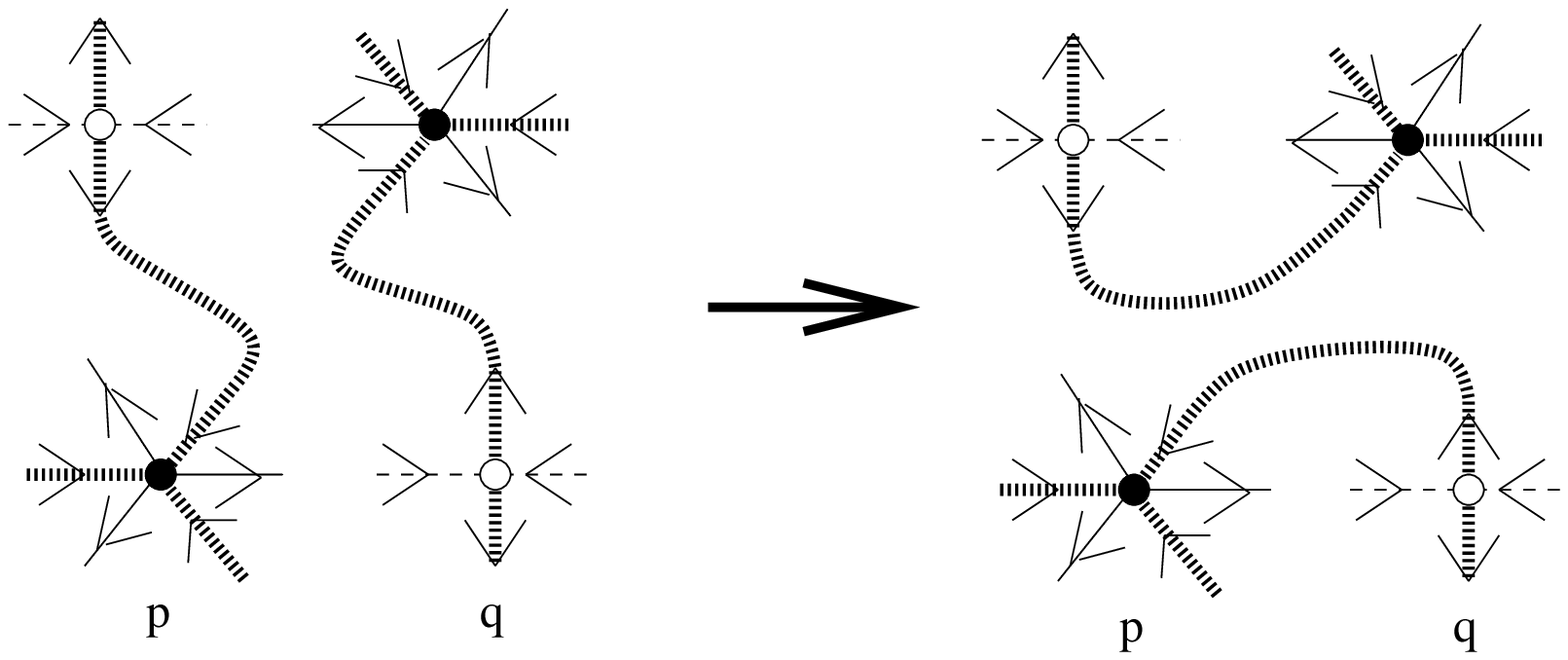}\\
 a)&&&&b) 
\end{tabular}
\caption{}
 \label{oper2}
\end{figure}

If $f^{-1}([\infty;0])$ is not connected, choose $p$ in
$f^{-1}(0)$ and $q$ in $f^{-1}(\infty)$  belonging to different
connected components of $f^{-1}([\infty;0])$. If $p$ and $q$ belong
to the same connected components of $\widetilde \Gamma$, choose $p$ and $q$ such
that they are connected in $\widetilde \Gamma$ by an arc of
$f^{-1}(]0;\infty[)$ and perform on $\widetilde \Gamma$ the operation depicted in Figure
\ref{oper2}a). Otherwise, choose $p$ and $q$ lying on
the boundary of one connected components of $\mathbb
CP^1\setminus\widetilde \Gamma$ and perform on $\widetilde \Gamma$ the operation depicted in Figure
\ref{oper2}b). Perform these operations until
$f^{-1}([\infty;0])$ is connected.\findemo

\textbf{Remark. }The connectedness of $\pi^{-1}([\infty;0])$ is not
necessary to the existence of the algebraic curve. We will use this
condition only in the next section.

\vspace{2ex}

\textbf{Example. }In Figure \ref{ex realgraph}b), we exhibit a real
trigonal graph which ensures the existence of a nonsingular real
trigonal curve in $\Sigma_1$ realizing the $\mathcal
L$-scheme depicted in Figure \ref{ex realgraph}a).

\vspace{2ex}

According to Proposition \ref{operations}, if an  $\mathcal L$-scheme
$A$ is realizable by a pseudoholomorphic curve, so is any $\mathcal
L$-scheme obtained from $A$ by some $\subset_j\supset_{j\pm
  1}\to\emptyset $ operations. Orevkov showed (\cite{O5}) that this is
not the case for algebraic curves. Nevertheless some operations are allowed
on real algebraic trigonal curves.

\begin{figure}[h]
      \centering
 \begin{tabular}{c}
\includegraphics[height=2cm, angle=0]{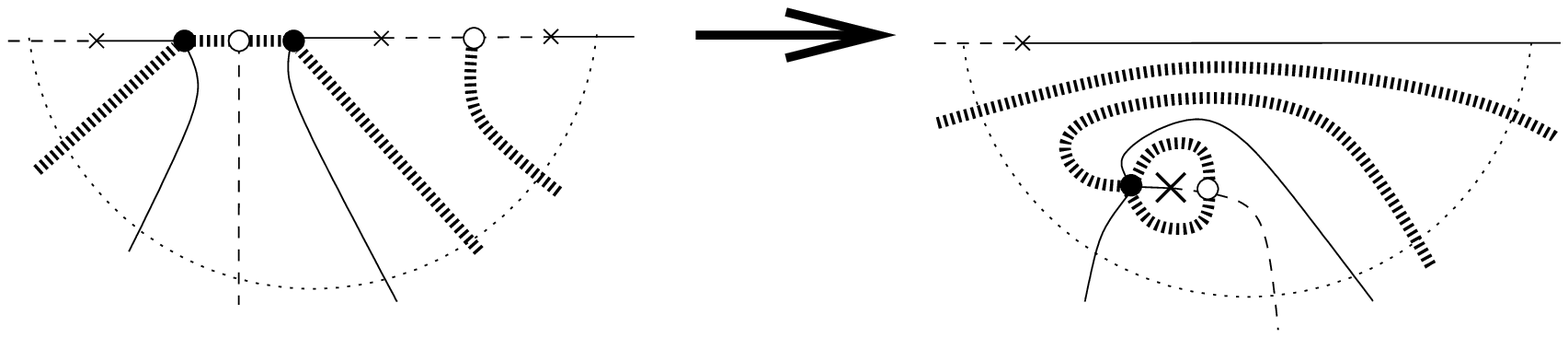}  
\end{tabular}
\caption{}
 \label{oper3}
\end{figure}

\begin{prop}\label{alg operations}
Let $A$ be an $\mathcal L-scheme$, and $A'$ be the $\mathcal L-scheme$ obtained
from $A$ by one of the following elementary operations~:
$$\begin{array}{c c c c c c }
\supset_{j}\subset_{j\pm 1}\supset_{j}\to \supset_{j}&, & &
\subset_{j}\supset_{j\pm 1}\subset_{j}\to \subset_{j}.&
\end{array}$$

Then if $A$ is realizable by a nonsingular real algebraic trigonal
curve in $\Sigma_n$, so is  $A'$.
\end{prop}
\textit{Proof. }Suppose that $A$ is realizable by a nonsingular real algebraic trigonal
curve in $\Sigma_n$. Let $\Gamma$ be a trigonal graph of degree $n$ as in
Theorem \ref{constr trig}. The first operation corresponds to the
operation depicted in Figure \ref{oper3} performed on
$\Gamma$. The
second operation is symmetric to the first one.
 
Proposition \ref{alg operations} follows now from Theorem \ref{constr trig}.
\findemo

\subsection{Comb theoretical method}
In this section, we reformulate  in an
algorithmic way the following problem : given a trigonal $\mathcal
L$-scheme $A$, does there exist a real trigonal graph $\Gamma$ of degree
$n$ such that $\Gamma\cap \mathbb RP^1$ is the real graph of $A$?
 \begin{figure}[h]
      \centering
 \begin{tabular}{ccccccc}
\includegraphics[height=0.7cm, angle=0]{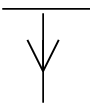}    &
\includegraphics[height=0.7cm, angle=0]{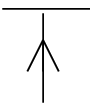}    &
\includegraphics[height=0.7cm, angle=0]{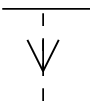}     &
\includegraphics[height=0.7cm, angle=0]{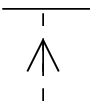}     &
\includegraphics[height=0.7cm, angle=0]{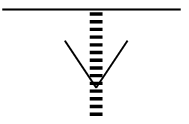}&
\includegraphics[height=0.7cm, angle=0]{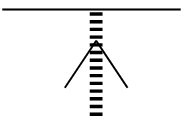} &
\includegraphics[height=0.7cm, angle=0]{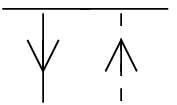}
\\ $g_1$&$g_2$ &$g_3$&$g_4$&$g_5$&$g_6$&$g_1g_4$
\end{tabular}
\caption{}
 \label{gener semi group}
\end{figure}

Denote by $\mathfrak M$ the semigroup generated by the elements $g_1,\ldots
g_6$ in $\mathbb R^2$ depicted in Figure \ref{gener semi group}. The
multiplication of two elements $\mathfrak m_1$ and $\mathfrak m_2$ in $\mathfrak M$ is the
attachment of the right endpoint of $\mathfrak m_1$ to the left endpoint of
$\mathfrak m_2$. 
\begin{defi} 
The elements of $\mathfrak M$ are called combs.
\end{defi}
For example, the comb $g_1g_4$ is depicted in Figure \ref{gener semi
group}. The unit element of $\mathfrak M$ is denoted by $1$.
\begin{defi}
A weighted comb is a quadruplet $(\mathfrak m,\alpha,\beta,\gamma)$ in
$\mathfrak M\times\mathbb Z^3$.
\end{defi}

Consider
a
trigonal $\mathcal L$-scheme
$A$ on $\Sigma_n$
which  intersects some fiber  $l_{\infty}$  in $3$
distinct real points. Consider also the encoding  $r_1\ldots
r_q$ of $A$ defined in
section \ref{ro sc}, using the symbols
$\subset$, $\supset $.  Define the weighted combs $(\mathfrak m_i,\alpha_i,\beta_i,\gamma_i)$ as
follows~:
\begin{itemize}
\item $(\mathfrak m_0,\alpha_0,\beta_0,\gamma_0)=(1,6n,3n,2n)$
\item $(\mathfrak m_1,\alpha_1,\beta_1,\gamma_1)=\left\{\begin{array}{ll}
(g_3,6n-1,3n,2n)& \textrm{ if }n\textrm{ is even and }r_1=\supset_k
\textrm{ and }r_q=\subset_k, \\ & \textrm{ or }n\textrm{ is odd and
}r_1=\supset_k\textrm{ and }r_q=\subset_{k\pm 1},\\ (g_5,6n-1,3n-1,2n)
& \textrm{ otherwise },\end{array}\right.$
\item for $i>1$,
\\$(\mathfrak m_i,\alpha_i,\beta_i,\gamma_i)=\left\{\begin{array}{ll}
(\mathfrak m_{i-1}g_2,\alpha_{i-1}-1,\beta_{i-1},\gamma_{i-1})& \textrm{ if
}r_i=\subset_k\textrm{ and }r_{i-1}=\supset_{k},\\
(\mathfrak m_{i-1}g_3,\alpha_{i-1}-1,\beta_{i-1},\gamma_{i-1})& \textrm{ if
}r_i=\supset_k\textrm{ and }r_{i-1}=\subset_{k},\\
(\mathfrak m_{i-1}g_5,\alpha_{i-1}-1,\beta_{i-1}-1,\gamma_{i-1}) &
\textrm{ if }r_i=\supset_k\textrm{ and }r_{i-1}=\subset_{k\pm 1},\\
(\mathfrak m_{i-1}g_6g_1g_4g_1g_6,\alpha_{i-1}-1,\beta_{i-1}-1,\gamma_{i-1}-2)& \textrm{ if
}r_i=\subset_k\textrm{ and
}r_{i-1}=\supset_{k\pm1}.\end{array}\right.$
\end{itemize}
\begin{defi}
The weighted comb
$(\mathfrak m_q,\frac{\alpha_q}{2},\frac{\beta_q}{2},\frac{\gamma_q}{2})$ is
said to be
associated to the $\mathcal L$-scheme $A$.
\end{defi}

\begin{figure}[h]
      \centering
 \begin{tabular}{ccc}
 \includegraphics[height=2cm, angle=0]{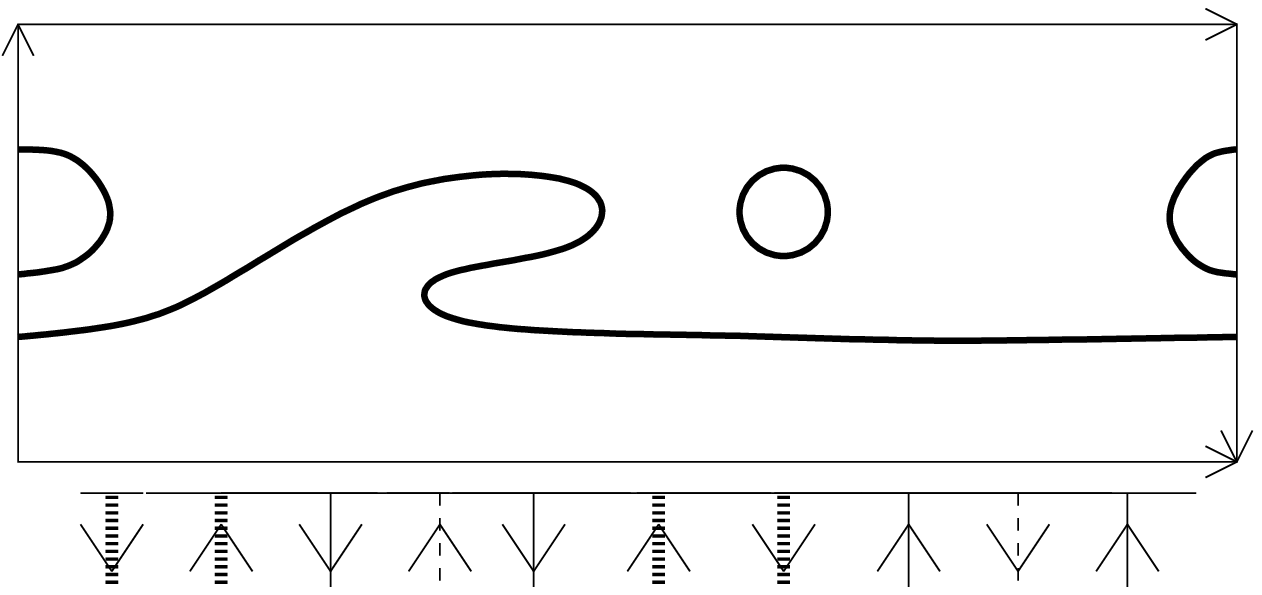} &
\includegraphics[height=2cm, angle=0]{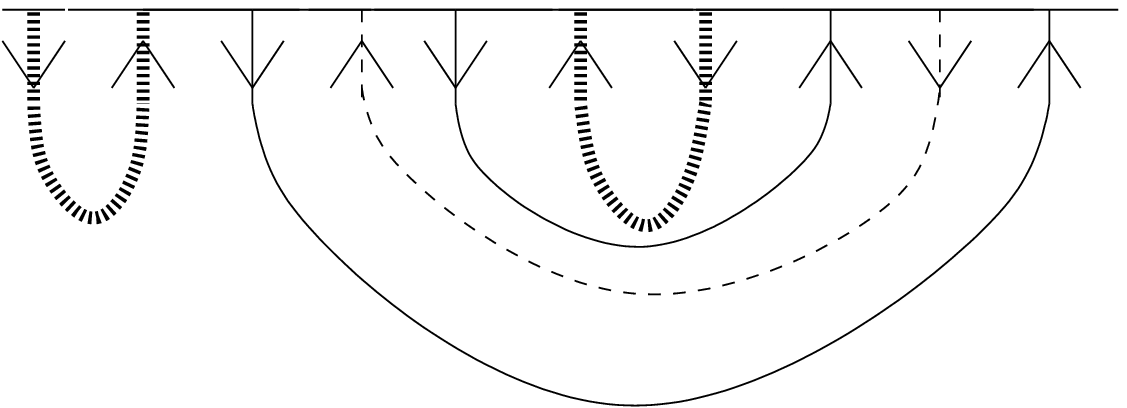} &
\includegraphics[height=2cm, angle=0]{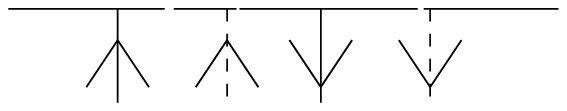}

\\ a)&b)&c)
 \end{tabular}
\caption{}
 \label{ex comb}
\end{figure}

\begin{defi}
Let $\mathfrak m$ be a comb. A closure of $\mathfrak m$ is a subset of $\mathbb R^2$
obtained by joining each generator $g_1$ (resp., $g_3$, $g_5$) in $\mathfrak m$
to a generator $g_2$ (resp., $g_4$, $g_6$) in $\mathfrak m$ by a path on
$\mathbb R^2$
in such a way that these paths do not intersect.

If there exists a closure of $\mathfrak m$, one says that $\mathfrak m$ is closed.
\end{defi}
\textbf{Example. } The weighted comb associated to the $\mathcal L$-scheme
on $\Sigma_1$ depicted in Figure \ref{ex comb}a) is
\\$(g_5g_6g_1g_4g_1g_6g_5g_2g_3g_2,0,0,0)$. A closure of this comb is
shown in Figure \ref{ex comb}b) (compare with Figure \ref{ex realgraph}). The comb depicted in Figure \ref{ex
comb}c) is not closed.

\begin{figure}[h]
      \centering
 \begin{tabular}{cccc}
\includegraphics[height=1.5cm, angle=0]{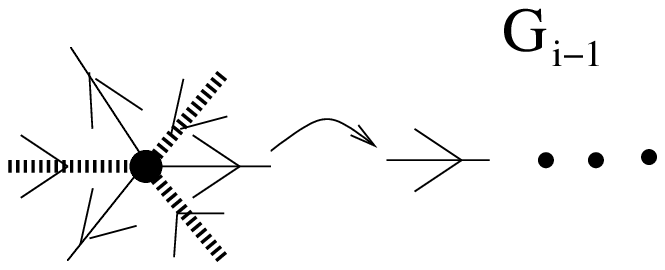}&
\includegraphics[height=1.5cm, angle=0]{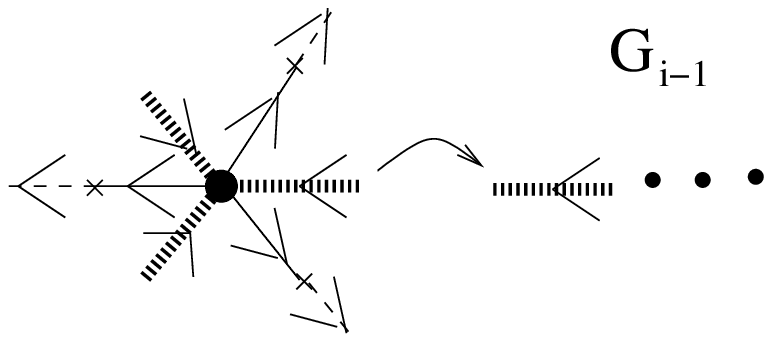}    &
\includegraphics[height=1.5cm, angle=0]{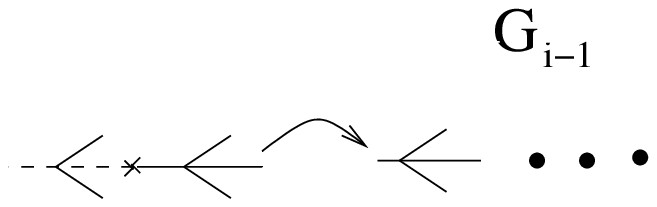}    &
\includegraphics[height=1.5cm, angle=0]{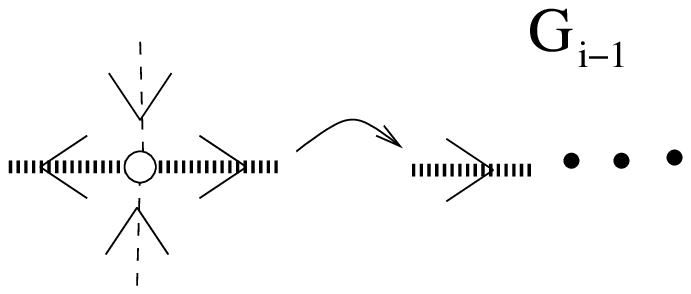}
\\ a)&b) &c)&d)
\end{tabular}
\caption{}
 \label{gluing2}
\end{figure}

\vspace{2ex}
\begin{defi}\label{chain}
A chain of weighted combs is a sequence $(w_i)_{1\le i\le k}$ of
weighted combs, with $w_i=(\mathfrak m_i,\alpha_i,\beta_i,\gamma_i),$ such
that~:
\begin{quote}
$w_k=(\mathfrak m,0,0,0)$,
where $\mathfrak m$ is a closed comb,

 $\forall i\in
\{1\ldots k-1\}$, the weighted comb $w_{i+1}$ is obtained from $w_i$
by one of the following operations :

 \begin{tabular}{llll}
(1) &if $\gamma_i>0$ : &&$g_2\to (g_6g_1)^2g_6$,
$\alpha_{i+1}=\alpha_i$, $\beta_{i+1}=\beta_i$,
$\gamma_{i+1}=\gamma_i-1$

\\ && or & $g_5\to (g_3g_6)^2g_3$,
$\alpha_{i+1}=\alpha_i-3$, $\beta_{i+1}=\beta_i$ ,
$\gamma_{i+1}=\gamma_i-1$ 

\\ (2) &else if $\alpha_i>0$ : & & $g_1\to g_3$,
$\alpha_{i+1}=\alpha_i-1$, $\beta_{i+1}=\beta_i$, $\gamma_{i+1}=0$ \\
(3) & else : && $g_5\to g_4g_5g_4$, $\alpha_{i+1}=0$,
$\beta_{i+1}=\beta_i-1$, $\gamma_{i+1}=0$
\end{tabular}
\end{quote}
where $g_j\to  a$ means ``replace one $g_j$ in $\mathfrak m_i$ by $a$''.
One says that the chain $(w_i)_{1\le i\le k}$ starts at $w_1$.
\end{defi}
\begin{defi}
Let $w$ be a weighted comb. The multiplicity of $w$, denoted by
$\mu(w)$, is defined as the number of chains of weighted combs which
start at $w$.
\end{defi}
\begin{thm}\label{prohib comb}
Let $A$ be a trigonal $\mathcal L$-scheme on $\Sigma_n$,
and $w$ its associated weighted comb. Then $A$ is realizable by
nonsingular real
algebraic trigonal curves  on $\Sigma_n$ if and only if
$\mu(w)>0$ or $w=(1,6n,3n,2n)$.
\end{thm}
\textit{Proof. }Let $w=(\mathfrak m,\alpha,\beta,\gamma)$,
and $G$
be the real graph associated to the $\mathcal L$-scheme $A$. If
$\mathfrak m=1$, it is well known that $A$ is realizable by a real algebraic
trigonal curve  on $\Sigma_n$. Otherwise, a chain of
weighted combs starting at $w$ is a reformulation of the statement~:

"choose
a half $D$ of $\mathbb CP^1\setminus\mathbb
RP^1$;
then there exists a finite sequence $(G_i)_{0\le i\le k}$ of
subsets
of $\mathbb CP^1$ such that:
\begin{itemize}
\item $G_0=G$,
\item for $i$ in $\{1,\ldots,\gamma\}$,
the subset
$G_i$ is obtained from
$G_{i-1}$ by one of the
gluings
 in $D$ depicted in Figures
 \ref{gluing2}a) and b);  denote by $c$ the number of times that  $G_i$ is 
obtained from
$G_{i-1}$ by the
gluing depicted
 in Figure \ref{gluing2}b) for $i$ in $\{1,\ldots,\gamma\}$,
\item for $i$ in $\{\gamma+1,\ldots,\alpha+\gamma-3c\}$,
the subset
$G_i$ is
obtained from $G_{i-1}$ by the gluing in $D$ depicted in Figure
\ref{gluing2}c),
\item for $i$ in $\{\alpha+\gamma-3c+1,\ldots,k-2\}$,
the subset
$G_i$ is obtained
from $G_{i-1}$ by the gluing in $D$ depicted in Figure \ref{gluing2}d),
\item $G_{k-1}$ has no boundary and contains $G_{k-2}$,
\item $k=\alpha+\beta+\gamma-3c+2$,
\item $G_k$ is the gluing of $G_{k-1}$ and its image under the complex
conjugation,
\item $G_k$ is a trigonal graph such that $G_k\cap\mathbb RP^1=G$''.

\end{itemize}
So, according to Theorem \ref{constr trig}, there exists a chain of
weighted combs starting at $w$ if and only if $A$ is realizable by a
real algebraic trigonal curve  on $\Sigma_n$.  \findemo

Theorem \ref{prohib comb} provides an algorithm to check whether an
$\mathcal L$-scheme is realizable by a real algebraic trigonal curve
 on $\Sigma_n$. In order to reduce computations, one
can use the following observations.
\begin{lemma}\label{l1}
Let $\mathfrak m$ be a closed comb,
and $\overline{\mathfrak m}$ one of its closures. Suppose that
$\mathfrak m=\mathfrak m_1g_{i_1}\mathfrak m_2g_{i_2}\mathfrak m_3$,
and that $g_{i_1}$ and $g_{i_2}$ are
joined in $\overline{\mathfrak m}$. Then the combs $\mathfrak m_1\mathfrak m_3$ and $\mathfrak m_2$ contain the same
number of
generators
$g_1$ (resp., $g_3$, $g_5$) and $g_2$ (resp., $g_4$, $g_6$).
\end{lemma}
\textit{Proof. }Straightforward.\findemo
\begin{lemma}\label{l2}
Let $(\mathfrak m,0,\beta,0)$ be an element of a chain of weighted combs. Then
it is possible to join each $g_1$ in $\mathfrak m$ to a $g_2$ in $\mathfrak m$ by
pairwise non-intersecting paths
on $\mathbb R^2$
such that if
$\mathfrak m=\mathfrak m_1g_{i_1}\mathfrak m_2g_{i_2}\mathfrak m_3$ with $g_{i_1}$ and $g_{i_2}$ joined, then
the combs $\mathfrak m_1\mathfrak m_3$ and $\mathfrak m_2$ contain the same number of
generators
$g_1$ and
$g_2$.
\end{lemma}
\textit{Proof. }Straightforward.\findemo
\begin{lemma}\label{l3}
Let $(\mathfrak m,\alpha,\beta,0)$ be an element of a chain of weighted combs,
where $\mathfrak m=\prod_{j=1}^k g_{i_j}$. Define the equivalence relation
$\sim$ on $\{j\textrm{ }|\textrm{ }i_j=1$ or $2\}$ as follows :
\begin{center}
$r\sim s$ if the cardinal of $\{j|r< j< s\textrm{ and } i_j=1,2,3\textrm{ or
    }4\}$ is odd.
\end{center}
Denote by $E_1^{\mathfrak m}$ and $E_2^{\mathfrak m}$ the two equivalence classes of
$\sim$. Then
$$|Card(E^{\mathfrak m}_1)-Card(E_2^{\mathfrak m})|\le \alpha.  $$
\end{lemma}
\textit{Proof. }Choose a chain of weighted combs $(w_i)_{1\le i\le k}$
which contains $(\mathfrak m,\alpha,\beta,0)$.
Let
$(\widetilde{\mathfrak m},0,\beta,0)$ be an element of this chain. Then there
exists $l\in\{1\ldots \alpha\}$ such that
\begin{center} $Card(E_1^{\widetilde{\mathfrak m}})=Card(E_1^{\mathfrak m})-l$ and $Card(E_2^{\widetilde{\mathfrak m}})=Card(E_2^{\mathfrak m})-\alpha +l$.
\end{center}
 It is obvious that in a closure of $\widetilde{\mathfrak m}$, an element of
 $E_1^{\widetilde{\mathfrak m}}$ has to be joined to an element of
 $E_{2}^{\widetilde{\mathfrak m}}$, hence  the cardinal of these two sets are
 equal.\findemo

The algorithm given by Theorem \ref{prohib comb} improved by Lemmas
\ref{l1}, \ref{l2}, and \ref{l3}, is quite efficient. It will allow us
in section \ref{algebraic state} to prohibit algebraically two
$\mathcal L$-schemes
pseudoholomorphically
realizable.

\section{Pseudoholomorphic statements}\label{pseudo state}

\subsection{Prohibitions for curves of bidegree $(3,1)$ on $\Sigma_2$}
We obtain all our results on symmetric curves of degree 7 via the
study of
the possible quotient curves. Hence, the first idea to prove that a
real scheme is not realizable by a symmetric curve is to show that no
quotient curve is admissible.  

\begin{lemma}\label{lem bezout1}
Let $X$ be a real curve  of bidegree $(3,1)$ on $\Sigma_2$, and consider  the encoding of the $\mathcal L$-scheme realized by $X$, as  defined in section \ref{braid} using the symbols 
$\subset$, $\supset $, $\times$ and $o$. Suppose that this encoding  contains a subsequence of the form $o_1^{k_1}o_2^{k_2}o_1^{k_3}o_2^{k_4}\ldots o_1^{k_{r-1}}o_2^{k_r}/$. Then $k_i=0$ for  $i\ge 4$ and $k_1k_3=0$.

Moreover, if $k_1+k_2+k_3\ge 3$, then $k_3=0$.
\end{lemma}
\textit{Proof.  }Suppose there exist $3$ ovals contradicting the
lemma. Then the base passing through these ovals intersects the curve
in at least $9$ points which contradicts the Bezout theorem.~\findemo

\begin{figure}[h]
\centering
\begin{tabular}{cccc}
\includegraphics[height=2cm, angle=0]{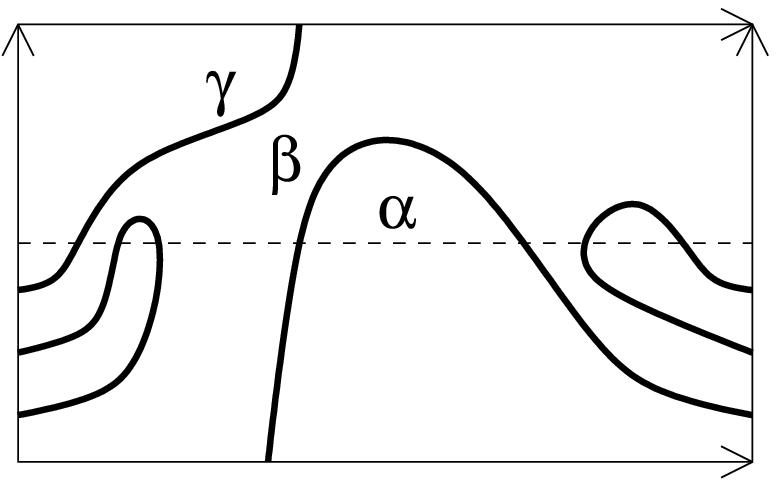}&
\includegraphics[height=2cm, angle=0]{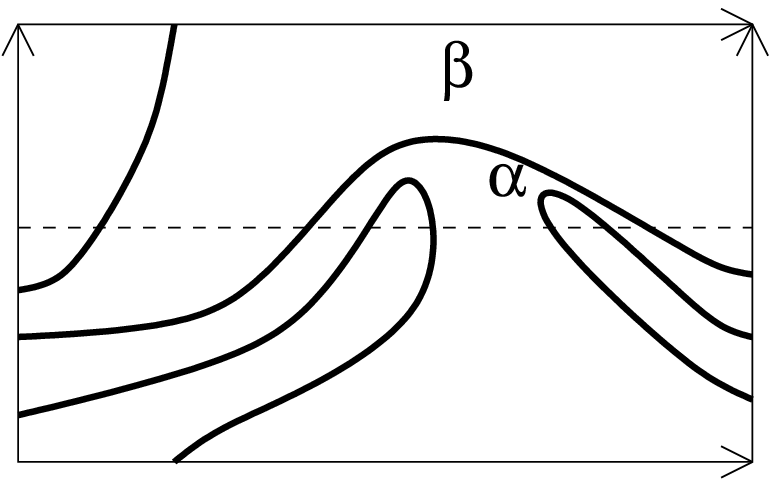}&
\includegraphics[height=2cm, angle=0]{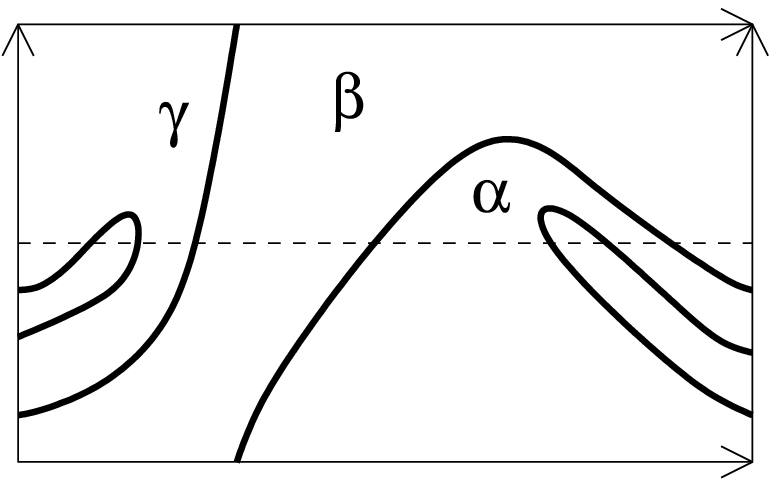}&
\includegraphics[height=2cm, angle=0]{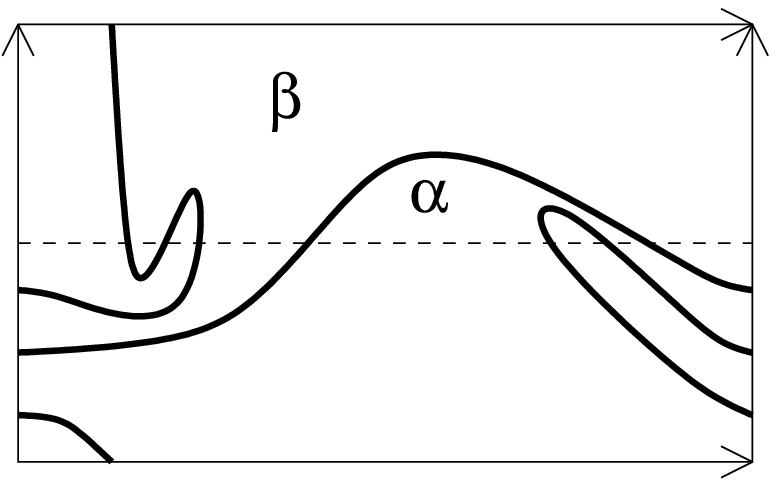}
\\ a)&b)&
c)&d)
\end{tabular}
\caption{}
\label{prohib even}
\end{figure}

\begin{figure}[h]
\centering
\begin{tabular}{ccc}
\includegraphics[height=2cm, angle=0]{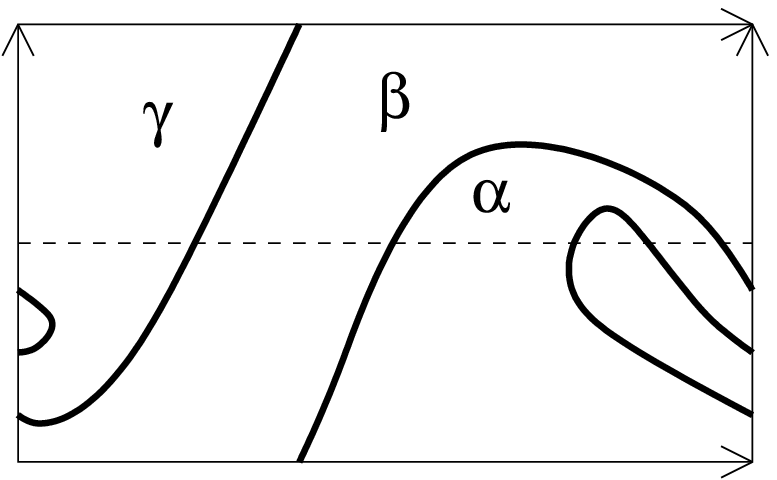}&
\includegraphics[height=2cm, angle=0]{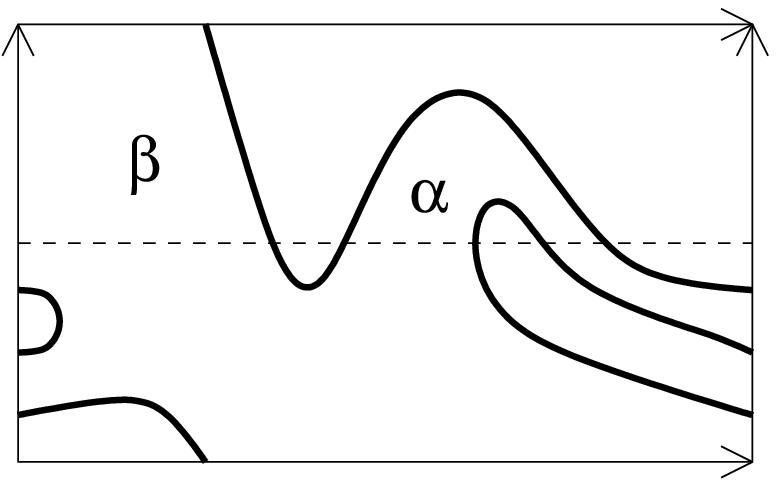}&
\includegraphics[height=2cm, angle=0]{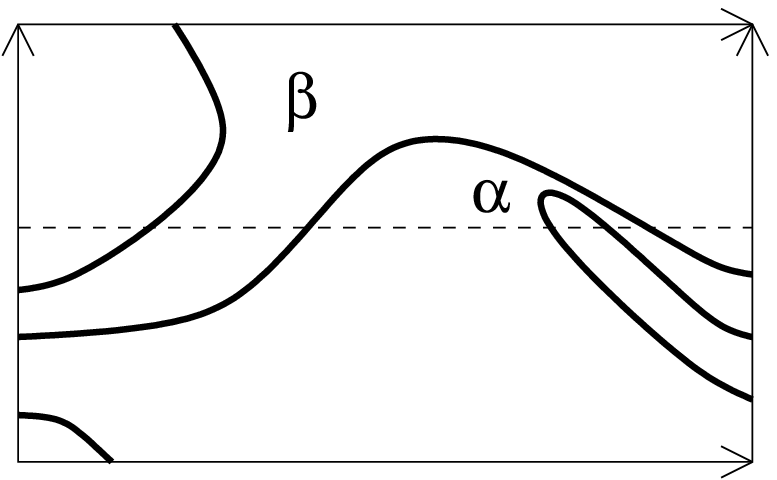}
\\ a)&b)&c)
\end{tabular}

\caption{}
\label{prohib odd2}
\end{figure}

In
Figures \ref{prohib even} and \ref{prohib odd2}, the dashed line
represents the base $\{y=0 \}$.

\begin{lemma}\label{pr1}
If there exists a pseudoholomorphic curve of bidegree $(3,1)$ on
$\Sigma_2$ realizing the  \mbox{$\mathcal L$-scheme} depicted in Figure
\ref{prohib even}a) with $\alpha+\beta+\gamma=6$, then
$(\alpha,\beta)=(0,1)$, $(0,5)$, $(1,5)$ or $(5,1)$.
\end{lemma}
\textit{Proof. }The corresponding braid is
$$b^1_{\alpha,\beta,\gamma}=\sigma_1^{-1}\sigma_2^{-1}\sigma_1\sigma_2^{-\gamma}\sigma_1^{-1}\sigma_2^2\sigma_1\sigma_2^{-\beta}\sigma_1^{-1}\sigma_2\sigma_1^{-\alpha}\sigma_2^{-1}
\sigma_1\Delta_3^2.$$
According to the Bezout theorem, $\alpha$ and $\gamma$
cannot be simultaneously non null. Moreover, it is clear by symmetry
of the $\mathcal L$-scheme that $b^1_{\alpha,\beta,\gamma}$ is
quasipositive if and only if $b^1_{\gamma,\beta,\alpha}$ is quasipositive.
Some calculations give the
following Alexander 
polynomials
\begin{center}
\begin{tabular}{lll}
$p^1_{0,6,0}=p^1_{6,0,0}=(t-1)(t^4-t^3+t^2-t+1),$&
$p^1_{2,4,0}=p^1_{4,2,0}=(t-1)(t^2-t+1),$
\\$p^1_{3,3,0}=(t-1)^3$,&
$p^1_{1,5,0}=p^1_{5,1,0}=0.$
\end{tabular}
\end{center}
We have  $e(b^1_{\alpha,\beta,\gamma})=7-(\alpha+\beta+\gamma)$ so according to
 Proposition \ref{alex}, the braid is not quasipositive
as soon as its Alexander polynomial is not identically zero, and the
 lemma is proved.~\findemo

\begin{lemma}\label{pr2}
If there exists a pseudoholomorphic curve of bidegree $(3,1)$ on
$\Sigma_2$ realizing one of the  $\mathcal L$-schemes depicted in Figures
\ref{prohib even}b) and d) or in Figure \ref{prohib odd2}c) with $\alpha+\beta=6$, then
$\alpha=1$ or $5$.
\end{lemma}
\textit{Proof. }First, note that the three  $\mathcal L$-schemes give
rise to the same braid which is
$$b^2_{\alpha,\beta}=
\sigma_1\sigma_2^{-\beta}\sigma_1^{-1}\sigma_2\sigma_1^{-\alpha}\Delta_3^2.$$
Some calculations show that $p^2_{\alpha,\beta}=p^1_{\alpha,\beta,0}$
as soon as $\alpha+\beta=6$. We have  $e(b^2_{\alpha,\beta})=7-(\alpha+\beta)$ so according to
 Proposition \ref{alex}, the braid is not quasipositive
as soon as its Alexander polynomial is not identically zero, and the
 lemma is proved.~\findemo

\begin{lemma}\label{pr3}
If there exists a pseudoholomorphic curve of bidegree $(3,1)$ on
$\Sigma_2$ realizing one of the  $\mathcal L$-schemes depicted in Figure
\ref{prohib even}c) or in Figure \ref{prohib odd2}a) with $\alpha+\beta+\gamma=6$, then
$(\alpha,\beta)=(0,4)$, $(0,0)$, $(2,4)$ or $(6,0)$.
\end{lemma}
\textit{Proof. }First, note the two  $\mathcal L$-schemes are the
same. The corresponding braid is
$$b^3_{\alpha,\beta,\gamma}= \sigma_2^{-(1+\gamma)}\sigma_1^{-1}\sigma_2^{2}\sigma_1\sigma_2^{-\beta}\sigma_1^{-1}\sigma_2\sigma_1^{-\alpha}\Delta_3^2.$$
According to the Bezout theorem, $\alpha$ and $\gamma$
cannot be simultaneously non null. Moreover, it is clear by symmetry
of the $\mathcal L$-scheme that $b^3_{\alpha,\beta,\gamma}$ is
quasipositive if and only if $b^3_{\gamma,\beta,\alpha}$ is quasipositive.
Some calculations give the
following Alexander 
polynomials
\begin{center}
\begin{tabular}{lll}
$p^3_{0,6,0}=(t^2+t+1)(t^2-t+1)(t-1)^3$,&
$p^3_{1,5,0}=(t-1)(t^4-t^3+t^2-t+1),$&
\\$p^3_{3,3,0}=p^3_{5,1,0}=(t-1)(t^2-t+1),$&
$p^3_{4,2,0}=(t-1)^3$,
\\
$p^3_{2,4,0}=p^3_{6,0,0}=0.$
\end{tabular}
\end{center}
We have  $e(b^3_{\alpha,\beta,\gamma})=7-(\alpha+\beta+\gamma)$ so according to
 Proposition \ref{alex}, the braid is not quasipositive
as soon as its Alexander polynomial is not identically zero, and the
 lemma is proved.~\findemo

\begin{lemma}\label{pr4}
There does not exist a pseudoholomorphic curve of bidegree $(3,1)$ on
$\Sigma_2$ realizing the  $\mathcal L$-scheme depicted in Figure
\ref{prohib odd2}b) with $\alpha+\beta=6$,
\end{lemma}
\textit{Proof. }The corresponding braid and its Alexander polynomial are
$$b^4= \sigma_2^{-7}\sigma_1\sigma_2\Delta_3^2\textrm{ and }
p^4=(t-1)(t^4-t^3+t^2-t+1).$$ 
We have $e(b^4)=1$, so according to  Proposition \ref{alex},
this braid is not
quasipositive.
~\findemo

\begin{prop}\label{first M}
The following real schemes cannot be realized by symmetric pseudoholomorphic curves of degree $7$ on $\mathbb RP^2$:
\begin{center}
$\langle J\amalg \beta\amalg 1\langle \alpha\rangle \rangle$ with $\alpha=6,7,8,9$ and $\beta=14-\alpha$,  $\alpha=7,9$ and $\beta=13-\alpha$.
\end{center}
\end{prop}
\textit{Proof. }According to Proposition \ref{max sym}, \ref{inter max}, the Bezout theorem, Lemmas \ref{lem bezout1} and \ref{convexe},
 the only possibilities for the $\mathcal L$-schemes of the quotient curves
  are depicted in Figure \ref{prohib even}a), c) and d) with
 $(\alpha,\beta)=(3,3)$ or $(4,2)$ and in Figure \ref{prohib even}b) with
 $(\alpha,\beta)=(3,3)$ or $(2,4)$  for the four $M$-curves, and
 in Figure \ref{prohib odd2} with $(\alpha,\beta)=(3,3)$, for the two $(M-1)$-curves.
Now the proposition follows from Lemmas \ref{pr1}, \ref{pr2}, \ref{pr3}
and \ref{pr4}.~\findemo

\subsection{Prohibitions for reducible curves of bidegree $(4,1)$ on $\Sigma_2$}\label{prohib1}
Here, we have to study more carefully the $\mathcal L$-schemes
realized  by the hypothetical quotient curves. Indeed, some of them are
realized by curves of bidegree $(3,1)$ in $\Sigma_2$. To prohibit the
symmetric curves, we have to take into account the mutual position of the
quotient curves and some base of $\Sigma_2$.

\begin{figure}[h]
      \centering
 \begin{tabular}{ccc}
\includegraphics[height=2cm, angle=0]{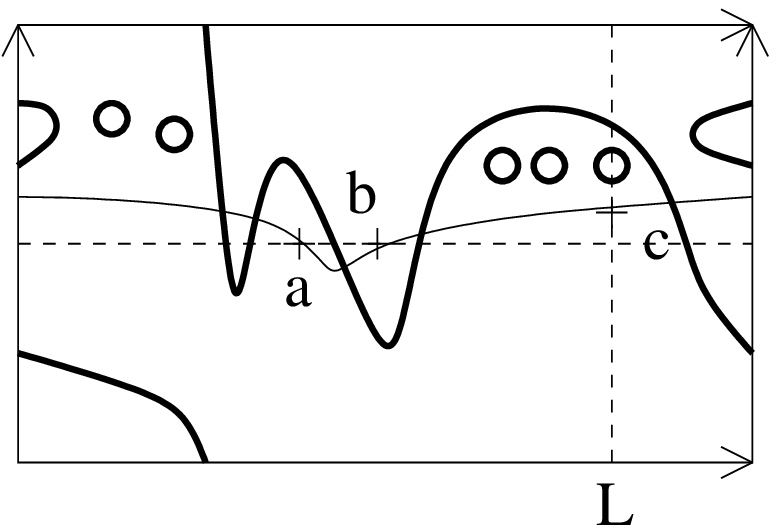}  &
\includegraphics[height=2cm, angle=0]{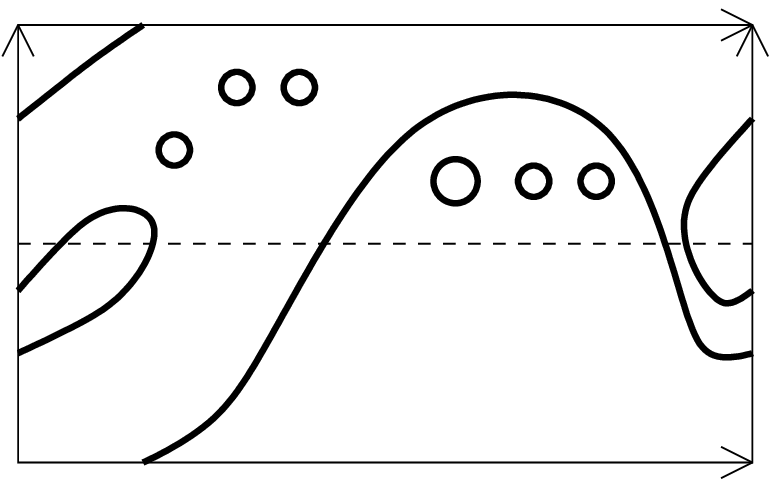}
 &\includegraphics[height=2cm,angle=0]{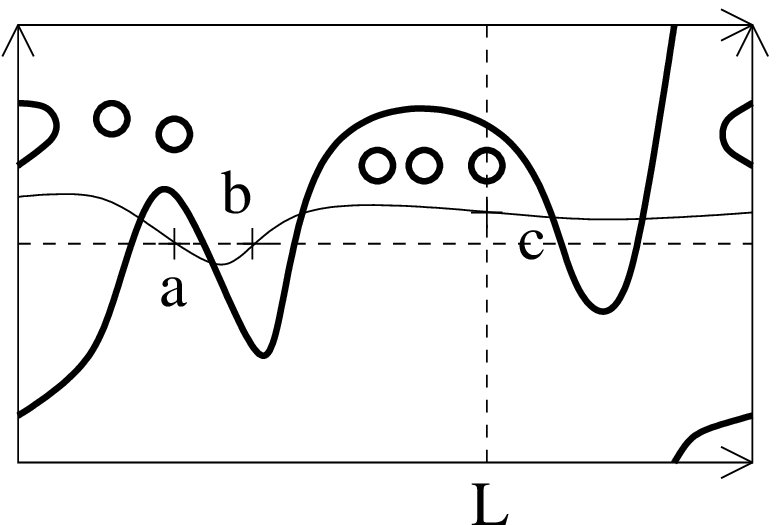}   \\
a)&b)&c)\\
 \end{tabular}
\caption{}
 \label{base1}
\end{figure}

\begin{prop}\label{prohib thm 1}
The real scheme $\langle J\amalg 7\amalg 1\langle 6\rangle \rangle$ is not realizable by a nonsingular symmetric pseudoholomorphic curve of degree $7$ on $\mathbb RP^2$.
\end{prop}
\textit{Proof. } The possible  $\mathcal L$-schemes realized by
the corresponding quotient curves with respect to the base $\{y=0\}$
are depicted in Figure \ref{base1}. 

The $\mathcal L$-scheme   in Figures \ref{base1}b) cannot be realized pseudoholomorphically. The braid corresponding to this
$\mathcal L$-scheme and its Alexander polynomial are~:
\begin{center}
\begin{tabular}{ll}
$ b^5=\sigma_1^{-4}\sigma_2^{2}\sigma_1^{-3}\sigma_2^{-1}\sigma_1\Delta_3^2$&
 and $p^5=(t-1)^3 $
\end{tabular}
\end{center}
The Alexander polynomial is not null although $e(b^5)=1$, so according to Proposition
\ref{alex}, the braid is  not
quasipositive.

\begin{figure}[h]
      \centering
 \begin{tabular}{cc cc}
\includegraphics[height=2cm, angle=0]{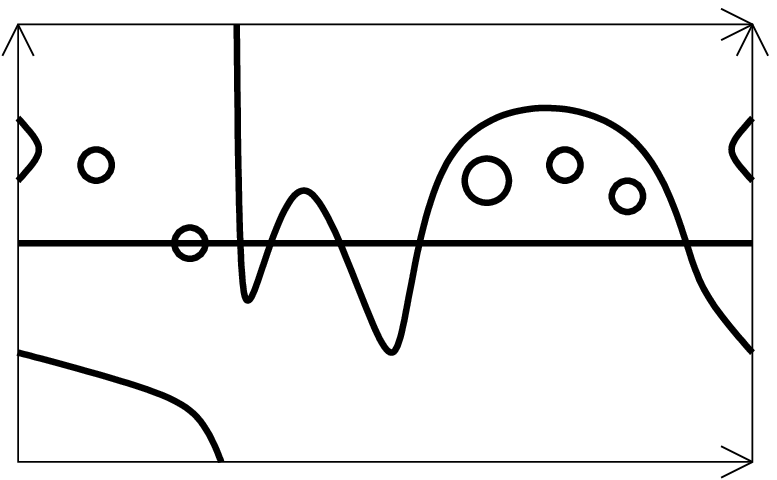}
 &\includegraphics[height=2cm,angle=0]{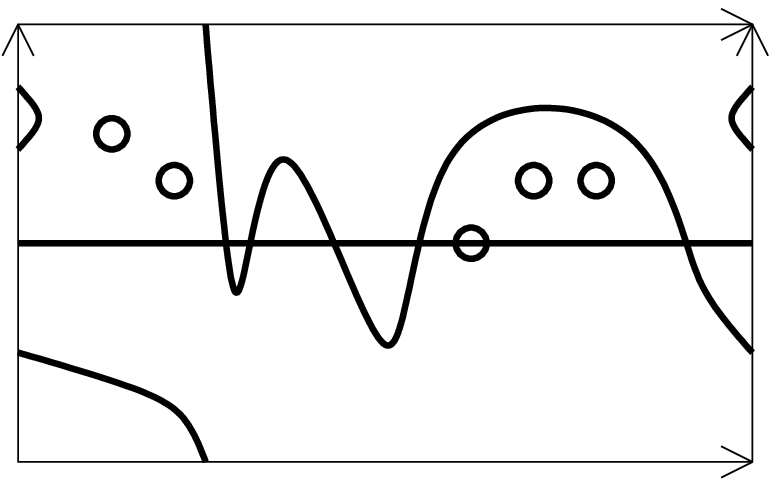}   &
\includegraphics[height=2cm, angle=0]{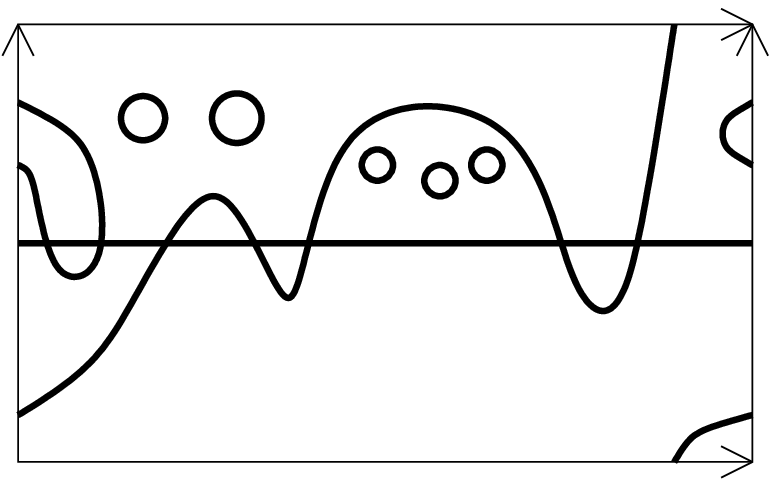} &
\includegraphics[height=2cm,angle=0]{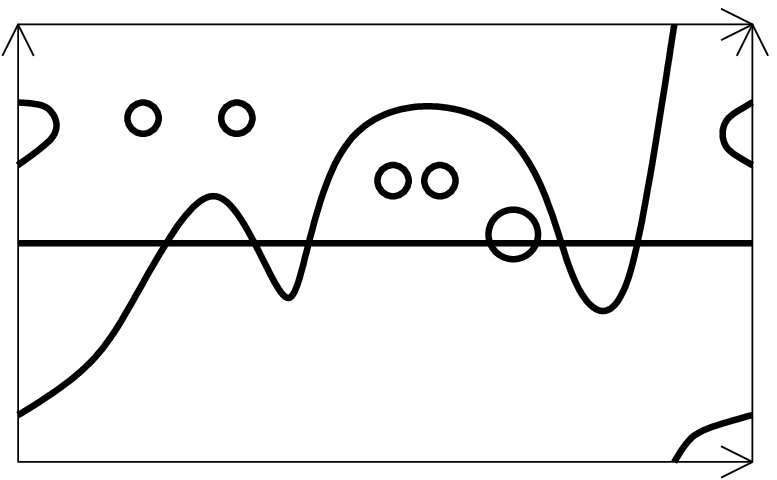}   \\
a)&b)&c)&d)\\
 \end{tabular}
\caption{}
 \label{base2}
\end{figure}

Consider the base $H$ passing through the points
$a$, $b$ and $c$ in Figures \ref{base1}a) and c), where $c$ is a point of the fiber $L$. If the point $c$ varies on $L$ from $0$
to $\infty$ in $\{y\ge0\}$, then, because of the choice of $L$, for
some $c$, the base $H$ passes through an oval. Let  $H$ be a base
which passes through the first oval
we meet as $c$ varies from $0$
to $\infty$. The only possible mutual arrangements for $H$
and the quotient curves which do not contradict the Bezout theorem are shown in Figure \ref{base2}.
The corresponding braids
are~:
\begin{quote}
$b^6= \sigma_3^{-2}\sigma_2^{-2}\sigma_3^{-1}\sigma_1\sigma_2^{2}\sigma_1^{-4}\sigma_2^{-1}\sigma_3
\sigma_2^{-3}\sigma_3^{-1}\sigma_2\sigma_1^{-1}\Delta_4^2$,
\\$b^7= \sigma_3^{-3}\sigma_1\sigma_2^{2}\sigma_1^{-4}\sigma_2^{-1}\sigma_3\sigma_1^{-2}
\sigma_2^{-3}\sigma_3^{-1}\sigma_2\sigma_1^{-1}\Delta_4^2$,
\\$b^8=  \sigma_2^{-2}\sigma_3^{-3}\sigma_1^{-3}\sigma_2^{-1}\sigma_3\sigma_2^{-3}\sigma_3^{-1}\sigma_2\sigma_1^{-2}
\sigma_2^{-1}\sigma_3^{2}\sigma_2\sigma_1\Delta_4^2$,
\\$b^{9}= \sigma_3^{-3}\sigma_1^{-3}\sigma_2^{-1}\sigma_3\sigma_2^{-2}\sigma_1^{-2}
\sigma_2^{-1}\sigma_3^{-1}\sigma_2\sigma_1^{-2}\sigma_2^{-1}\sigma_3^2\sigma_2\sigma_1\Delta_4^2$.
\end{quote}
The computation of the corresponding Alexander polynomials gives~:
\begin{center}
\begin{tabular}{ll}
$p^6=(t^2-t+1)(t^6-3t^5+6t^4-5t^3+6t^2-3t+1)(-1+t)^3 $,&
$p^8=(-1+t)^7  $,\\
$p^7=(2t^4-2t^3+3t^2-2t+2)(t^2-t+1)^2(-1+t)^3 $,&
$p^{9}=(t^2-t+1)(-1+t)^3 $.
\end{tabular}
\end{center}
In each case, $e(b)=2$, so according to Proposition
\ref{alex}, none of these braids is quasipositive.
~\findemo

\subsection{Restrictions for dividing symmetric curves}

In this section, we use extensively the observation made at the end of
section \ref{general}~: useful information on the complex
orientations of a symmetric curve can be extracted from the topology of
its quotient curve. 

 \begin{figure}[h]
      \centering
 \begin{tabular}{cccc}
\includegraphics[height=2cm, angle=0]{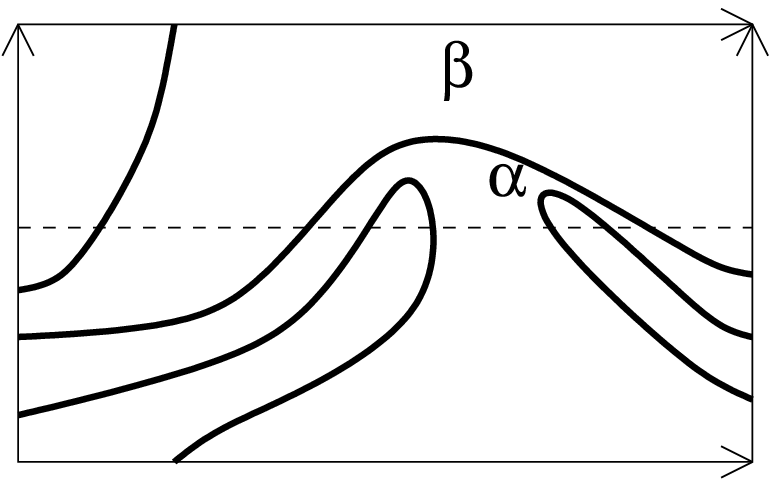}    &
 \includegraphics[height=2cm, angle=0]{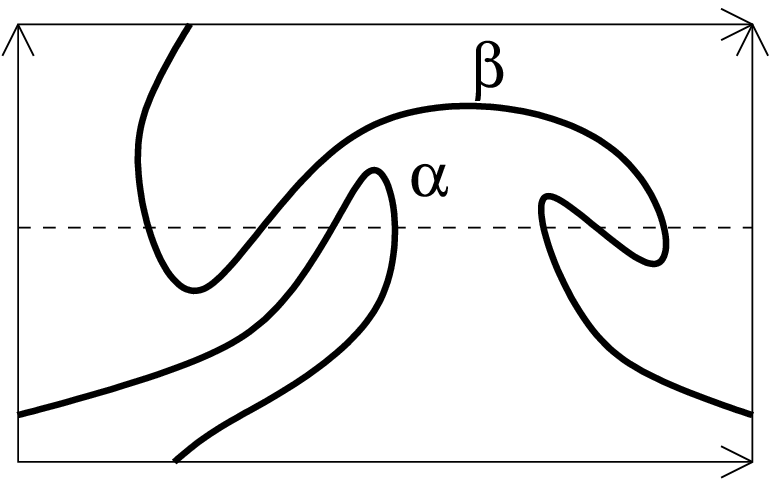} &
    \includegraphics[height=2cm, angle=0]{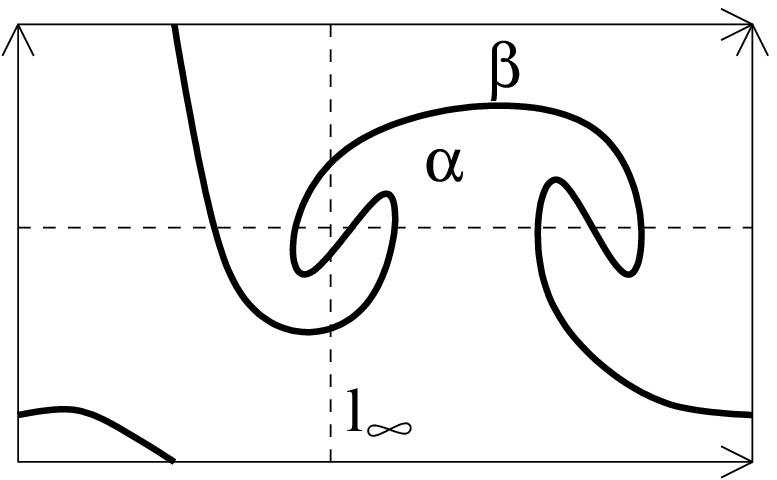}&
\includegraphics[height=2cm, angle=0]{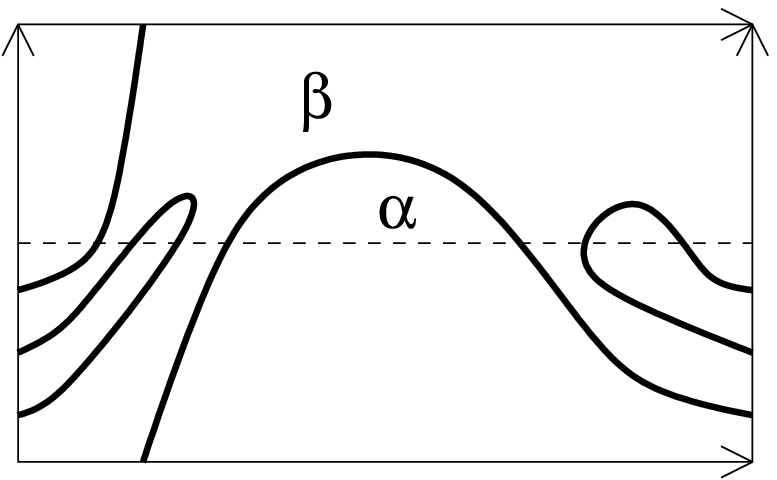}
\\ a)&b)&c)&d)\\
 \includegraphics[height=2cm, angle=0]{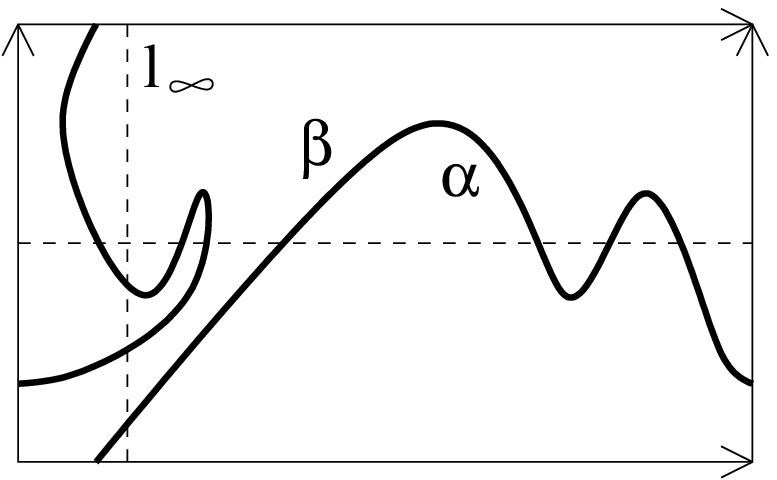} &
    \includegraphics[height=2cm, angle=0]{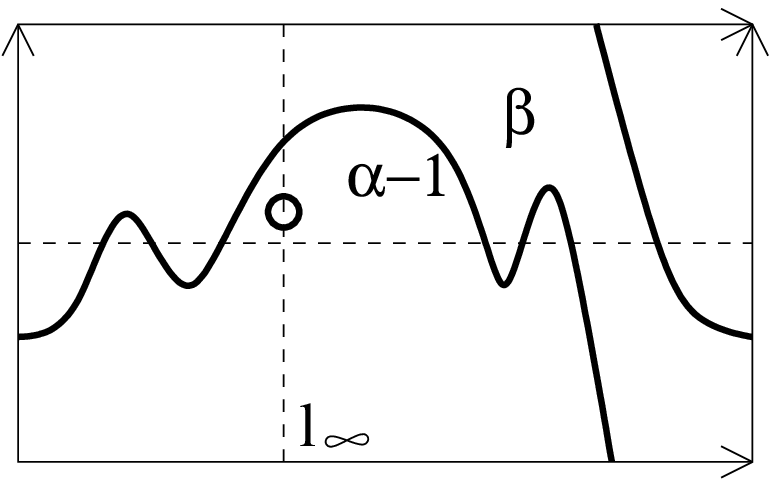}&
\includegraphics[height=2cm, angle=0]{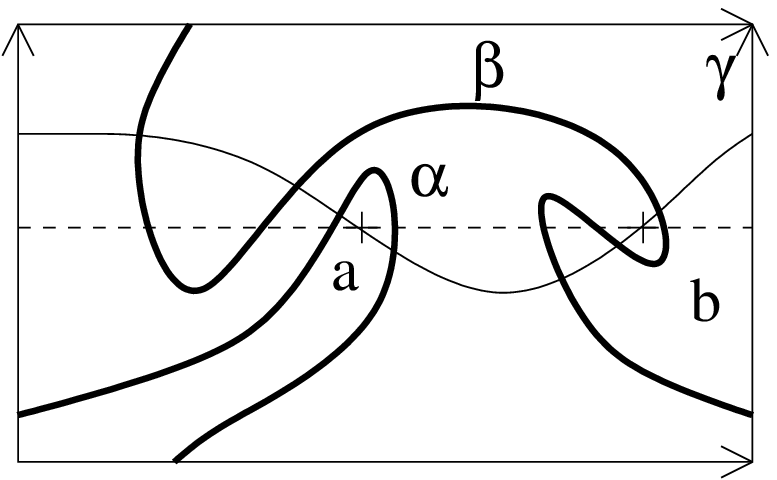}    &
 \includegraphics[height=2cm, angle=0]{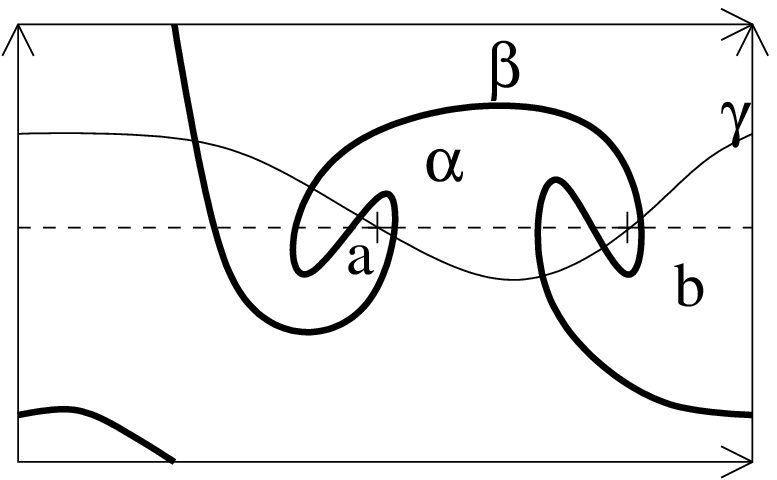}
\\ e)&f)&g)&h)
 \end{tabular}
\caption{}
 \label{quotient type I}
\end{figure}
\begin{lemma}\label{lem1}
There does not exist a symmetric dividing pseudoholomorphic curve of degree $7$ on $\mathbb RP^2$ with a quotient curve
realizing the $\mathcal L$-scheme depicted in Figures \ref{quotient type I}a) and d) with $\alpha+\beta$ odd.
\end{lemma}
\textit{Proof. }Such a quotient curve is of type $II$ because it is an $(M-2i-1)$-curve. The mirror curve of the
initial symmetric curve is a nest of depth $3$ with an odd component, and so is of type $I$. Thus, according to Proposition
\ref{quotient dividing}, the initial symmetric curve cannot be of type $I$.
~\findemo

\begin{lemma}\label{lem2}
There does not exist a symmetric dividing pseudoholomorphic curve of degree $7$ on $\mathbb RP^2$ with a quotient curve
realizing the $\mathcal L$-scheme depicted in Figure \ref{quotient type I}b) with $\alpha+\beta$ odd.
\end{lemma}
\textit{Proof. }According to the Fiedler orientation alternating rule on symmetric curves corresponding to these quotient curves,
 the symmetric curves
cannot be of type $I$ if $\alpha+\beta$ is odd, as the two invariant empty ovals have opposite orientations.
~\findemo

\begin{lemma}\label{lem3}
There does not exist a symmetric dividing pseudoholomorphic curve of degree $7$ on $\mathbb RP^2$ with a quotient curve
realizing the $\mathcal L$-scheme depicted in Figure \ref{quotient type I}c) with
$(\alpha,\beta)=(4,1)$, $(3,2)$, $(2,3)$, and $(1,2).$
\end{lemma}
\textit{Proof. }According to the Fiedler orientation alternating rule on symmetric curves corresponding to these quotient curves,
the three invariant ovals are positive and  we have
$\Lambda_+-\Lambda_-=1$,
$\Pi_+-\Pi_-=0$ if $\alpha$ is odd, and
$\Pi_+-\Pi_-=-2$ if $\alpha$ is even.
Thus, the Rokhlin-Mischachev orientation formula is fulfilled only for $(\alpha,\beta)=(3,2)$. Choose $l_\infty$ as depicted in Figure \ref{quotient type I}c).
The braid corresponding to the quotient curve and its Alexander
polynomial are~: 
\begin{center}
\begin{tabular}{ll}
$\sigma_1^{-4}\sigma_2^{-1}\sigma_1\sigma_2^{-2}\sigma_1\Delta_3^2$&
and $(-1+t)^3 $.
\end{tabular}
\end{center}
Since $e(b)=1$, according to Proposition \ref{alex}, the braid is not quasipositive.
~\findemo
\begin{lemma}\label{lem4}
There does not exist a symmetric dividing pseudoholomorphic curve of degree $7$ on $\mathbb RP^2$ with a quotient curve
realizing the $\mathcal L$-scheme depicted in Figure \ref{quotient type I}e) with
$(\alpha,\beta)=(5,0)$, $(4,1)$, $(3,2)$, $(2,3)$, and $(2,1).$
\end{lemma}
\textit{Proof. }According to the Fiedler orientation alternating rule on symmetric curves corresponding to these quotient curves,
 the two invariant empty ovals have opposite orientations, the non-empty oval is negative, and we have
$\Lambda_+-\Lambda_-=1$,
$\Pi_+-\Pi_-=0$ if $\alpha$ is even, and
$\Pi_+-\Pi_-=-2$ if $\alpha$ is odd.
Hence, the Rokhlin-Mischachev orientation formula is fulfilled only
 for $(\alpha,\beta)=(4,1)$ or $(2,3)$.  Choose $l_\infty$ as depicted
 in Figure \ref{quotient type I}e). Then the braids corresponding to
 the union of the quotient curves and of the base $\{y=0\}$ 
are~:
\begin{quote}
$b^{10}_{\alpha,\beta}=\sigma_3^{-1}\sigma_2^{-1}\sigma_1^{-1}\sigma_3^{-(1+\beta)}\sigma_2^{-1}\sigma_3
\sigma_2^{-\alpha}\sigma_3^{-1}
\sigma_2\sigma_1^{-3}\sigma_3\sigma_2\sigma_1\sigma_3^{-1}\Delta_4^2$.
\end{quote}
The computation of the corresponding Alexander polynomials gives
\begin{center}
\begin{tabular}{ll}
$p^{10}_{4,1}=(t^2+1)(t^2-t+1)(-1+t)^3$,
&$p^{10}_{2,3}=(t^4-2t^3+4t^2-2t+1)(-1+t)^3$
\end{tabular}
\end{center}
In each case we have $e(b)=2$, so according to Proposition \ref{alex}, both braids are not quasipositive.~\findemo

\begin{lemma}\label{lem5}
There does not exist a symmetric dividing pseudoholomorphic curve of degree $7$ on $\mathbb RP^2$ with a quotient curve
realizing the $\mathcal L$-scheme depicted in Figure \ref{quotient type I}f) with $(\alpha,\beta)=(5,0)$, $(3,2)$, and  $(2,3)$.
\end{lemma}
\textit{Proof. }Choose $l_\infty$ as depicted in Figure \ref{quotient
  type I}f).
Then the braids corresponding to
 the union of the quotient curves and of the base $\{y=0\}$ 
are~:
\begin{quote}
$b^{11}_{\alpha,\beta}=\sigma_2^{-\alpha}\sigma_3^{-1}\sigma_2\sigma_3^{-\beta}\sigma_1^{-3}
\sigma_1\sigma_2^{2}\sigma_1^{-4}
\sigma_2^{-1}\sigma_3\Delta_4^2$.
\end{quote}
The computation of the determinant gives
$ 976$ for $b^{11}_{5,0}$
and
$592$ for $b^{11}_{3,2}$ which
are not squares in $\mathbb Z$ although $e(b)=3$. So according to Proposition \ref{double alex}, these two braids are not quasipositive.

The computation of the Alexander polynomial of $b^{11}_{2,3}$ gives
\begin{quote}
$(t^2+1)(t^6-5t^5+12t^4-14t^3+12t^2-5t+1)(-1+t)^2$.
\end{quote}
The number $i$ is a simple root of this polynomial and $
e(b)=3$. Thus, according to Proposition \ref{square}, this braid is not quasipositive.
~\findemo

\begin{prop}\label{restr thm 4 3}
The real schemes $\langle J\amalg 2\amalg 1\langle 10\rangle \rangle$ and  $\langle J\amalg 6\amalg 1\langle 6\rangle \rangle$ are not realizable by nonsingular
 symmetric dividing pseudoholomorphic curves of degree $7$ on $\mathbb RP^2$.
\end{prop}
\textit{Proof. }According to the Bezout theorem, Proposition \ref{oval
  quotient}, Lemma \ref{lem bezout1} and Lemma \ref{convexe},
the only possibilities for the  $\mathcal L$-scheme of
the quotient curve of such a dividing symmetric curve of degree $7$ on $\mathbb RP^2$ are depicted
in Figures \ref{quotient type I}a), g) and h) with $(\alpha,\beta+\gamma)=(4,1)\textrm{ and }(2,3)$
and in Figures  \ref{quotient type I}d), e) and f) with
$(\alpha,\beta)=(5,0),\textrm{ and }(3,2)$.
If a curve of bidegree $(3,1)$ realizes one of the two $\mathcal L$-schemes depicted in Figures \ref{quotient type I}g) and h) then $\gamma=0$. Otherwise, the base passing through
the points $a$ and $b$ and through an oval $\gamma$ intersects
the quotient curve in more than $7$ points, which contradicts the Bezout theorem.

The remaining quotient curves have been prohibited in Lemmas \ref{lem1}, \ref{lem2}, \ref{lem3}, \ref{lem4}, and \ref{lem5}~\findemo

\begin{prop}\label{only possibility}
If  a nonsingular symmetric dividing pseudoholomorphic curve of degree $7$ on $\mathbb RP^2$ realizes
the real scheme $\langle J\amalg 8\amalg 1\langle 4\rangle \rangle$, then the  $\mathcal L$-scheme of
its quotient curve is as depicted in Figure \ref{quotient type I}c)
with $(\alpha,\beta)=(1,4)$.

If  a nonsingular symmetric dividing pseudoholomorphic curve of degree $7$ on $\mathbb RP^2$ realizes
the real scheme $\langle J\amalg 2\beta+2\amalg 1\langle 2\alpha\rangle \rangle$ with $(\alpha,\beta)=(4,1)$ or $(2,1)$,
 then the  $\mathcal L$-scheme of
its quotient curve is as depicted in Figure \ref{quotient type I}f).
\end{prop}
\textit{Proof. }The proof is the same as for the previous proposition.
~\findemo

\begin{figure}[h]
      \centering
 \begin{tabular}{ccc}
\includegraphics[height=2.5cm, angle=0]{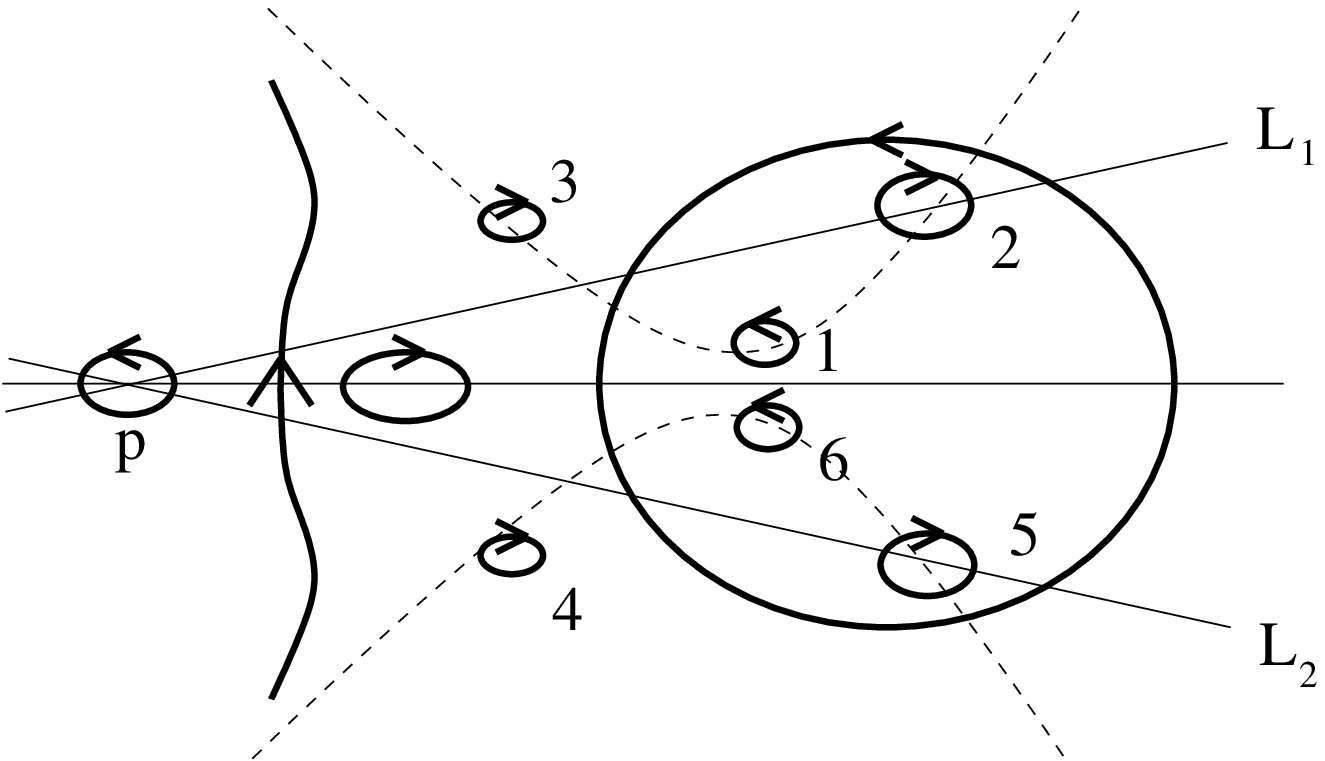}    &
\includegraphics[height=2cm, angle=0]{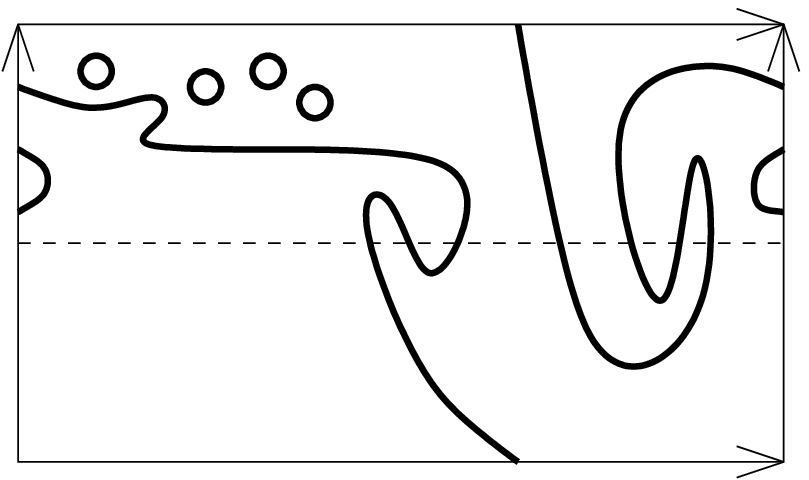}    &
\includegraphics[height=2cm, angle=0]{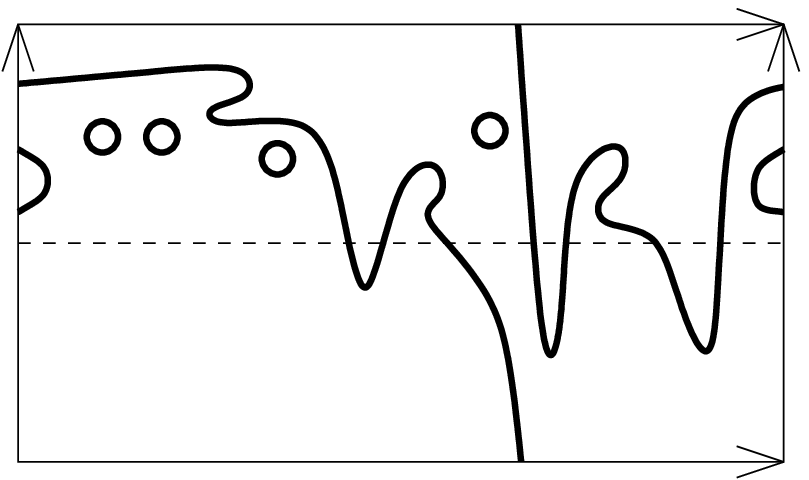}
\\ a)&b)&c)
 \end{tabular}
\caption{}
 \label{orientations 4 4}
\end{figure}

\begin{prop}\label{non 4 4 I}
The real scheme $\langle J\amalg 4\amalg 1\langle 4\rangle \rangle$ is not realizable by
 nonsingular symmetric dividing pseudoholomorphic curves of degree $7$ on $\mathbb RP^2$.
\end{prop}
\textit{Proof. }Suppose that there exists a pseudoholomorphic curve $C$
which contradicts Proposition \ref{non 4 4 I}.
Then, according to Proposition
\ref{only possibility}, its quotient curve is as depicted in Figure
\ref{quotient type I}f) with $(\alpha,\beta)=(2,1)$. Using the Fiedler
orientation alternating rule, and denoting by $\epsilon$ the sign of
the two non-invariant outer ovals of $C$, we have 
$\Lambda_+-\Lambda_-=-1+2\epsilon$ and
$\Pi_+-\Pi_-=0$.
Thus, the Rokhlin-Mischachev orientation formula is fulfilled only if
$\epsilon=-1$. Hence, one of the two complex orientations of the curve
is as depicted in Figure \ref{orientations 4 4}a). Using once again
the Fiedler orientation alternating rule, we see that the pencil of
lines through the point $p$ induces a cyclic  order on the $6$ non-invariant
ovals of $C$ as depicted in Figure \ref{orientations 4
  4}a). So the ovals $4$ and $1$ are not on the same connected
component of $\mathbb RP^2\setminus (L_1\cup L_2)$. A symmetric conic
passing through all the non-invariant ovals (in bold line in Figure
\ref{orientations 4 4}a)) intersects $C$ in at least $18$ points what
contradicts the Bezout Theorem.~\findemo

\subsection{Constructions}
\begin{prop}\label{pseudo}
There exist nonsingular real pseudoholomorphic curves of bidegree $(3,1)$ on $\Sigma_2$ such that the $\mathcal L$-scheme realized by the union of this curve and a base is
as shown in Figures \ref{orientations 4 4}b) and c). In particular, all the real tangency points of the curve with the pencil $\mathcal L$ are above
the base $\{y=0\}$.
\end{prop}
\textit{Proof. }The braids associated to these $\mathcal L$-schemes are
\begin{quote}
$b^{12}=\sigma_2^{-1}\sigma_3^{-1}\sigma_2 \sigma_3^{-1}\sigma_2^{-1}\sigma_3^{-3} \sigma_2^{-1}\sigma_3
 \sigma_1^{-1}\sigma_2^{-2}\sigma_3^{-1}\sigma_1\sigma_2^{2}\sigma_1^{-1}\sigma_2^{-2}\sigma_1^{-1}
\sigma_2^{-1}\Delta_4^2$,
\\$b^{13}=\sigma_2^{-3}\sigma_3^{-1}\sigma_2^{-1} \sigma_3\sigma_2^{-1}\sigma_3^{-1} \sigma_2\sigma_1^{-2} \sigma_2^{-1}\sigma_1^{-1}\sigma_3^{-1}\sigma_1\sigma_2^{2}\sigma_1^{-2}\sigma_2^{-1}\sigma_1^{-2}\sigma_2^{-1}\sigma_3
\Delta_4^2$.
\end{quote}
Using the Garside normal form (see \cite{Gar} or \cite{Jaq}), we see that these braids are trivial, so quasipositive.~\findemo

\begin{figure}[h]
      \centering
 \begin{tabular}{cc}
\includegraphics[height=4cm, angle=0]{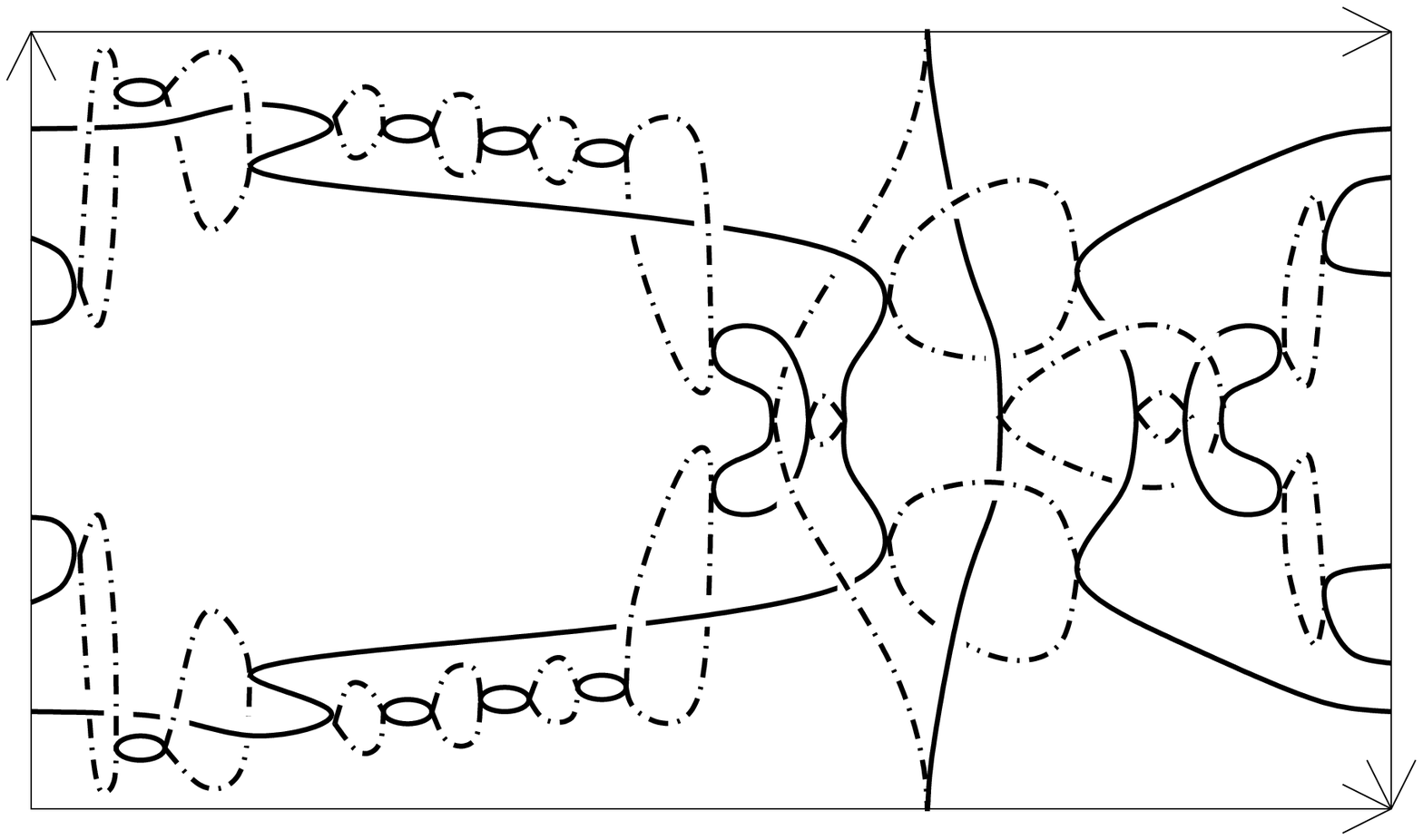}    &
\includegraphics[height=4cm, angle=0]{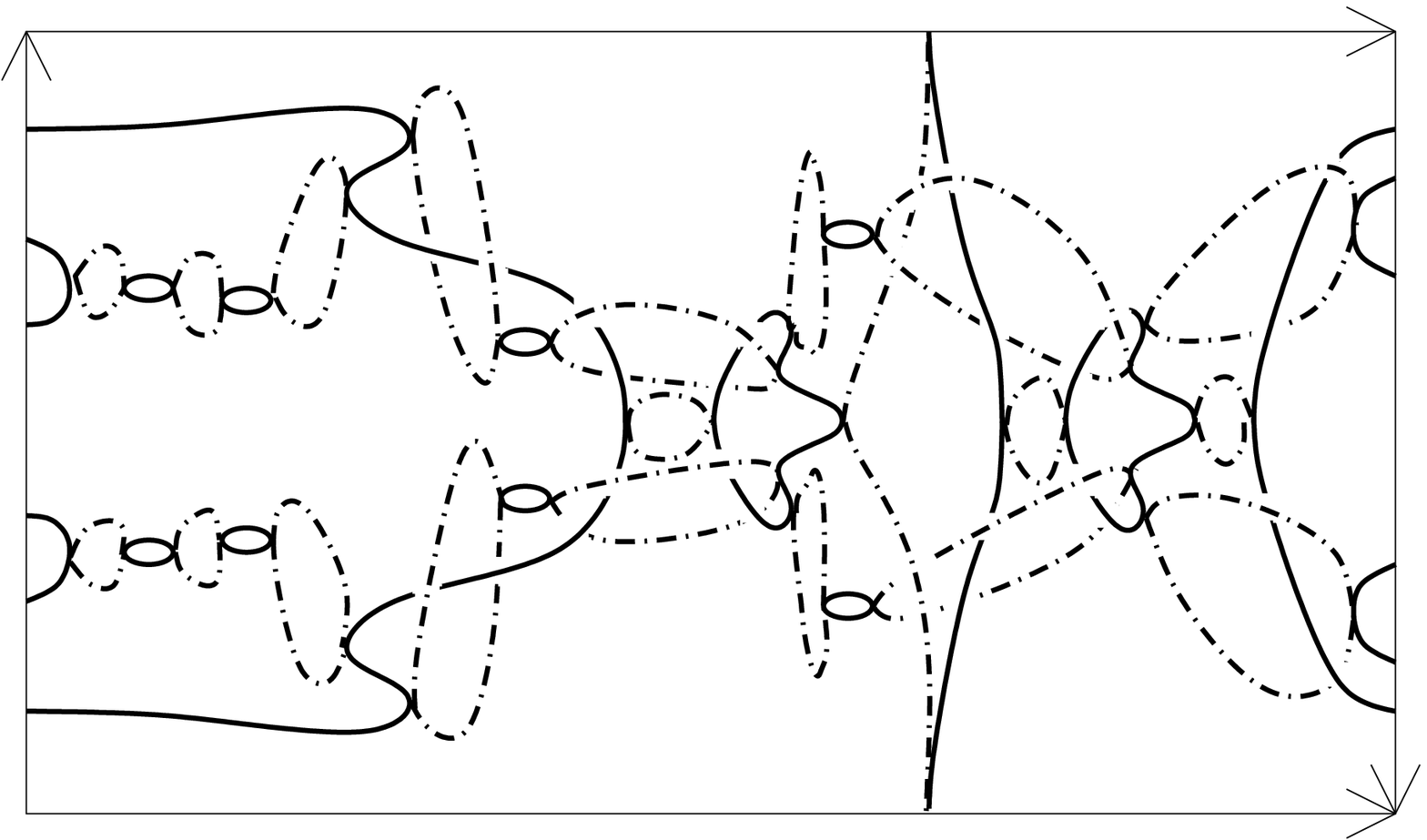}
\\ a)&b)
 \end{tabular}
\caption{}
 \label{pseudo constr}
\end{figure}
\begin{figure}[h]
      \centering
\begin{tabular}{c}
\includegraphics[height=3cm, angle=0]{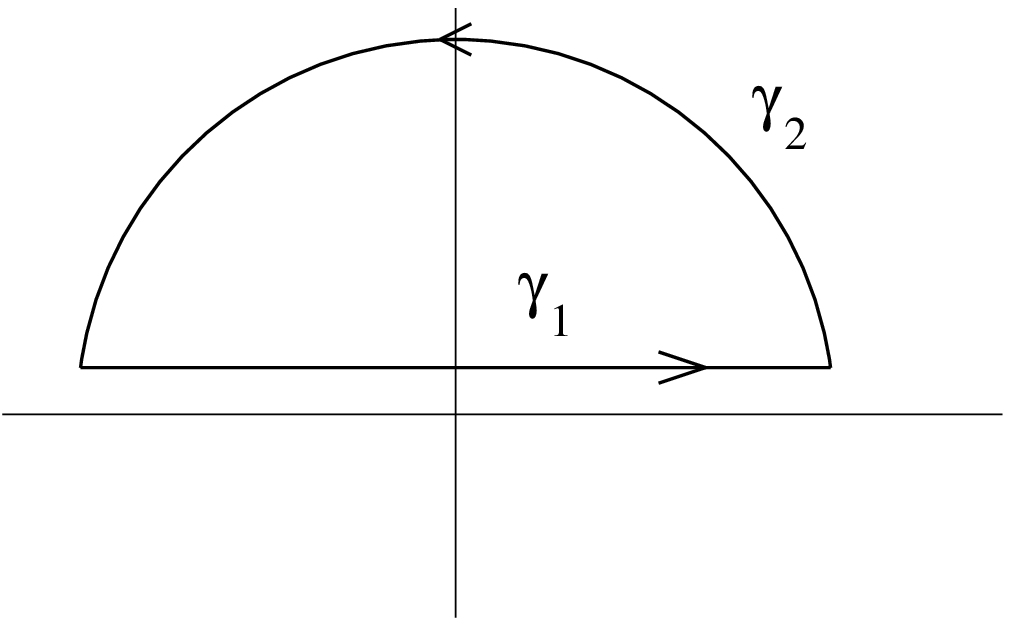}  

\end{tabular}
\caption{}
 \label{perturb p}
\end{figure}
Denote by $C$ (resp. $C'$) the strict transform by the blow up of $\mathbb CP^2$ at $[0:1:0]$ of
a symmetric nonsingular pseudoholomorphic curve of degree $7$ on $\mathbb RP^2$ corresponding to the
quotient curve depicted in Figure \ref{orientations 4 4}b)
(resp. c)). If we denote by $p$ the intersection point of $C$ with the
exceptional section $E$ 
and $F_p$ the fiber of $\Sigma_1$ through $p$, the curve $C$ has
a tangency point of order $2$ with $F_p$ at $p$.

Let us introduce some notations. For $\epsilon>0$, denote by $\gamma_{1,\epsilon}$ and
 $\gamma_{2,\epsilon}$ the following paths
\begin{center}
\begin{tabular}{cccc}
$\gamma_{1,\epsilon}~:$&$[-1,1]$&$\to$&$\mathbb C$\\
   &   $t$&$\mapsto$&$\frac{1}{\epsilon}t+i\epsilon$
\end{tabular} and
\begin{tabular}{cccc}
$\gamma_{2,\epsilon}~:$&$[0,1]$&$\to$&$\mathbb C$\\
  &  $ t$&$\mapsto$&$\frac{1}{\epsilon}e^{i\pi t}+i\epsilon$
\end{tabular}.
\end{center}
Let $\gamma_\epsilon$  be the union of the images of $\gamma_{1,\epsilon} $
and $\gamma_{2,\epsilon}$ (see Figure \ref{perturb p}). Denote also by $\pi$ the projection
$\Sigma_1\to\mathbb CP^1$ on the base $\{y=0\}$, 
$S_C=\pi^{-1}(\mathbb RP^1)$, $D_\epsilon$ the compact region
of $\mathbb C$ bounded by $\gamma_\epsilon$,
$b_{\epsilon}=\pi^{-1}(\gamma_{\epsilon})\cap C$, 
$\hat b_\epsilon$ the closure of the braid $b_{\epsilon}$,
and $N_\epsilon=\pi^{-1}(D_\epsilon)\cap C$. 

\begin{figure}[h]
      \centering
\begin{tabular}{cc}
\includegraphics[height=4.5cm, angle=0]{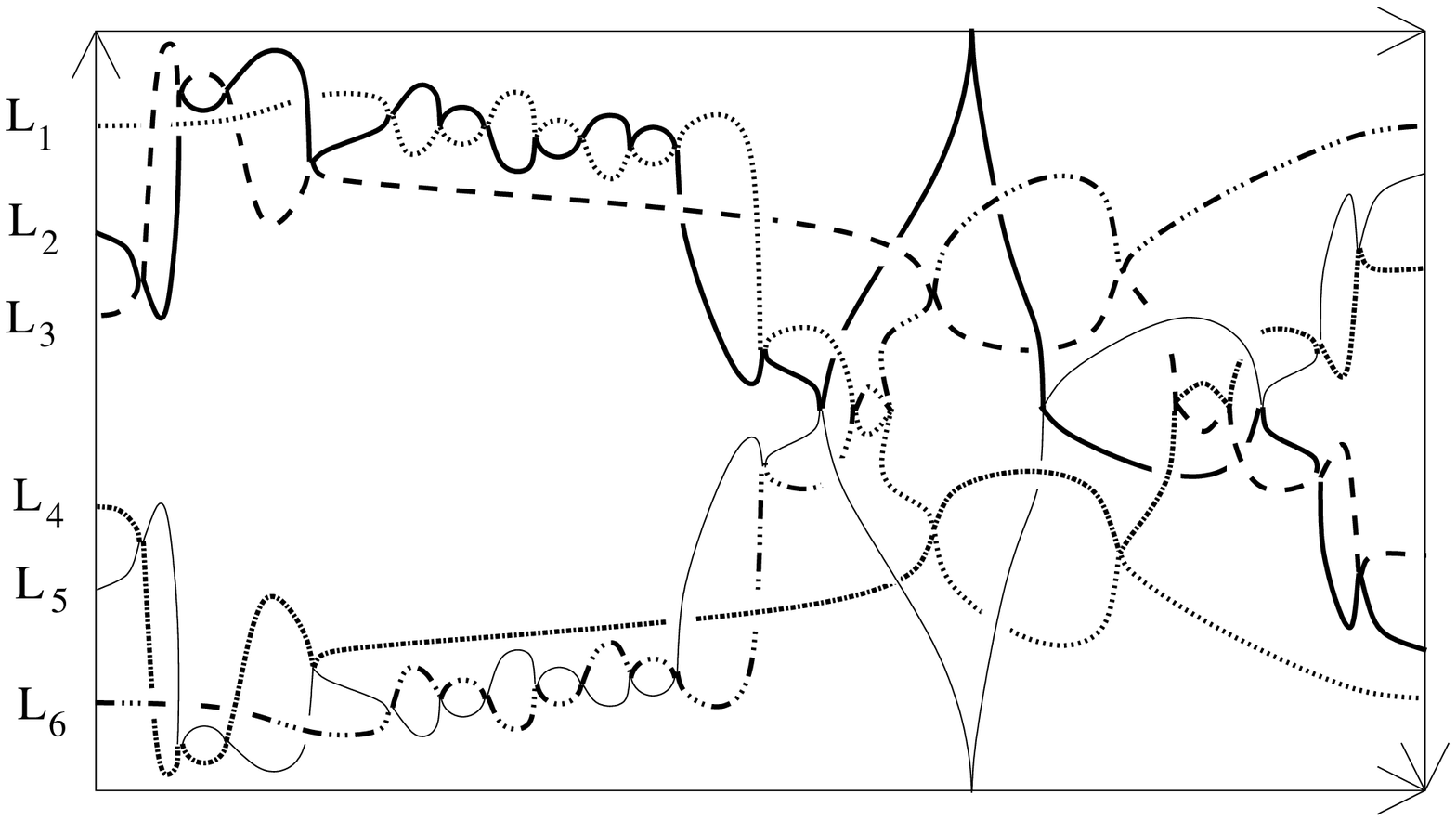}&
     \includegraphics[height=4.5cm, angle=0]{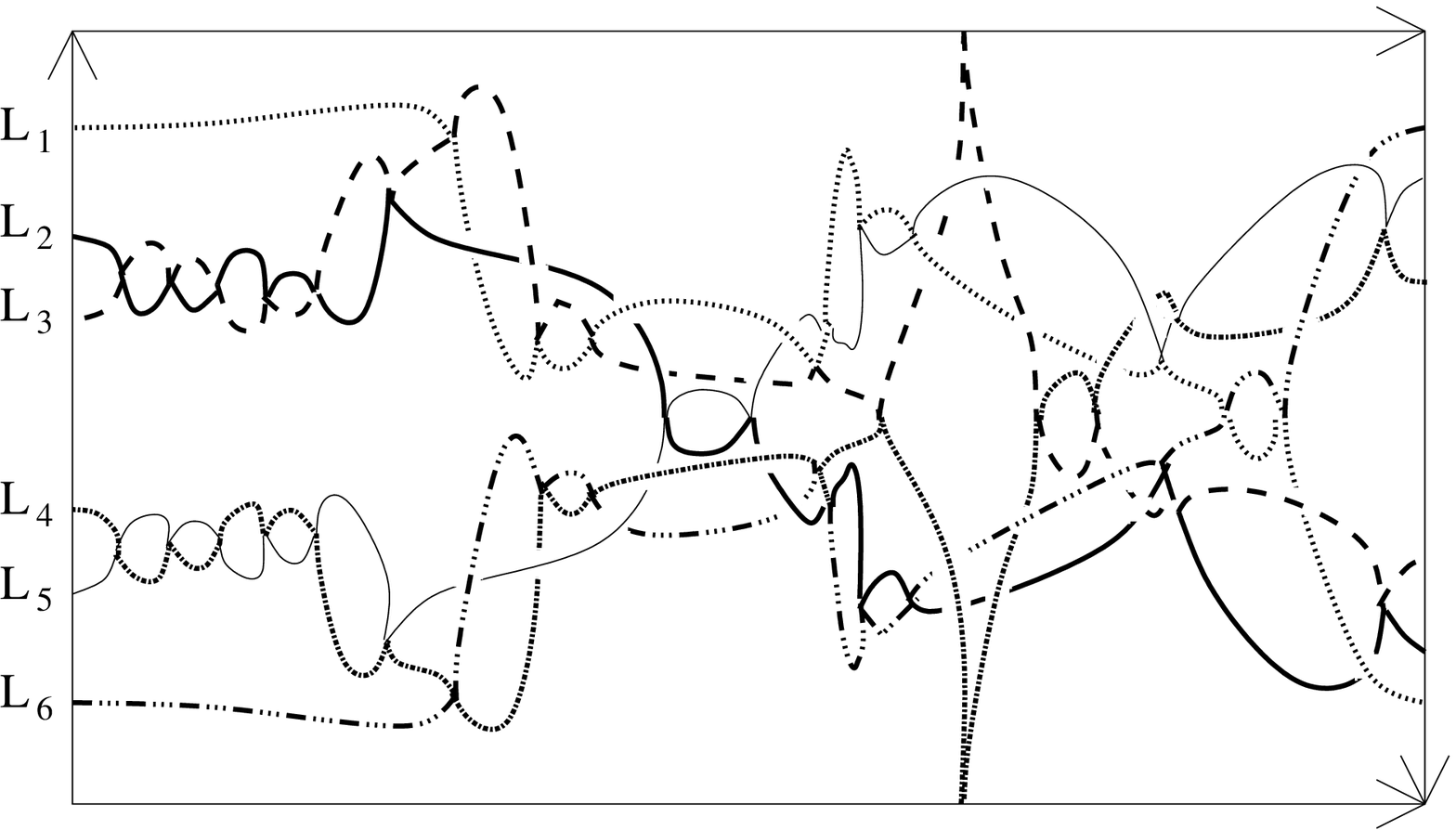}
\\a)&b)
\end{tabular}
\caption{}
 \label{gluing}
\end{figure}

As $C$ is a real curve, $S_C$ is formed by a real part which is $\mathbb RC$, and by a non-real part.
This latter space has several connected components  which are globally invariant by
the complex conjugation. One can  deduce $S_C$ (resp. $S_{C'}$) from the quotient curve of $C$. The curve $S_C$ (resp. $S_{C'}$) is depicted
in Figure \ref{pseudo constr}a) (resp. c)), where the bold lines are used to draw $\mathbb RC$ ant the dashed lines are used to draw
$S_C\setminus\mathbb RC$.
\begin{prop}\label{difference 1}
The real pseudoholomorphic curve $C$ constructed above is a dividing curve.
\end{prop}
\textit{Proof. }
By the Riemann-Hurwitz formula, we have
$$\mu(N_\epsilon)=g(N_\epsilon)+\frac{\mu(\hat b_\epsilon) + 6-e(b_\epsilon)}{2}, $$
where $\mu(N_\epsilon)$ is the number of
connected components of $N_\epsilon$ and $g(N_\epsilon)$ the sum of
the genus of the connected components of $N_\epsilon$. 
We have   $\mu(\hat b_\epsilon)=6$,  $e(b_\epsilon)=0$ and
$\mu(N_\epsilon)\le 6$, so $N_\epsilon$ is 
composed by $6$ disks. Denote these disks
$D_{1,\epsilon},\ldots,D_{6,\epsilon}$ and  their boundaries 
$\partial D_{i,\epsilon}=L_{i,\epsilon}$.
Define also $\overline{D_{i,\epsilon}}=\textrm{conj}(D_{i,\epsilon})$. As
$\epsilon\to 0$, these $12$ disks glue together along $S_C$ as
depicted in Figure \ref{gluing}a), and $C$ is the
result of this gluing. Moreover, $C\setminus\mathbb RC$ is the result of the gluing of these $12$ disks along
$S_C\setminus\mathbb RC$. Hence, to find the type of $C$, we just have to study how the $12$ disks
glue along $S_C\setminus\mathbb RC$.

Denote by $D_{i,\epsilon}|| D_{j,\epsilon}$ the relation
``$D_{i,\epsilon}$ glues  with $D_{j,\epsilon}$ along a connected component of $S_C\setminus\mathbb RC$
as $\epsilon\to 0$''.
Using the fact that each connected component of $S_C\setminus\mathbb RC$ is globally invariant by the complex
conjugation, we have  (see Figure \ref{gluing}a)) :
\begin{center}
\begin{tabular}{lllll}
 $D_{1,\epsilon}|| \overline{D_{2,\epsilon}}$,&
$D_{1,\epsilon}|| \overline{D_{6,\epsilon}}$,&
$D_{1,\epsilon}|| \overline{D_{4,\epsilon}}$,&
$D_{2,\epsilon}|| \overline{D_{3,\epsilon}}$,&
$D_{2,\epsilon}|| \overline{D_{5,\epsilon}}$,\\
$D_{3,\epsilon}|| \overline{D_{6,\epsilon}}$,&
$D_{3,\epsilon}|| \overline{D_{4,\epsilon}}$,&
$D_{4,\epsilon}|| \overline{D_{5,\epsilon}}$,&
$D_{5,\epsilon}|| \overline{D_{6,\epsilon}}$,&
$D_{i,\epsilon}|| \overline{D_{j,\epsilon}}\implies D_{j,\epsilon}|| \overline{D_{i,\epsilon}}.$
\end{tabular}
\end{center}
The curve $C$ is a dividing curve
if and only if there exist two equivalence classes for $||$.
Here the equivalence classes are
$\{D_{1,\epsilon},D_{3,\epsilon},D_{5,\epsilon},\overline{D_{2,\epsilon}},
\overline{D_{4,\epsilon}},\overline{D_{6,\epsilon}} \} $
and
$\{D_{2,\epsilon},D_{4,\epsilon},D_{6,\epsilon},\overline{D_{1,\epsilon}},
\overline{D_{3,\epsilon}},\overline{D_{5,\epsilon}} \} ,$
Hence, $C$ is a dividing curve.~\findemo
\begin{prop}\label{difference 2}
The real pseudoholomorphic curve $C'$ constructed above is a dividing curve.
\end{prop}
\textit{Proof. }We keep the same notations than in Proposition \ref{difference 1}. As in this proposition, the closure of the braid $b_\epsilon$ has $6$ components, the surface $N_\epsilon$ is
composed by $6$ disks and
the two  equivalence classes for the relation $||$ are
$\{D_{1,\epsilon},D_{2,\epsilon},D_{4,\epsilon},\overline{D_{3,\epsilon}},
\overline{D_{5,\epsilon}},\overline{D_{6,\epsilon}} \} $
and
$\{D_{3,\epsilon}, D_{5,\epsilon}, D_{6,\epsilon}, \overline{D_{1,\epsilon}},
\overline{D_{2,\epsilon}}, \overline{D_{4,\epsilon}} \}$  (see Figures
\ref{pseudo constr}b) and \ref{gluing}b)).
Hence, $C'$ is a dividing curve.
~\findemo

\begin{cor}\label{pseudo non alg 2}
The complex schemes $\langle J\amalg 4_+\amalg 4_-\amalg  1_+\langle 2_+\amalg 2_-\rangle\rangle_I$ and $\langle J\amalg 2_+\amalg 2_-\amalg  1_+\langle 4_+\amalg 4_-\rangle\rangle_I$ are realizable by
 nonsingular symmetric real pseudoholomorphic
 curves of degree $7$ on $\mathbb RP^2$.
\end{cor}

\section{Algebraic statements}\label{algebraic state}

\subsection{Prohibitions}

We prove in this section the algebraic prohibitions stated in Theorem
\ref{Main thm4}. The main tools are the real trigonal graphs (see
section \ref{comb theo}) and the cubic resolvent of an algebraic curve of bidegree
$(4,0)$ in $\Sigma_n$ (see \cite{OS2}).

\begin{figure}[h]
      \centering
 \begin{tabular}{cccc}
\includegraphics[height=2cm, angle=0]{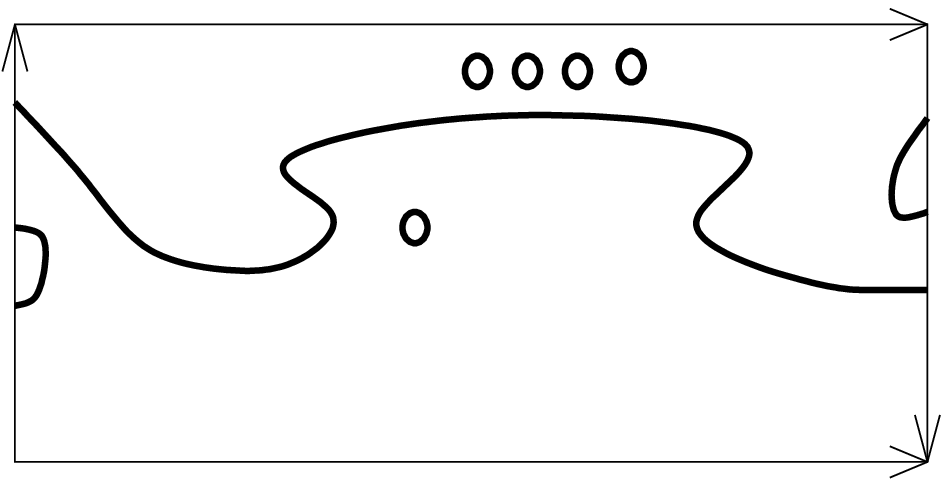}    &
\includegraphics[height=2cm, angle=0]{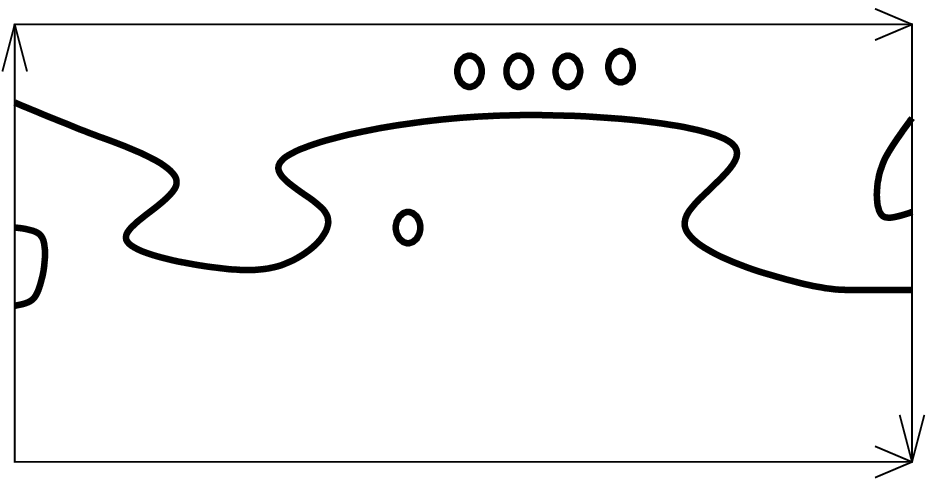}    &
 \includegraphics[height=2cm, angle=0]{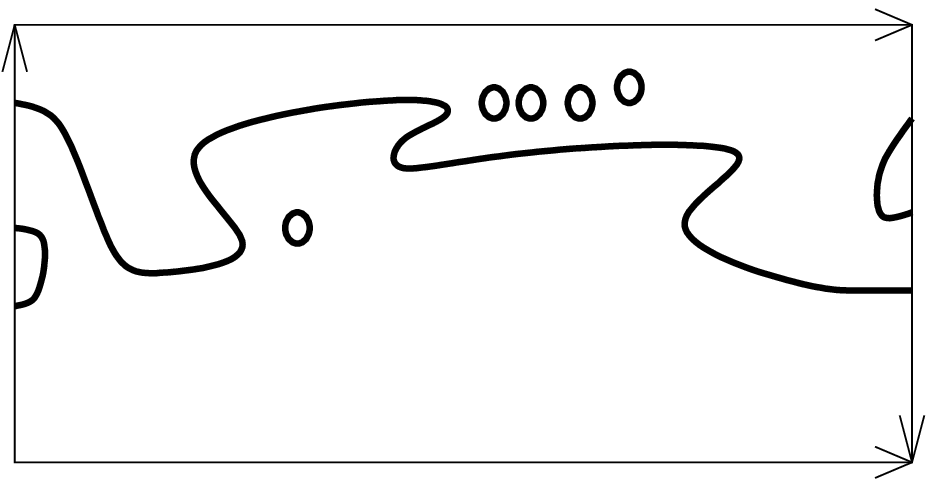} &
 \includegraphics[height=2cm, angle=0]{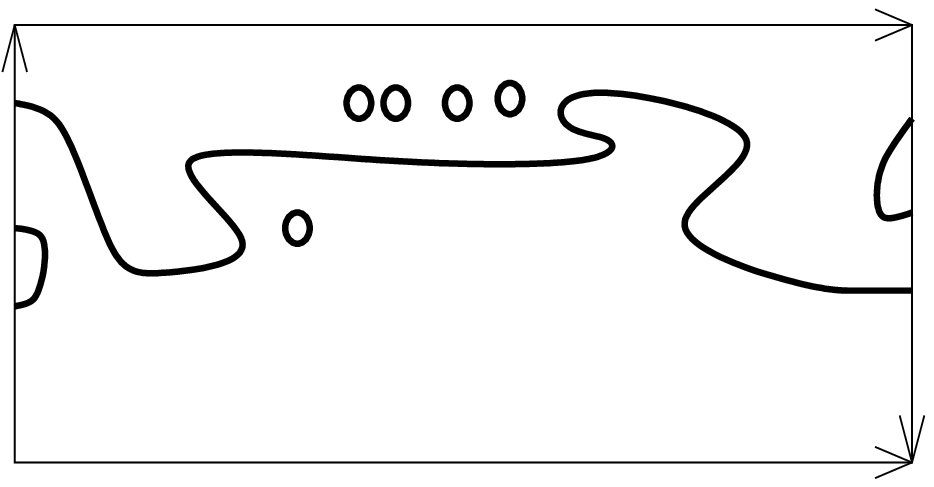}
\\ a)&b) &c)&d)\\
 \end{tabular}
\caption{}
 \label{8 4 I alg}
\end{figure}
\begin{lemma}\label{lem I 1}
The $\mathcal L$-schemes depicted in Figures \ref{8 4 I alg}b), c) and
d) are not realizable by   trigonal nonsingular real pseudoholomorphic curves on
$\Sigma_3$.
\end{lemma}
\textit{Proof. }Compute the braids associated to these $\mathcal L$-schemes~:
\begin{center}
\begin{tabular}{ll} $b^{14}=\sigma_1^{-1}\sigma_2^{-1}\sigma_1^{-2}\sigma_2^{-1}\sigma_1\sigma_2^{-4}\sigma_1^{-1}\Delta_3^3$,
& $b^{15}=\sigma_1^{-1}\sigma_2^{-1}\sigma_1^{-1}\sigma_2^{-5}\sigma_1^{-1}\Delta_3^3$,
\\$b^{16}=\sigma_1^{-1}\sigma_2^{-1}\sigma_1^{-1}\sigma_2^{-1}\sigma_1\sigma_2^{-4}\sigma_1^{-1}\sigma_2^{-1}\Delta_3^3$.
\end{tabular}
\end{center}
These braids verify $e(b)=0$, so they are quasipositive if and only if they are trivial. Computing
their Garside normal form (see \cite{Gar} or \cite{Jaq}), we find
\begin{center}
\begin{tabular}{lll}
$b^{14}=\sigma_2^{3}\sigma_1^{2}\sigma_2^{2}\sigma_1^{2}\Delta_3^{-3}$,
&$b^{15}=\sigma_1\sigma_2^{2}\sigma_1^2\sigma_2^{2}\sigma_1^{2}\Delta_3^{-3}$,
&$b^{16}=\sigma_1\sigma_2^{3}\sigma_1^{2}\sigma_2^{2}\sigma_1\Delta_3^{-3}$.
\end{tabular}
\end{center}
Thus, no one of these braids is quasipositive.
~\findemo

\begin{lemma}\label{lem I 3}
The $\mathcal L$-scheme depicted in Figure \ref{8 4 I alg}a) is not realizable by
trigonal nonsingular real algebraic curves on
$\Sigma_3$.
\end{lemma}
\textit{Proof. }The weighted comb associated to this $\mathcal L$-scheme is
$$w_1=(g_3g_6g_1g_4g_1g_6g_5g_2g_3g_6g_1g_4g_1g_6(g_3g_2)^3g_3g_6g_1g_4g_1g_6g_5g_2,1,2,0).$$ We have $\mu(w_1)=0$ so, according to Proposition \ref{prohib comb} the lemma is proved.
\begin{prop}\label{prop 84}
The real scheme $\langle J\amalg 8 \amalg 1\langle 4 \rangle\rangle$ is not realizable by nonsingular symmetric dividing real algebraic curves of degree $7$ on $\mathbb RP^2$.
\end{prop}
\textit{Proof. }Suppose that there exists a dividing symmetric curve
which contradicts Proposition \ref{prop 84}. Denote by $X$ its quotient curve.
Blow up $\Sigma_2$ at the intersection point of $X$ and $E$ and then
blow down  the strict transform of the 
fiber. The obtained surface is $\Sigma_3$ and  the strict transform of
$X$ is 
a trigonal curve which has a double
point with non-real tangents.
Smooth this double point in order to obtain an oval. Then,
according to Propositions \ref{only possibility} and \ref{alg
  operations}, one should obtain one of the 
 $\mathcal L$-schemes   depicted in Figures
\ref{8 4 I alg}a), b), c) and d). However, according to Lemmas \ref{lem I 1} and
\ref{lem I 3}, these $\mathcal L$-schemes are not algebraically realizable so there is a contradiction.~\findemo

\begin{figure}[h]
      \centering
\begin{tabular}{cc}
\includegraphics[height=1.8cm, angle=0]{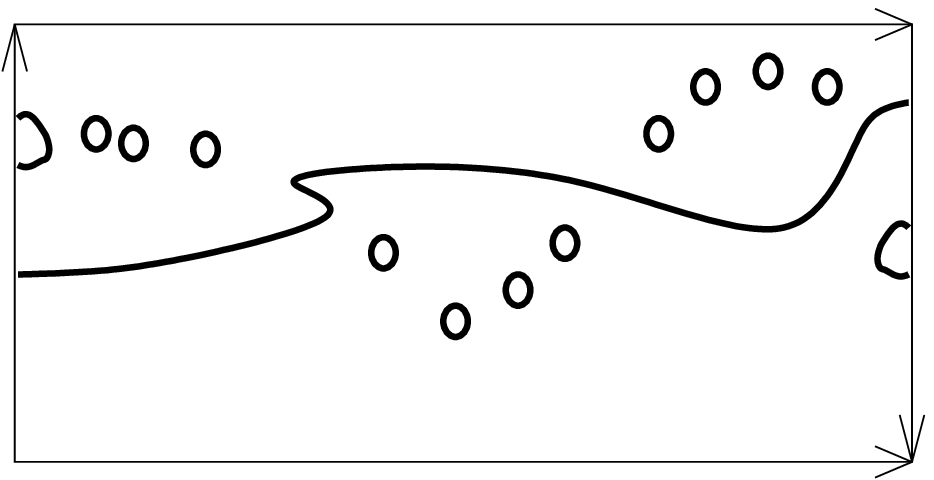}    &
\includegraphics[height=1.8cm, angle=0]{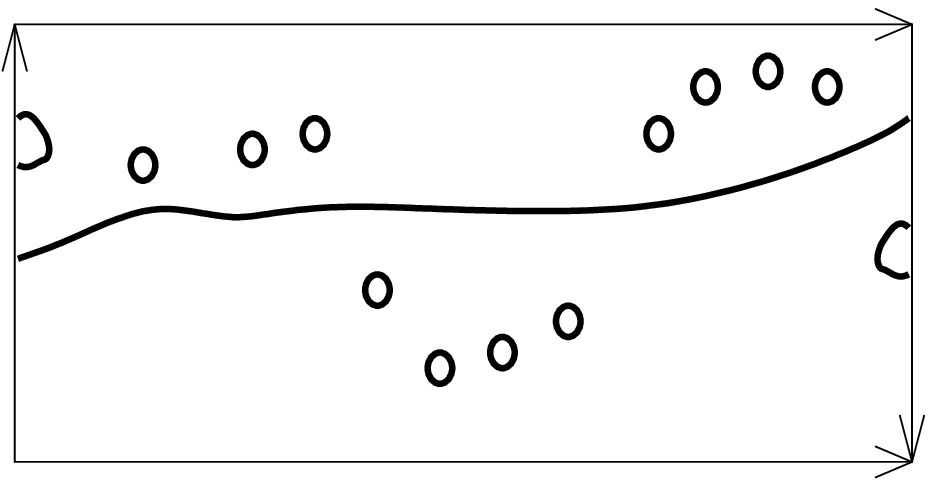}
\\ a)&b)
 \end{tabular}
\caption{}
 \label{4 8 I alg}
\end{figure}

\begin{lemma}\label{lem 48 1}
The $\mathcal L$-scheme depicted in Figures \ref{4 8 I alg}a)is not realizable by  trigonal nonsingular real pseudoholomorphic curves on
$\Sigma_5$.
\end{lemma}
\textit{Proof. }
The braid  associated to this $\mathcal L$-schemes and its determinant
are~:
\begin{center}
\begin{tabular}{lll}

$b^{17}=\sigma_2^{-4}\sigma_1^{-5}\sigma_2^{-1}\sigma_1\sigma_2^{-4}\sigma_1^{-1}\sigma_2\Delta_3^5$&
  and & 
$301.$
\end{tabular}

\end{center}
We have $e(b^{17})=2$, so according to Proposition \ref{alex}, this
braid is not quasipositive.

\begin {lemma}\label{lem 48 3}
The $\mathcal L$-scheme depicted in Figures \ref{4 8 I alg}b) is not realizable by
  trigonal nonsingular real algebraic curves on
$\Sigma_5$.
\end{lemma}
\textit{Proof. }The weighted comb associated to this $\mathcal
L$-scheme is
\begin{center}
$w_2=((a^3g_3g_6g_1g_4g_1g_6)^3,3,6,2)$
where $a=g_3g_2$.
\end{center}
We have $\mu(w_2)=0$, so according to Proposition \ref{prohib comb},
  this $\mathcal L$-scheme is not realizable by
  trigonal nonsingular real algebraic curves on
$\Sigma_5$.
\findemo

\begin {prop}\label{prop 48}
The real scheme $\langle J\amalg 4 \amalg 1\langle 8 \rangle\rangle$ is not realizable by nonsingular symmetric dividing real algebraic curves of degree $7$ on $\mathbb RP^2$.
\end{prop}
\textit{Proof. }Suppose that there exists a dividing symmetric curve
which contradicts the lemma. Denote by $X$ its quotient curve. 
Blow up $\Sigma_2$ at the intersection point of $X$ and $E$ and then
blow down  the strict transform of the 
fiber. The obtained surface is $\Sigma_3$ and  the strict transform of
$X$, still denoted by $X$,  is 
a trigonal curve which has a double
point with non-real tangents on the base $\{y=0\}$. Let $\widetilde X$ be
the cubic resolvent (see \cite{OS2}) of the union of the base $\{y=0\}$ and
$X$. The curve $\widetilde X$ is a trigonal real algebraic curve  on
$\Sigma_6$ with a triple point coming from the triple point of $X\cup\{y=0\}$.
Blow up $\Sigma_6$ at this triple point  and then  blow down  the
strict transform of the 
fiber. The obtained surface is $\Sigma_5$ and  the strict transform of
$\widetilde X$ is 
a trigonal curve with ordinary double points. Smooth all the double
points in order to obtain ovals. Then,
according to Propositions \ref{only possibility} and \ref{alg
  operations}, one should obtain one of the 
 $\mathcal L$-schemes   depicted in Figure \ref{4 8 I alg}. It
has been proved in Lemmas \ref{lem 48 1} and \ref{lem 48 3} that these $\mathcal L$-schemes
 are not algebraically realizable, so there is a contradiction.~\findemo

\subsection{Perturbation of a reducible symmetric curve}\label{viro}
The standard method to construct a lot of different isotopy types of nonsingular algebraic curves is to perturb a singular curve in many ways.
 So the first idea to construct symmetric algebraic curves is to perturb in many symmetric ways a singular symmetric
algebraic curve.

To perturb real algebraic singular curves, we use the Viro method. The
unfamiliar reader can refer to \cite{V1}, \cite{V4},   \cite{V2}, \cite{Ris} and \cite{IS2}.

\begin{prop}\label{main constr}
All the real schemes listed in table \ref{list scheme} are realizable by nonsingular symmetric real algebraic curves of degree $7$ on $\mathbb RP^2$. Moreover, those marked with a $^*$ are realized
by a dividing symmetric curve and those marked with a $^\circ$ are realized by a non-dividing curve.
\begin{table}[h]
\centering
\begin{tabular}{lllll}
$\langle J\rangle^\circ$ &$\langle J\amalg 10 \amalg 1\langle 1 \rangle\rangle^\circ$ &
 $\langle J\amalg 8 \amalg 1\langle 3 \rangle\rangle^\circ$ &$\langle  J \amalg 1\langle 6 \rangle\rangle^\circ$ &
 $\langle  J \amalg 1\langle 9 \rangle\rangle^\circ$

\\

$\langle J\amalg 1\rangle^\circ$ &$\langle J\amalg 11 \amalg 1\langle 1 \rangle\rangle$$^{\circ,*}$&
  $\langle J\amalg 9 \amalg 1\langle 3 \rangle\rangle$$^{\circ,*}$ &  $\langle J\amalg 1 \amalg 1\langle 6 \rangle\rangle^\circ$ &
  $\langle J\amalg 1 \amalg 1\langle 9 \rangle\rangle$$^{\circ,*}$

\\

$\langle J\amalg 2\rangle^\circ$ &$\langle J\amalg 12 \amalg 1\langle 1 \rangle\rangle^\circ$&
  $\langle J\amalg 10 \amalg 1\langle 3 \rangle\rangle^\circ$ &$\langle J\amalg 2 \amalg 1\langle 6 \rangle\rangle^\circ$  &
    $\langle J\amalg 2 \amalg 1\langle 9 \rangle\rangle^\circ$

\\

$\langle J\amalg 3\rangle^\circ$ &$\langle J\amalg 13 \amalg 1\langle 1 \rangle\rangle$$^{*}$&
 $\langle J\amalg 11 \amalg 1\langle 3 \rangle\rangle$$^{*}$  & $\langle J\amalg 3 \amalg 1\langle 6 \rangle\rangle^\circ$ &
    $\langle J\amalg 3 \amalg 1\langle 9 \rangle\rangle$$^{\circ,*}$

\\

$\langle J\amalg 4\rangle^\circ$ & &
  &$\langle J\amalg 4 \amalg 1\langle 6 \rangle\rangle$$^{\circ,*}$  &

\\

$\langle J\amalg 5\rangle^\circ$ & $\langle  J \amalg 1\langle 2 \rangle\rangle^\circ$  &
 $\langle J \amalg 1\langle 4\rangle\rangle^\circ$ &  &

\\

$\langle J\amalg 6\rangle^\circ$ &  $\langle J\amalg 1 \amalg 1\langle 2 \rangle\rangle^\circ$  &
 $\langle J\amalg 1 \amalg 1\langle 4 \rangle\rangle^\circ$ &$\langle J\amalg 6 \amalg 1\langle 6 \rangle\rangle^\circ$  &

\\

$\langle J\amalg 7\rangle$$^{\circ,*}$ &  $\langle J\amalg 2 \amalg 1\langle 2 \rangle\rangle$$^{\circ,*}$ &
 $\langle J\amalg 2 \amalg 1\langle 4 \rangle\rangle$$^{\circ,*}$   &  &
  $\langle  J \amalg 1\langle 10 \rangle\rangle$$^{\circ,*}$

\\

$\langle J\amalg 8\rangle^\circ$ & $\langle J\amalg 3 \amalg 1\langle 2 \rangle\rangle^\circ$  &
 $\langle J\amalg 3 \amalg 1\langle 4 \rangle\rangle^\circ$  &  &
  $\langle J\amalg 1 \amalg 1\langle 10 \rangle\rangle^\circ$

\\

$\langle J\amalg 9\rangle$$^{\circ,*}$ & $\langle J\amalg 4 \amalg 1\langle 2 \rangle\rangle$$^{\circ,*}$  &
 $\langle J\amalg 4 \amalg 1\langle 4 \rangle\rangle^\circ$   &  &
   $\langle J\amalg 2 \amalg 1\langle 10 \rangle\rangle^\circ$

\\

$\langle J\amalg 10\rangle^\circ$ & $\langle J\amalg 5 \amalg 1\langle 2 \rangle\rangle^\circ$  &
 $\langle J\amalg 5 \amalg 1\langle 4 \rangle\rangle^\circ$   & $\langle J \amalg 1\langle 7 \rangle\rangle^\circ$  &
  $\langle J\amalg 3 \amalg 1\langle 10 \rangle\rangle^\circ$

\\

$\langle J\amalg 11\rangle$$^{\circ,*}$ & $\langle J\amalg 6 \amalg 1\langle 2 \rangle\rangle^\circ$  &
 $\langle J\amalg 6 \amalg 1\langle 4 \rangle\rangle$$^{\circ,*}$  & $\langle J\amalg 1 \amalg 1\langle 7 \rangle\rangle$$^{\circ,*}$  &
  $\langle J\amalg 4 \amalg 1\langle 10 \rangle\rangle$$^{*}$

\\

$\langle J\amalg 12\rangle^\circ$ & $\langle J\amalg 7 \amalg 1\langle 2 \rangle\rangle^\circ$  &
  &$\langle J\amalg 2 \amalg 1\langle 7 \rangle\rangle^\circ$   &

\\

$\langle J\amalg 13\rangle$$^{\circ,*}$ & $\langle J\amalg 8 \amalg 1\langle 2 \rangle\rangle$$^{\circ,*}$  &
 $\langle J\amalg 8 \amalg 1\langle 4 \rangle\rangle^\circ$  &$\langle J\amalg 3 \amalg 1\langle 7 \rangle\rangle$$^{\circ,*}$   &
  $\langle  J \amalg 1\langle 11 \rangle\rangle^\circ$

\\

$\langle J\amalg 14\rangle^\circ$ & $\langle J\amalg 9 \amalg 1\langle 2 \rangle\rangle^\circ$  &
  &  &
 $\langle J\amalg 1 \amalg 1\langle 11 \rangle\rangle$$^{\circ,*}$

\\

$\langle J\amalg 15\rangle$$^{*}$ & $\langle J\amalg 10 \amalg 1\langle 2 \rangle\rangle^\circ$  &
 $\langle J\amalg 10 \amalg 1\langle 4 \rangle\rangle$$^{*}$  &$\langle J\amalg 5 \amalg 1\langle 7 \rangle\rangle$$^{*}$   &
  $\langle J\amalg 2 \amalg 1\langle 11 \rangle\rangle^\circ$

\\

 & $\langle J\amalg 11 \amalg 1\langle 2 \rangle\rangle^\circ$ &
   &  &
  $\langle J\amalg 3 \amalg 1\langle 11 \rangle\rangle$$^{*}$

\\

$\langle J \amalg 1\langle 1 \rangle\rangle^\circ$ & $\langle J\amalg 12 \amalg 1\langle 2 \rangle\rangle$$^{*}$ &
 $\langle J \amalg 1\langle 5 \rangle\rangle^\circ$  &  &

\\

$\langle J\amalg 1 \amalg 1\langle 1 \rangle\rangle^\circ$ &   &
 $\langle J\amalg 1 \amalg 1\langle 5 \rangle\rangle^\circ$ &  &
    $\langle  J \amalg 1\langle 12 \rangle\rangle^\circ$

\\

$\langle J\amalg 2 \amalg 1\langle 1 \rangle\rangle^\circ$ &  $\langle J \amalg 1\langle 3 \rangle\rangle^\circ$ &
  $\langle J\amalg 2 \amalg 1\langle 5 \rangle\rangle^\circ$ & $\langle J \amalg 1\langle 8 \rangle\rangle^\circ$ &

\\

$\langle J\amalg 3 \amalg 1\langle 1 \rangle\rangle$$^{\circ,*}$& $\langle J\amalg 1 \amalg 1\langle 3 \rangle\rangle^\circ$  &
  $\langle J\amalg 3 \amalg 1\langle 5 \rangle\rangle$$^{\circ,*}$   &  $\langle J\amalg 1 \amalg 1\langle 8 \rangle\rangle^\circ$ &
 $\langle J\amalg 2 \amalg 1\langle 12 \rangle\rangle$$^{*}$

\\

$\langle J\amalg 4 \amalg 1\langle 1 \rangle\rangle^\circ$ & $\langle J\amalg 2 \amalg 1\langle 3 \rangle\rangle^\circ$  &
   $\langle J\amalg 4 \amalg 1\langle 5 \rangle\rangle^\circ$  &$\langle J\amalg 2 \amalg 1\langle 8 \rangle\rangle$$^{\circ,*}$  &

\\

$\langle J\amalg 5 \amalg 1\langle 1 \rangle\rangle^\circ$ & $\langle J\amalg 3 \amalg 1\langle 3 \rangle\rangle$$^{\circ,*}$  &
  $\langle J\amalg 5 \amalg 1\langle 5 \rangle\rangle$$^{\circ,*}$   &$\langle J\amalg 3 \amalg 1\langle 8 \rangle\rangle^\circ$  &
  $\langle  J \amalg 1\langle 13 \rangle\rangle^\circ$

\\

$\langle J\amalg 6 \amalg 1\langle 1 \rangle\rangle^\circ$ &  $\langle J\amalg 4 \amalg 1\langle 3 \rangle\rangle^\circ$  &
  $\langle J\amalg 6 \amalg 1\langle 5 \rangle\rangle^\circ$   &$\langle J\amalg 4 \amalg 1\langle 8 \rangle\rangle^\circ$  &
  $\langle J\amalg 1 \amalg 1\langle 13 \rangle\rangle$$^{*}$

\\

$\langle J\amalg 7 \amalg 1\langle 1 \rangle\rangle$$^{\circ,*}$ &  $\langle J\amalg 5 \amalg 1\langle 3 \rangle\rangle$$^{\circ,*}$ &
   $\langle J\amalg 7 \amalg 1\langle 5 \rangle\rangle$$^{\circ,*}$  &  &

\\
 $\langle J\amalg 8 \amalg 1\langle 1 \rangle\rangle^\circ$  &  $\langle J\amalg 6 \amalg 1\langle 3 \rangle\rangle^\circ$  &
   $\langle J\amalg 8 \amalg 1\langle 5 \rangle\rangle^\circ$   &   &
 $\langle J\amalg 1\amalg 1\langle 1 \langle 1\rangle\rangle\rangle$$^{*}$

\\
 $\langle J\amalg 9 \amalg 1\langle 1 \rangle\rangle$$^{\circ,*}$  &   $\langle J\amalg 7 \amalg 1\langle 3 \rangle\rangle$$^{\circ,*}$ &
  $\langle J\amalg 9 \amalg 1\langle 5 \rangle\rangle$$^{*}$    &   
\end{tabular}
\caption{}
\label{list scheme}
\end{table}
\end{prop}
\textit{Proof. }In order to apply the Viro method without change of coordinates, we consider here symmetry with respect to the line $\{y-z=0\}$. Consider the union of the line $\{x=0\}$ and three symmetric conics on $\mathbb RP^2$ tangent
to each other in the two symmetric points $[0:0:1]$ and $[0:1:0]$. Using the
 Viro method and the classification, up to isotopy, of the
curves of degree $7$ on $\mathbb RP^2$ with the only singular point $Z_{15}$ established by A. B. Korchagin in \cite{K1}, we perturb these reducible symmetric curves.
In order to obtain nonsingular symmetric curves, we have to perturb symmetrically the two singular points. That is to say, if we
perturb the singular point at $[0:0:1]$ gluing the chart of a polynomial $P(x,y)$, we have to perturb the singular point at
$[0:1:0]$ gluing the chart of the polynomial $y^7P(\frac{x}{y},\frac{1}{y})$.
~\findemo

\textbf{Remark. }Using this method, we constructed nonsingular symmetric algebraic curves of degree $7$ on $\mathbb RP^2$ realizing
the complex schemes $\langle J\amalg  4_+\amalg 5_-\amalg 1_+\langle  1_+\rangle\rangle_I$ and
 $\langle J\amalg 3_+\amalg 6_-\amalg 1_-\langle 1_+\rangle\rangle_I$. So unlike in the $M$-curves case, the real scheme of
a nonsingular symmetric curves of degree $7$ on $\mathbb RP^2$ does not determine its complex scheme.

\begin{prop}\label{restr thm 4 5}
The complex schemes 
$$\langle J\amalg  1_-\langle  1_+\amalg 3_-\rangle\rangle_I\textrm{,
 }\langle J\amalg 2_+\amalg 3_-\amalg 1_-\langle
 1_-\rangle\rangle_I\textrm{ and }\langle J\amalg 1_-\amalg 1_-\langle 2_+\amalg 3_-\rangle\rangle_I$$
 are realizable by nonsingular symmetric real algebraic dividing curves of degree $7$ on $\mathbb RP^2$.
\end{prop}
\textit{Proof. }In  \cite{It}, symmetric sextics realizing the complex schemes
 $\langle 1_-\langle  1_+\amalg 3_-\rangle\rangle_I$, $\langle
 5\amalg 1_-\langle 1_-\rangle\rangle_I$ and 
 $\langle 1 \amalg 1_-\langle 2_+\amalg 3_-\rangle\rangle_I$
  are constructed. Consider the union of each of these
curves and a real line oriented and disposed on $\mathbb RP^2$ such that the (symmetric) perturbations according to the orientations satisfies the Rokhlin-Mischachev orientation formula. So, according
to Theorem $4.8$ in \cite{V3}, the obtained real algebraic symmetric curves of degree $7$ on $\mathbb RP^2$ are of type $I$ and
realize the announced complex schemes.
~\findemo

\subsection{Parametrization of a rational curve}\label{param constr}
Here we apply the method used in \cite{O6} and \cite{O7}.  Namely, we construct a singular rational curve and perturb it using Shustin's results on the independent perturbations of generalized semi-quasihomogeneous singular points of a curve keeping the same Newton polygon
 (see \cite{Sh1} and \cite{Sh2}).

\begin{figure}[h]
      \centering
 \begin{tabular}{ccc}
\includegraphics[height=4cm, angle=0]{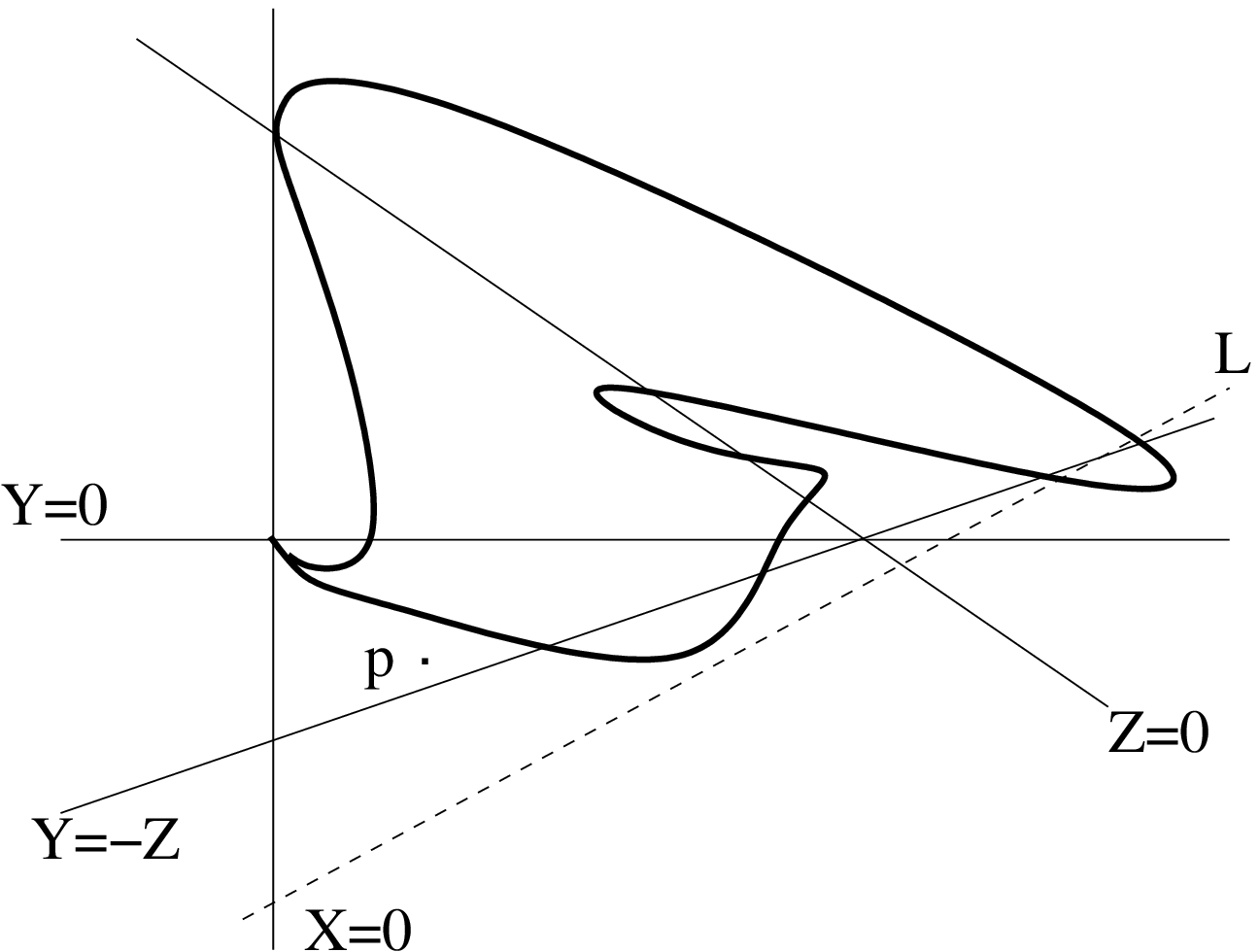}    &
    \includegraphics[height=2cm, angle=0]{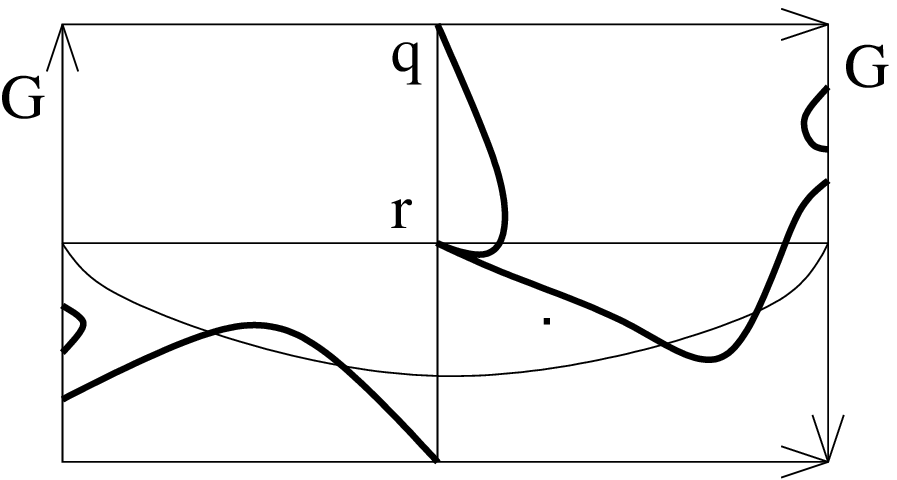}&
\includegraphics[height=2cm, angle=0]{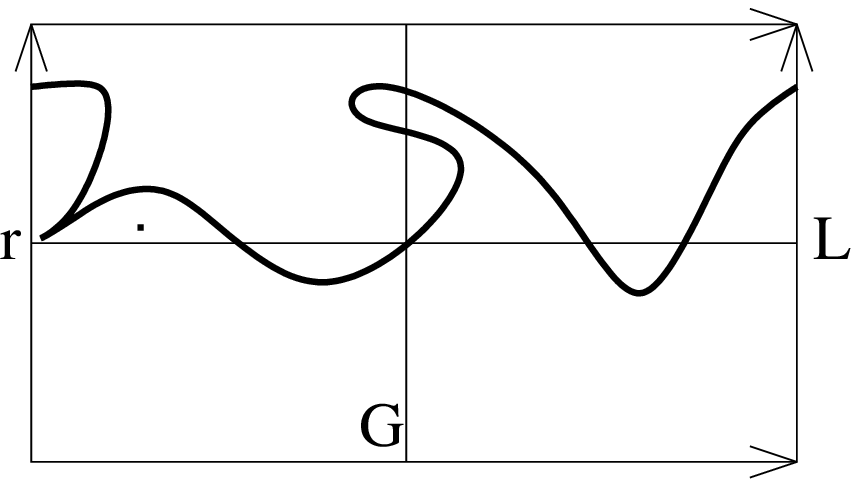}
\\ a)&b)&c)
 \end{tabular}
\caption{}
 \label{rational}
\end{figure}

\begin{prop}\label{par}
There exists a rational real algebraic curve of degree $4$ on $\mathbb RP^2$ situated with respect to the lines
$\{x=0\}$, $\{y=0\}$, $\{z=0\}$ and $\{y=-z\}$ as shown in Figure \ref{rational}a), with a singular point of type $A_4$ at $[0:0:1]$, a singular point of type $A_1$ at $p$
and a tangency point of order $2$ with the line $\{x=0\}$ at $[0:1:0]$.
\end{prop}
\textit{Proof. }Consider the  map from $\mathbb C$ to
$\mathbb CP^2$ given by $t\mapsto [x(t):y(t):z(t)]$ where
\begin{center}
\begin{tabular}{lcl}
$\left\{ \begin{array}{l}
x(t)=t^2
\\y(t)=\alpha t^2(t-\gamma)(t+\gamma)
\\z(t)=-((t-1)(\epsilon t-1)(\delta t-1))
\end{array}\right .$
& &with $\alpha =\frac{1}{11}$, $\gamma =\delta =\frac{9}{10}$ and $\epsilon =\frac{99}{100}$.
\end{tabular}
\end{center}
The curve defined by this map has a singular point of type $A_4$ at $[0:0:1]$, as we can see using the following identity~:
$$y(t)z(t)+\alpha\gamma ^2x(t)z(t) - \alpha x(t)^2=- {\displaystyle \frac {81}{1000}} \,t^{7} + {\displaystyle
\frac {2781}{11000}} \,t^{6} - {\displaystyle \frac {289}{1100}}
\,t^{5}. $$
Moreover, it is clear from the equations that the curve has a tangency point of order $2$ with the line $\{x=0\}$ at $[0:1:0]$.
This map define an algebraic curve of degree $4$ on $\mathbb CP^2$,
$C=\textrm{Res}_t(x(t)Y-y(t)X,x(t)Z-z(t)X)/X^3$.
Considering $C$ on the affine plane $\{Z=1\}$, we obtain
\begin{center}
\begin{tabular}{ll}
$C=$&${\displaystyle -\frac{1}{11}Y^2-\frac{13468021579}{1331000000000000}X^4+\frac{1666467523}{1210000000000}YX^3+\frac{275261}{3025000}X^2Y-\frac{65485017}{1100000000}X^2Y^2}$\\
\\&${\displaystyle +\frac{2781}{5500}XY^2-\frac{81}{6050}XY+\frac{793881}{1000000}XY^3+\frac{26346141}{6655000000}X^3-\frac{6561}{13310000}X^2}.$
\end{tabular}
\end{center}
Projecting this curve to the line $\{y=0\}$ and using the Budan Fourier theorem (see \cite{BPR}), one can check that $C$ verifies the hypothesis of the proposition.
~\findemo

\begin{figure}[h]
      \centering
 \begin{tabular}{cc}
\includegraphics[height=2cm, angle=0]{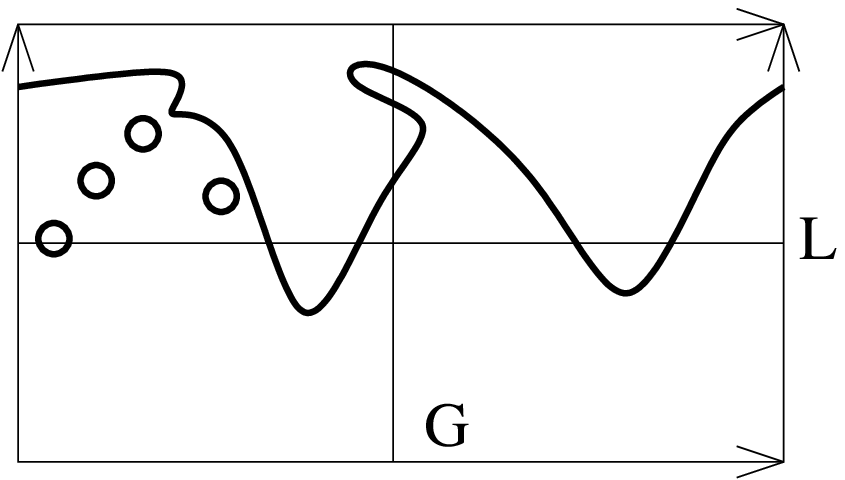}    &
 \includegraphics[height=2cm, angle=0]{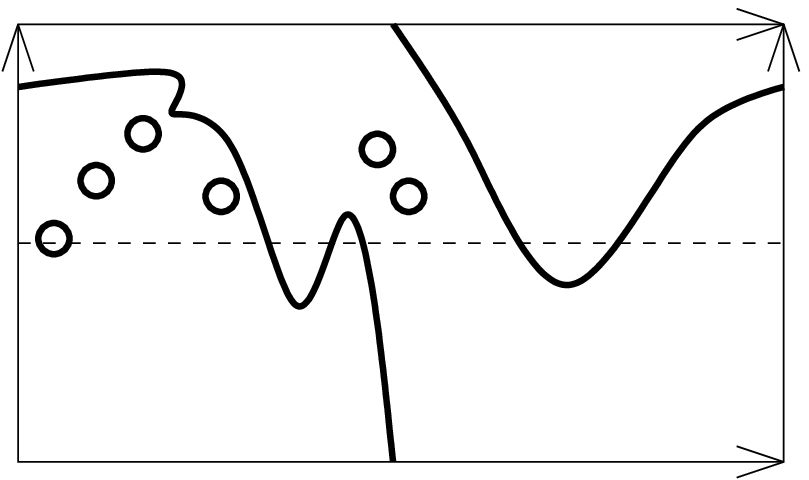}
\\ a)&b)
 \end{tabular}
\caption{}
 \label{rational3}
\end{figure}

\begin{cor}\label{restr thm 4 4}
The complex scheme $\langle J\amalg 5\amalg 1\langle 7\rangle \rangle_{II}$ is realizable by  nonsingular non-dividing symmetric real algebraic curves of degree $7$ on $\mathbb RP^2$.
\end{cor}
\textit{Proof. }The strict transform of the curve constructed in Proposition \ref{par} under the blow up
of $\mathbb CP^2$ at  the point $[0:1:0]$ is the rational real algebraic curve of bidegree $(3,1)$ on $\Sigma_1$
depicted in Figure \ref{rational}b). Blowing up the point $q$ and
blowing down the strict transform of the fiber, we obtain the rational
real algebraic trigonal curve
on $\Sigma_2$ depicted on Figure \ref{rational}c), with a singular point of type $A_6$ at the point $r$. Then according to Shustin's results
(see \cite{Sh1} and \cite{Sh2}), it is possible to smooth this curve
as depicted in Figure \ref{rational3}a). Perturbing the union of this curve and the fiber $G$, we obtain the nonsingular curve of
bidegree $(3,1)$ on $\Sigma_2$ arranged with the base $\{y=0\}$ as shown in Figure \ref{rational3}b). The corresponding symmetric curve realizes the
 real scheme $\langle J\amalg 5\amalg 1\langle 7\rangle \rangle$ and according to Proposition \ref{oval quotient} this is a non-dividing symmetric curve.
~\findemo

\begin{cor}\label{constr thm1 1}
The real schemes $\langle J\amalg 5\amalg 1\langle 8\rangle \rangle$ and
$\langle J\amalg 4\amalg 1\langle 7\rangle \rangle$ are realizable by nonsingular symmetric real algebraic curves of degree $7$ on $\mathbb RP^2$.
\end{cor}
\textit{Proof. }One obtains these two curves modifying slightly the previous construction. To obtain the real scheme $\langle J\amalg 5\amalg 1\langle 8\rangle \rangle,$ one can keep all the ovals above the base depicted in Figure \ref{rational3}b). To obtain the real scheme $\langle J\amalg 4\amalg 1\langle 7\rangle \rangle$, one can consider the line $L$ instead of the line $\{y+z=0\}$ in Figure \ref{rational}a).
~\findemo

\subsection{Change of coordinates in $\Sigma_2$}

\begin{figure}[h]
      \centering
 \begin{tabular}{ccc}
\includegraphics[height=2cm, angle=0]{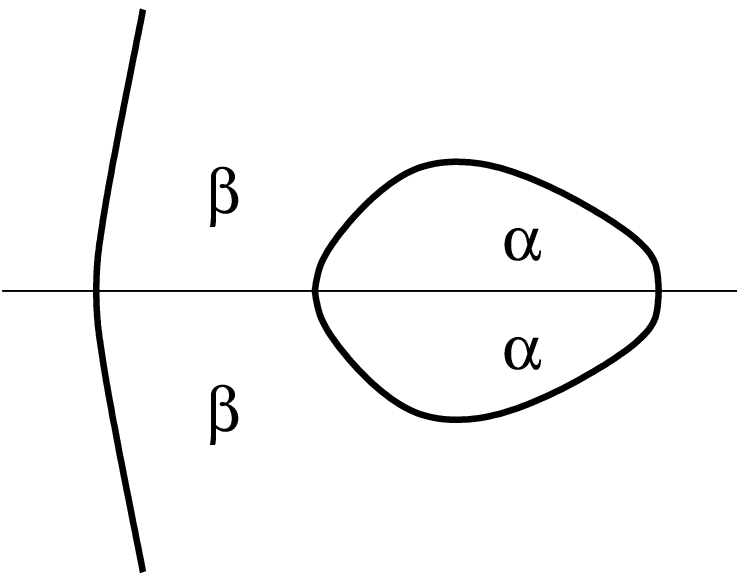}    &
\includegraphics[height=2cm, angle=0]{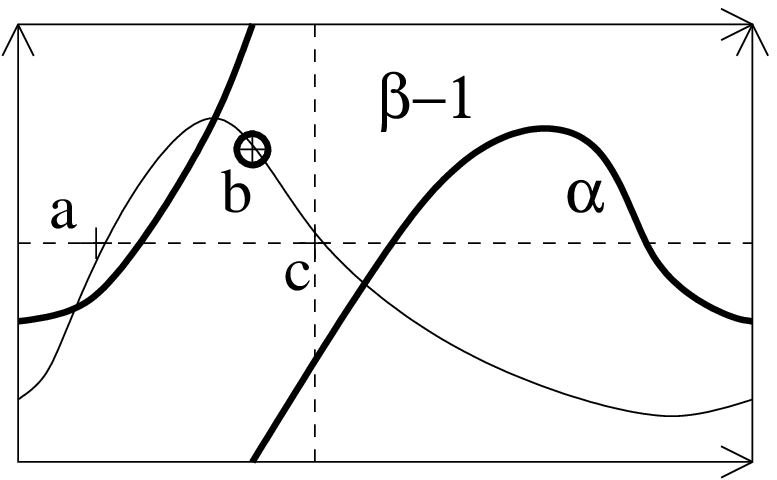}    &   \includegraphics[height=2cm, angle=0]{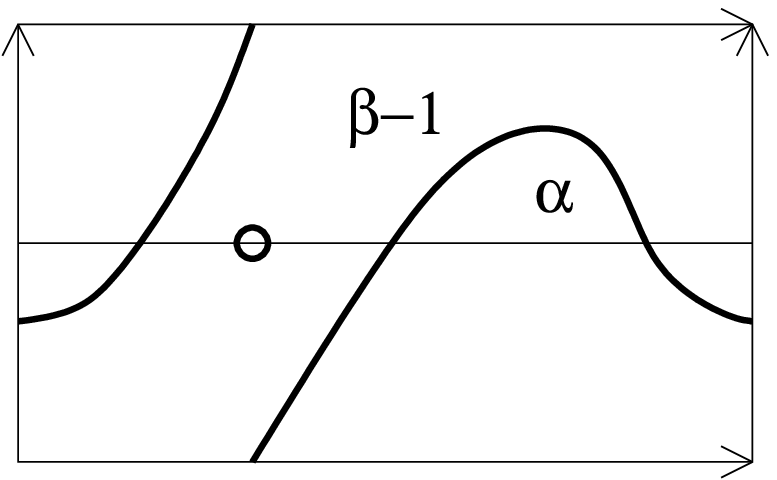}
\\ a)&b)&c)
\\&$(\alpha,\beta)=(2,4)$ or $(3,3)$&
 \end{tabular}
\caption{}
 \label{constr1}
\end{figure}

\begin{prop}\label{change1}
The real schemes $\langle J\amalg 7\amalg 1\langle 4\rangle \rangle$ and $\langle J\amalg 5\amalg 1\langle 6\rangle \rangle $ are realizable by nonsingular symmetric real algebraic curves of degree $7$ on $\mathbb RP^2$.
\end{prop}
\textit{Proof. }
In section \ref{viro}, we constructed symmetric curves on $\mathbb RP^2$  shown in Figure \ref{constr1}a).
According to Lemma \ref{lem bezout1}, their quotient curve $X$ is depicted in Figure \ref{constr1}b).
\begin{figure}[h]
      \centering
 \begin{tabular}{ccc}
\includegraphics[height=2cm, angle=0]{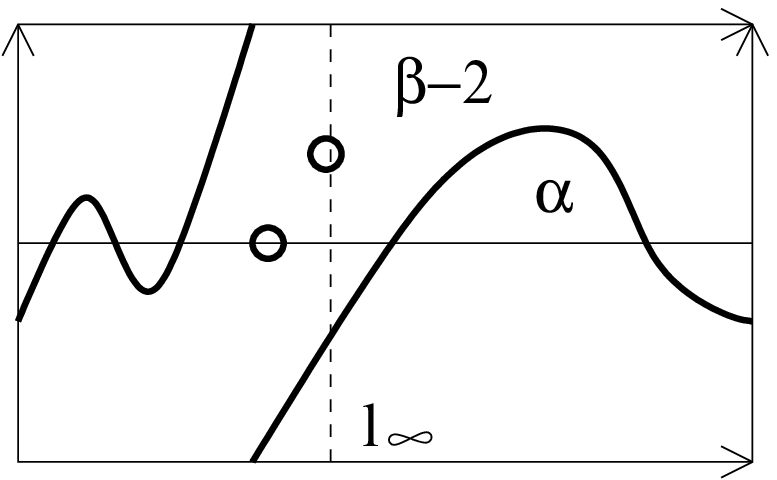}    &
\includegraphics[height=2cm, angle=0]{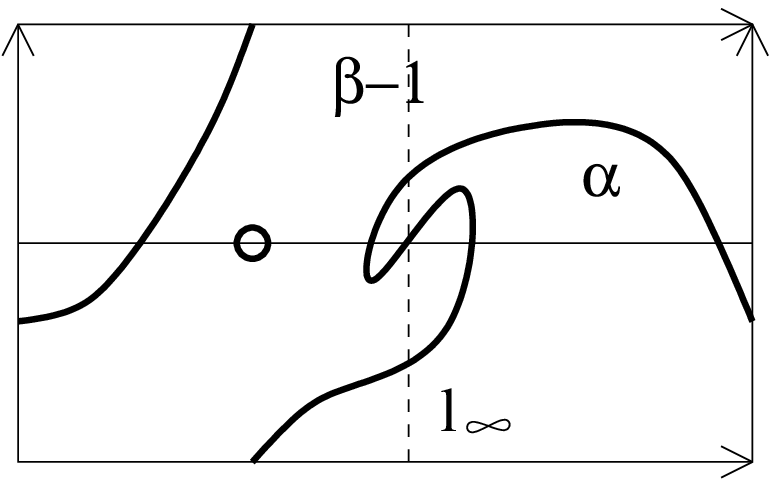}    &   \includegraphics[height=2cm, angle=0]{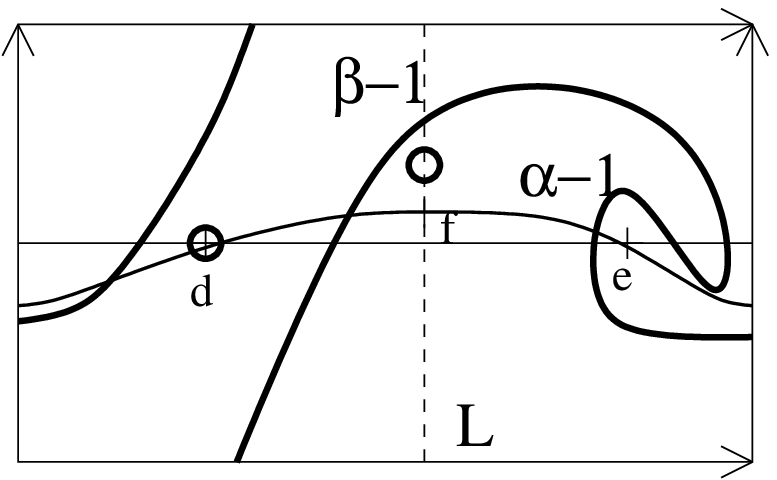}
\\ a)&b)&c)
\\&$(\alpha,\beta)=(2,4)$ or $(3,3)$&
 \end{tabular}
\caption{}
 \label{constr2}
\end{figure}

Let $H$ be the base which passes through the points $a$, $b$ and $c$
as depicted in Figure \ref{constr1}b).
All possible mutual arrangements for $H$ and the quotient curves which do not contradict the Bezout theorem and Lemma \ref{convexe} are depicted in Figures \ref{constr1}c) and \ref{constr2}.

First, we prohibit pseudoholomorphically the $\mathcal L$-schemes realized by the union of $X$ and $H$ in Figure \ref{constr2}a) and by $X$ in Figure \ref{constr2}b). Choose $l_\infty$ as depicted, and the braid corresponding to these $\mathcal L$-schemes are~:
\begin{quote}
$ b^{18}_{\alpha,\beta}= \sigma_3^{-(\beta-1)}\sigma_1^{-1}\sigma_2^{-1}\sigma_3\sigma_2^{-\alpha}\sigma_3^{-1}\sigma_2\sigma_1^{-4}
\sigma_2^{-1}\sigma_3^{2}\sigma_2\sigma_1\sigma_2^{-2}\sigma_3^{-1}\Delta_4^2
  \textrm{ with  }
(\alpha,\beta)=(3,3),(2,4)$,
\\$b^{19}=  \sigma_1^{-7}\sigma_2\sigma_1\Delta_3^2$.
\end{quote}
The braid $b^{18}_{3,3}$ was already shown to be not quasipositive in section \ref{prohib1}.
The computation of the Alexander polynomials of the remaining braids gives
\begin{center}
\begin{tabular}{ll}
$p^{19}=(-1+t)(t^4-t^3+t^2-t+1)$,&
$p^{18}_{4,2}=(t^2-t+1)(-1+t)^3 $.
\end{tabular}
\end{center}
Since $e(b^{19})=1$ and $e(b^{18}_{\alpha,\beta})=2$,  according to Proposition \ref{alex},  none of these braids is quasipositive.

Thus, the two remaining possibilities for the mutual arrangement of $X$ and $H$ are depicted in Figures \ref{constr1}c) and \ref{constr2}c).

In the first case, let $H'=H$.

In the second case,
consider the base $G$
passing through the points $d$, $e$ and $f$ where $f$ is some point on the fiber $L$. For some $f$,
the base $G$ passes through two ovals of $X$. Since $G$ cannot have more than $7$ common points with $X$, there exist $f$ for
which the mutual arrangement of $G$ and $X$ is as shown in Figure \ref{constr1}c). Let $H'$ be the base corresponding to such an $f$.

The symmetric curves of degree $7$ on $\mathbb RP^2$ corresponding to the mutual arrangement of $H'$ and $X$ realize the real schemes
$\langle J\amalg 7\amalg 1\langle 4\rangle \rangle$ and $\langle J\amalg 5\amalg 1\langle 6\rangle \rangle$.
~\findemo

\begin{figure}[h]
      \centering
 \begin{tabular}{ccc}
\includegraphics[height=2cm, angle=0]{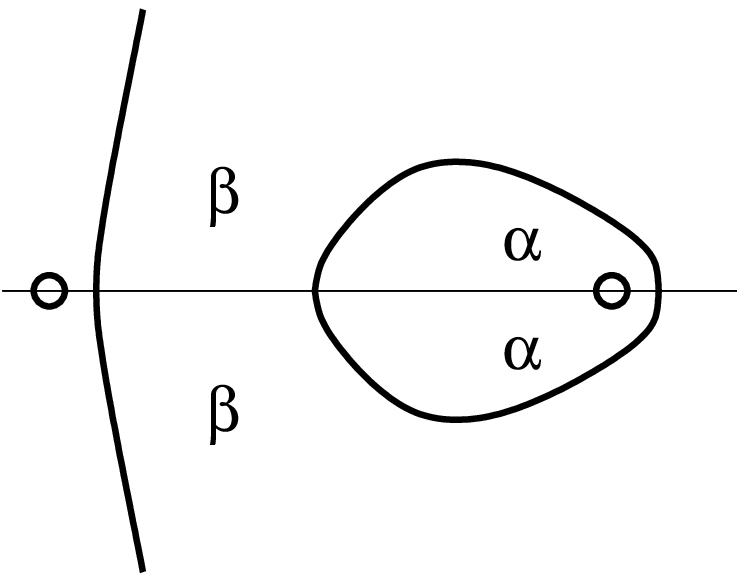}    &
\includegraphics[height=2cm, angle=0]{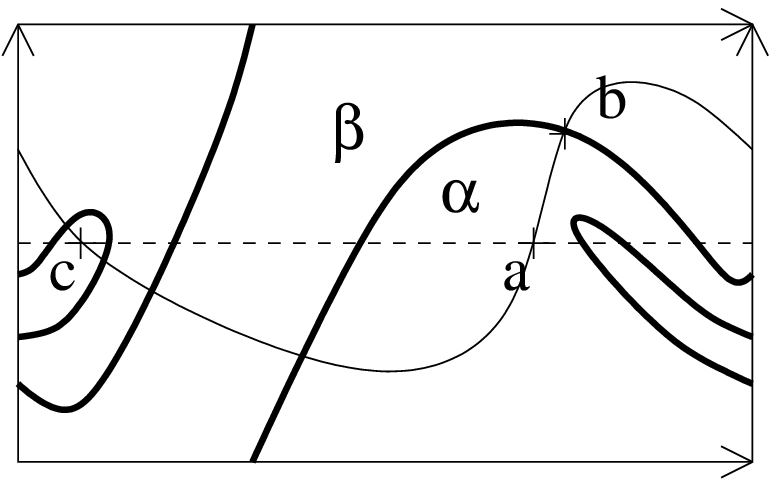}    &   
\includegraphics[height=2cm, angle=0]{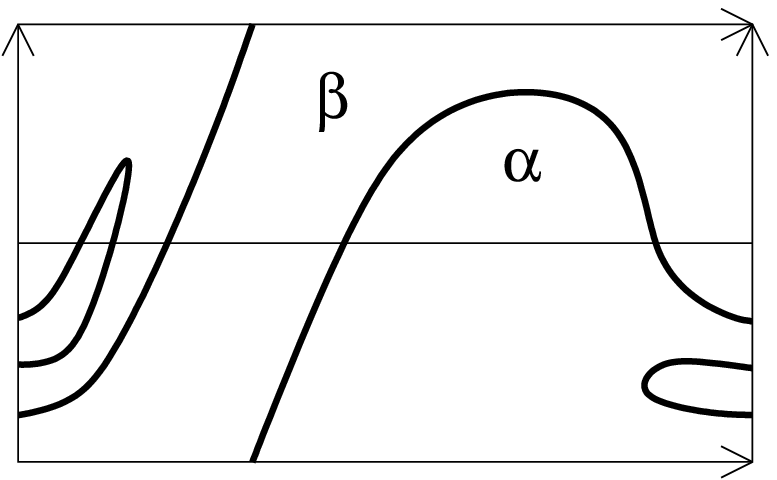}
\\ a)&b)&c)
\end{tabular}

$(\alpha,\beta)=(6,0)$ or $(2,4)$
\caption{}
 \label{constr5}
\end{figure}

\begin{prop}\label{constr thm1 3}
The real schemes $\langle J\amalg 1\amalg 1\langle 12\rangle \rangle$ and $\langle J\amalg 9\amalg 1\langle 4\rangle \rangle $ are realizable by nonsingular symmetric real algebraic curves of degree $7$ on $\mathbb RP^2$.
\end{prop}
\textit{Proof. }In section \ref{viro},
 we constructed symmetric curves of degree $7$ on $\mathbb RP^2$
depicted in Figure \ref{constr5}a).
These are $M$-curves, so according to Proposition \ref{max sym}, \ref{inter max}, the Bezout theorem, Lemmas \ref{lem bezout1} and \ref{convexe},
  the  $\mathcal L$-scheme realized
by their quotient curve $X$ is   depicted in Figure
 \ref{constr5}b).
Let $H$ be the base which passes through the points $a$, $b$ and $c$
as depicted in Figure \ref{constr5}b).
The only possible mutual arrangement for $H$ and $X$ which does not contradict the Bezout theorem, Proposition \ref{max sym} and Lemma \ref{convexe} is depicted in Figure \ref{constr5}c). The corresponding symmetric curves realize the real schemes
$\langle J\amalg 1\amalg 1\langle 12\rangle \rangle$ and
$\langle J\amalg 9\amalg 1\langle 4\rangle \rangle$.~\findemo

\subsection{Construction of auxiliary curves}
Here we construct charts of some symmetric curves. We will use these
charts in
section \ref{last constr} to perturb symmetrically  some singular symmetric
real algebraic curves.

\begin{figure}[h]
      \centering
 \begin{tabular}{cccc}
\includegraphics[width=2.5cm, angle=0]{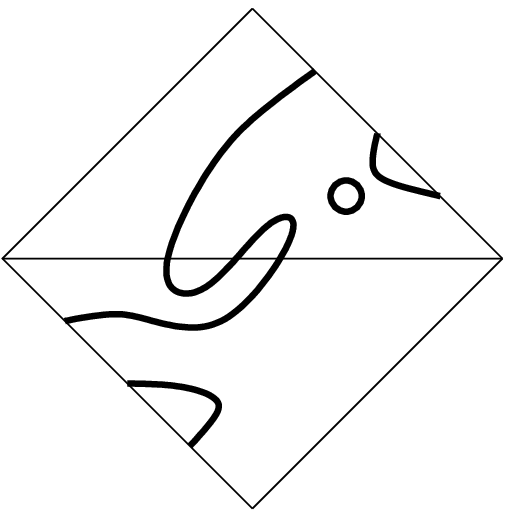}
 &\includegraphics[width=2.5cm,angle=0]{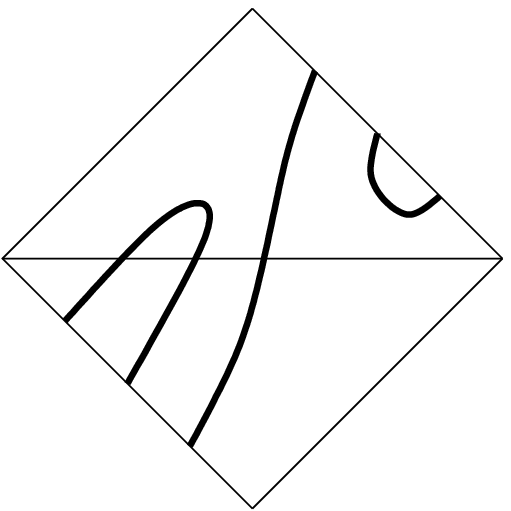}&
 \includegraphics[width=4cm, angle=0]{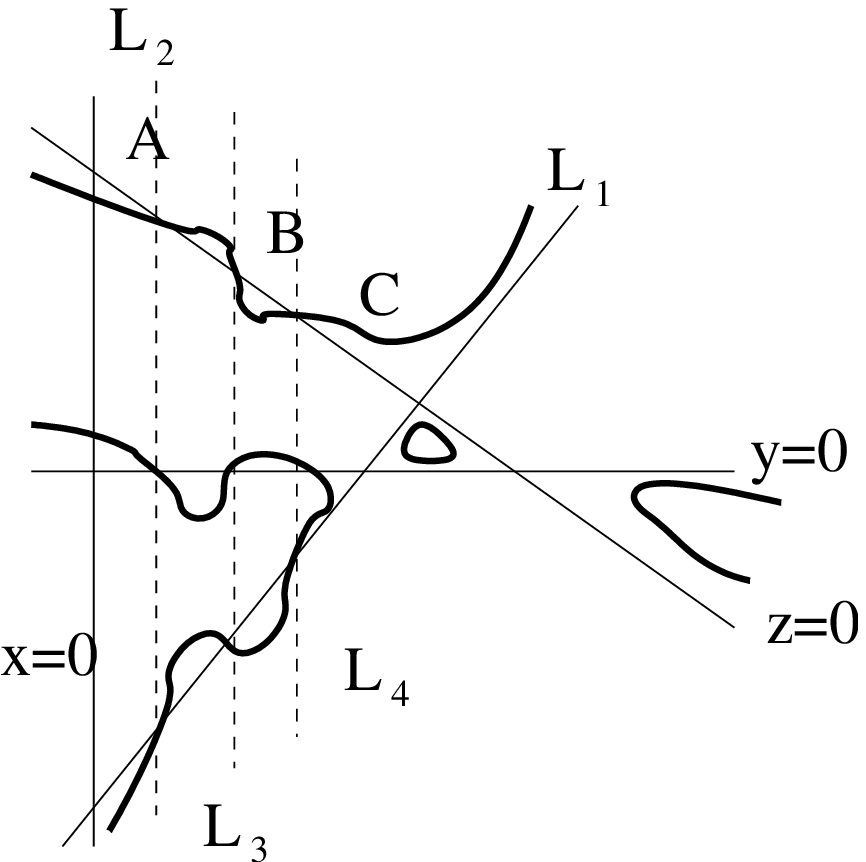}
 &\includegraphics[width=4cm,angle=0]{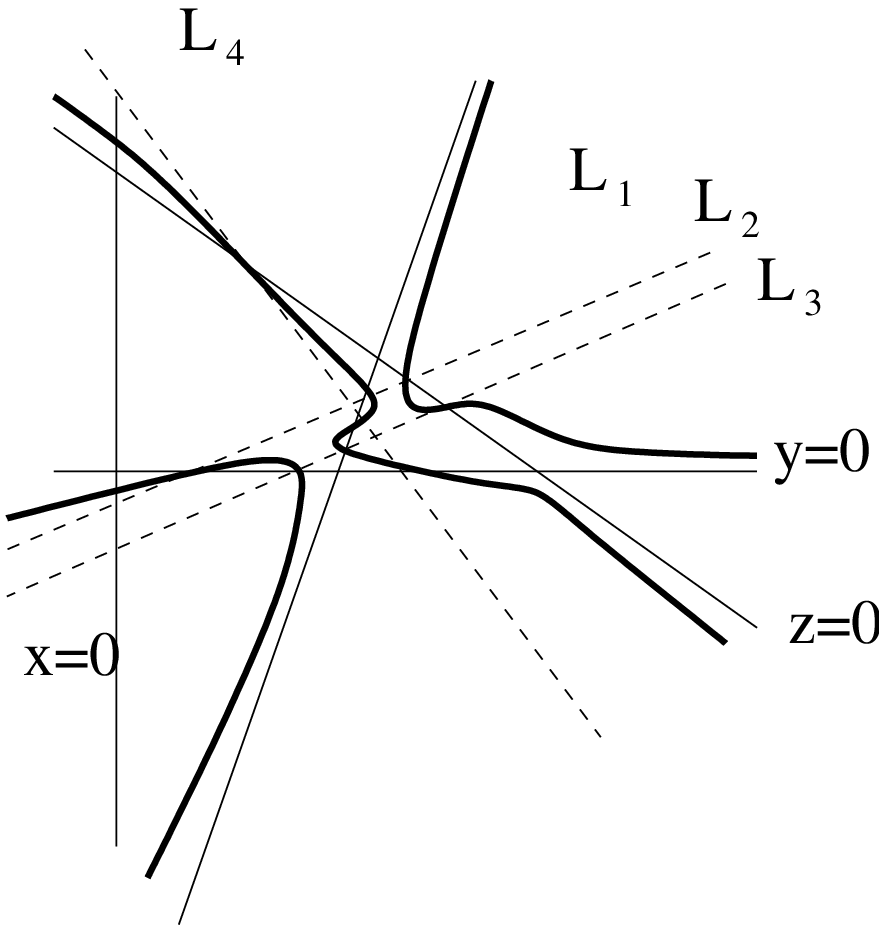}  \\
a)&b)&c)&d)\\
 \end{tabular}
\caption{}
 \label{degree 3 1}
\end{figure}

\begin{lemma}\label{curve 3 1}
For any real positive numbers $\alpha,\beta,\gamma$, there exist real curves of degree $3$  on $\mathbb RP^2$
having the charts and the arrangement
with respect to the axis $\{y=0\}$ shown in Figures \ref{degree 3 1}a) and b) with truncation on the segment $[(0,3),(3,0)]$ equal to
$(x-\alpha y)(x-\beta y)(x-\gamma y)$.
\end{lemma}
\textit{Proof. }Consider the points $A=[\alpha:1:0],B=[\beta:1:0],C=[\gamma:1:0]$ and
four lines $L_1,L_2,L_3,L_4$  as shown in Figure  \ref{degree 3 1}c). For $t$ small enough and of suitable sign,
 the curve $yzL_1+tL_2L_3L_4$
is arranged with respect to the coordinate axis and the lines $L_1,L_2,L_3,L_4$ as shown in Figure  \ref{degree 3 1}c).

To construct the curve with the chart depicted in Figure \ref{degree 3 1}b), we perturb the third degree curve $yzL_1$ as shown in Figure \ref{degree 3 1}d).
~\findemo

\begin{figure}[h]
      \centering
 \begin{tabular}{ccccc}
\includegraphics[width=2cm,  angle=0]{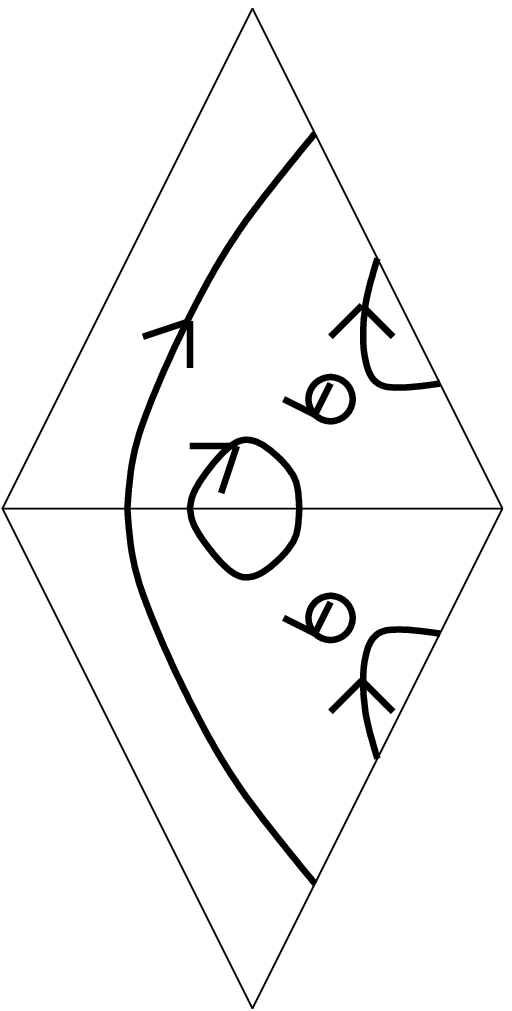}&
\includegraphics[width=2cm,  angle=0]{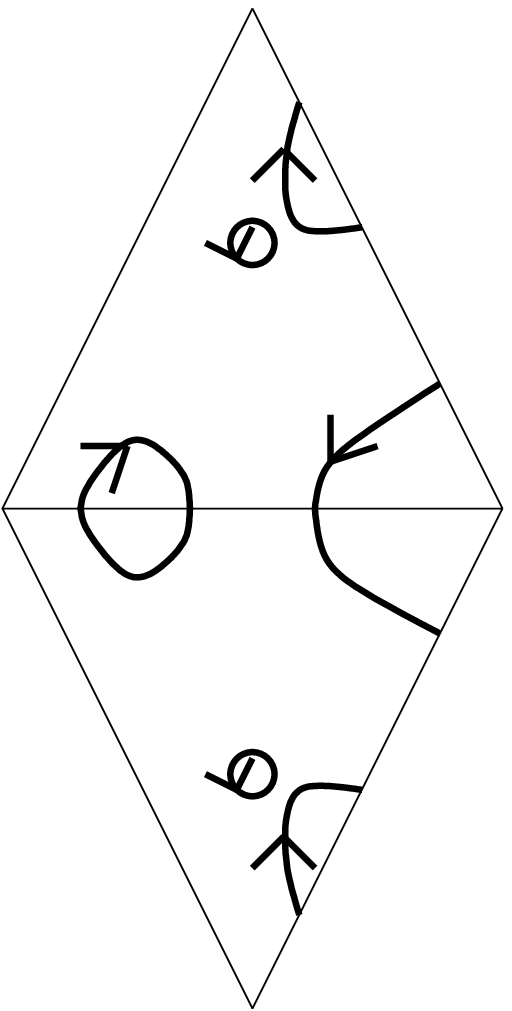}&
\includegraphics[width=2cm,  angle=0]{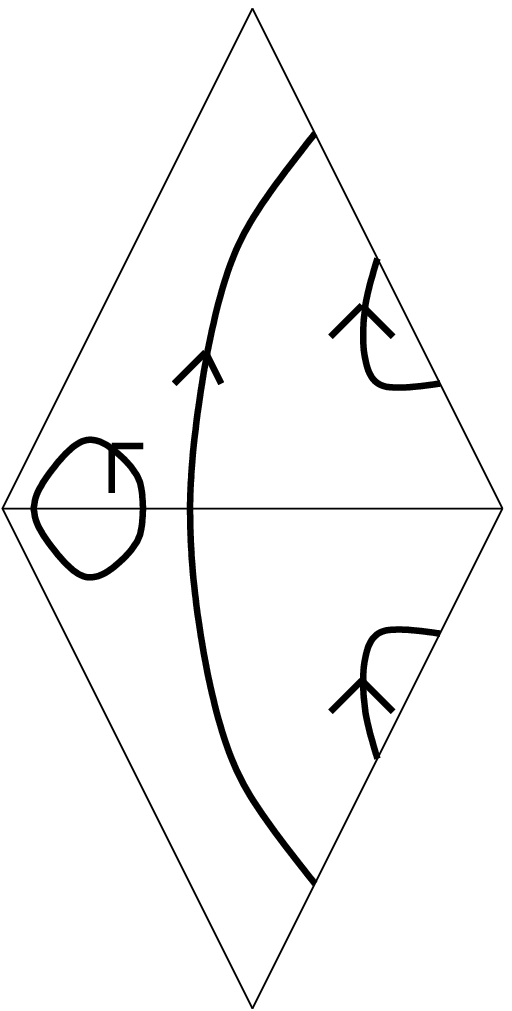}&
\includegraphics[width=2cm,  angle=0]{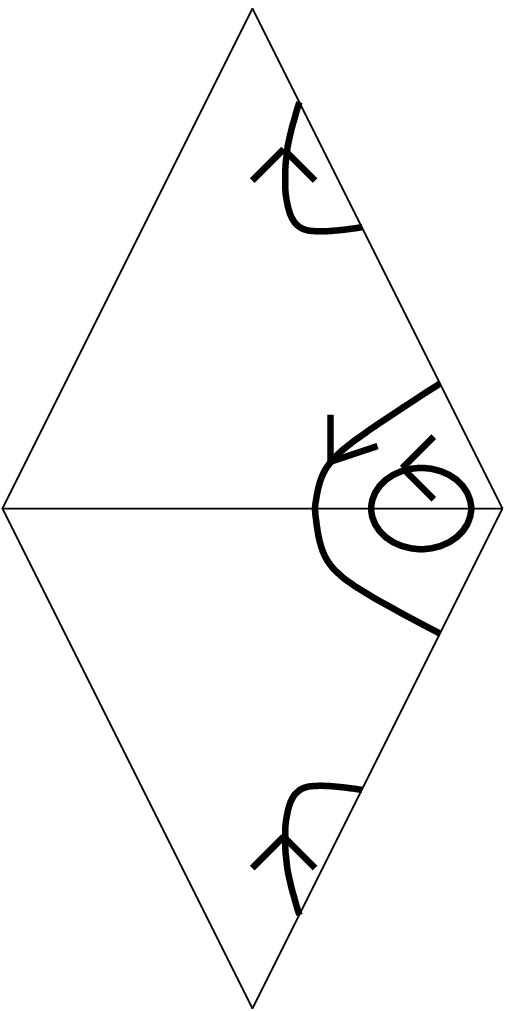}&
\includegraphics[height=4cm,angle=0]{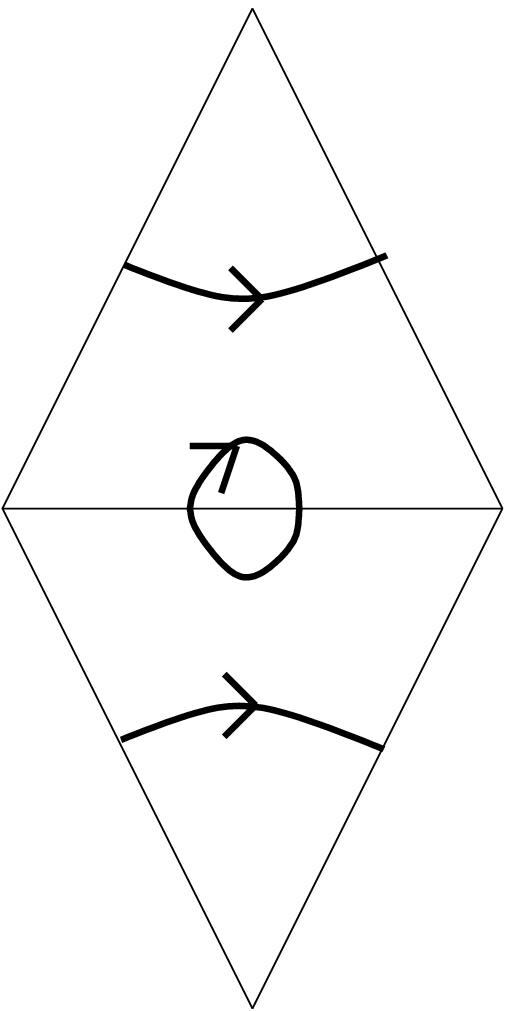}
\\a)&b)&c)&d)&e)
\end{tabular}
\caption{}
\label{degree 6}
\end{figure}

\begin{cor}
For any real positive numbers $\alpha,\beta,\gamma$, there exist  real symmetric dividing curves of degree $6$ on $\mathbb RP^2$
with a singular point of type $J_{10}$ at $[1:0:0]$ having  the charts, the arrangement
with respect to the axis $\{y=0\}$ and the complex orientations
shown in Figures \ref{degree 6}a), b) and c) with truncation on the segment $[(0,3),(6,0)]$ equal to
$(x-\alpha y^2)(x-\beta y^2)(x-\gamma y^2)$.
\end{cor}
\textit{Proof. } The Newton polygon of the third degree curves constructed in Lemma \ref{curve 3 1} lies inside the triangle with vertices
$(0,3),(0,0)$ and $(6,0)$, so these curves can be seen as a (singular)
trigonal curve on $\Sigma_2$. The corresponding symmetric
 curves are of degree $6$ and has the chart and the arrangement with respect to the axis $\{y=0\}$ shown in Figures \ref{degree 6}a) and c). Moreover,
it is well known that such curves are of type I, and we deduce their complex orientations from their quotient curve.
\\ If we perform the coordinate changes $(x,y)\mapsto (-x+\delta
y^2,y)$  with $\delta\in \mathbb R$ to the curves with chart depicted
in Figure \ref{degree 6}a) (resp. \ref{degree 6}c)), we obtain curves
with the chart depicted in Figure \ref{degree 6}b) (resp. \ref{degree 6}d)).
~\findemo

The following lemma can be proved using the same technique.
\begin{lemma}
For any real positive numbers $\alpha$ and $\beta$, there exist real symmetric dividing curves of degree $4$ on $\mathbb RP^2$
with a singular point of type $A_3$ at $[1:0:0]$ having the charts, the arrangement
with respect to the axis $\{y=0\}$ and the complex orientations
shown in Figure \ref{degree 6}e) with truncation on the segment $[(0,2),(4,0)]$ equal to
$(x-\alpha y^2)(x-\beta y^2)$.\findemo
\end{lemma}

\subsection{Perturbation of irreducible singular symmetric
  curves}\label{last constr}

\begin{figure}[h]
      \centering
 \begin{tabular}{cccc}
\includegraphics[width=3.7cm, angle=0]{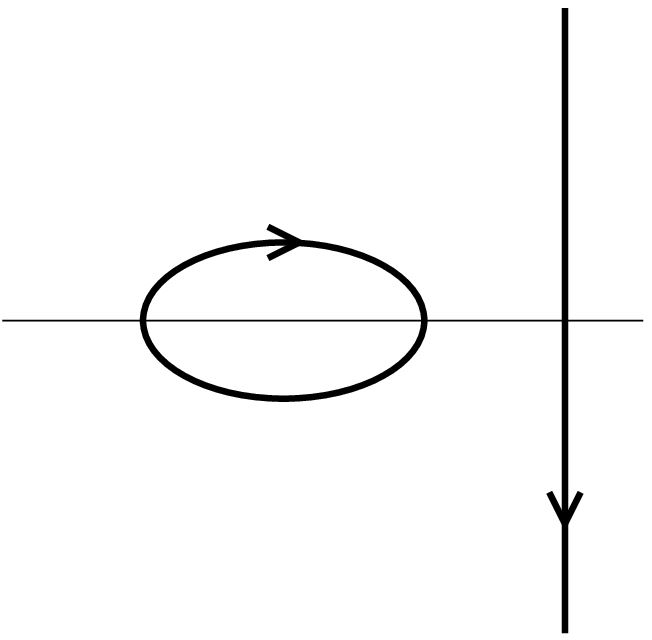}    &
 \includegraphics[width=3.7cm, angle=0]{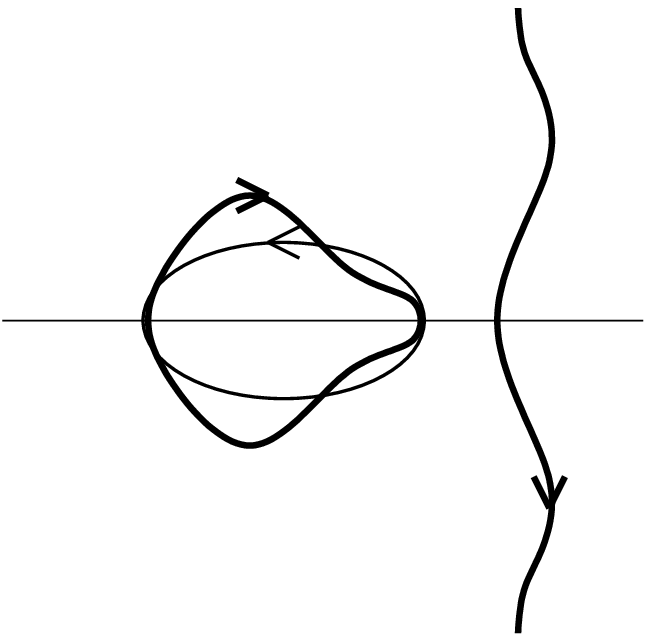} &
\includegraphics[width=3.7cm, angle=0]{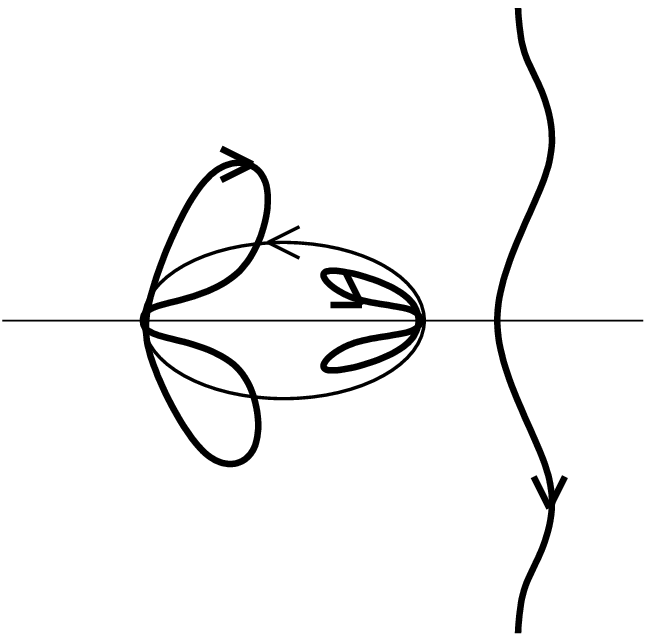}
&\includegraphics[width=3.7cm,angle=0]{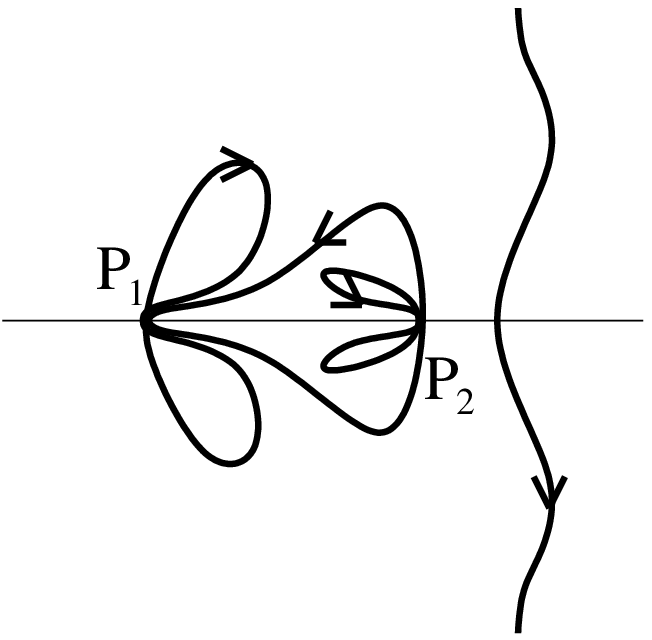}   \\
a)&b)&c)&d)\\
 \end{tabular}
\caption{}
 \label{2_6_I}
\end{figure}

\begin{prop}\label{restr thm 4 6}
The complex schemes $\langle J\amalg 2_-\amalg 1_+\langle 4_+\amalg 2_-\rangle\rangle_I$ and $\langle J\amalg 2_+\amalg 4_-\amalg 1_+\langle 2_+\rangle\rangle_I $ are realizable by nonsingular symmetric real algebraic curves of degree $7$ on $\mathbb RP^2$.
\end{prop}
\textit{Proof. }First, we construct the symmetric singular dividing
curve of degree $7$ with two singular points $J_{10}$ depicted in
Figure \ref{2_6_I}d). 
To construct such a curve, we use the Hilbert method as in
\cite{V1}. Let $C_0$ be a symmetric conic.
We symmetrically perturb the union of $C_0$  and a disjoint
symmetric line 
(Figure \ref{2_6_I}a)) keeping the tangency points with $C_0$.
 Next, we symmetrically perturb the union of the third degree curve
 obtained and $C_0$ (Figure \ref{2_6_I}b)) 
keeping the tangency points of order $4$ with $C_0$. Perturbing
in a symmetric way the union of the curve of degree $5$ 
obtained and $C_0$ (Figure \ref{2_6_I}c)) keeping the
tangency points of order $6$ with $C_0$,  we obtain 
a symmetric singular dividing
curve $C$ of degree $7$ with two singular points $J_{10}$ as depicted in
Figure \ref{2_6_I}d).

Now we symmetrically perturb the singular points of $C$  using the
chart shown in Figure \ref{degree 6}b) (resp. d)) in $P_1$ and
\ref{degree 6}c) (resp. a)) in $P_2$ 
and obtain the desired curves. 
~\findemo

\begin{figure}[h]
      \centering
 \begin{tabular}{cccc}
\includegraphics[height=3cm, angle=0]{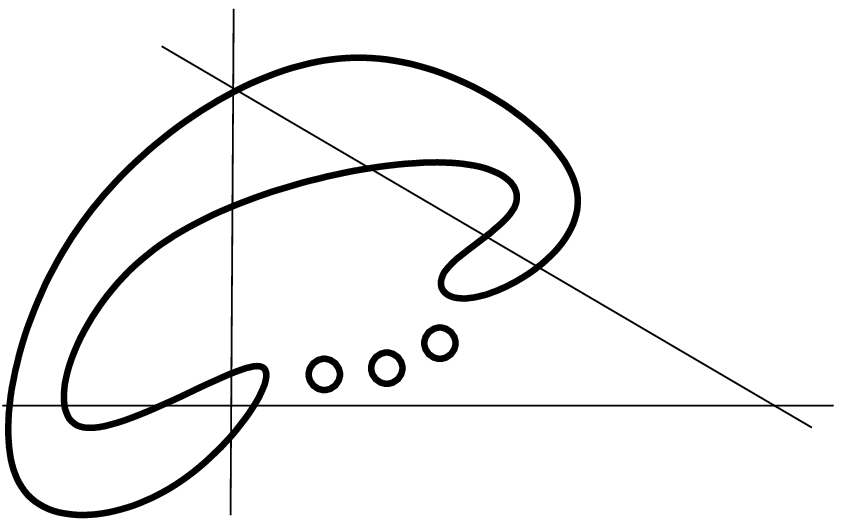}    &
 \includegraphics[height=3cm, angle=0]{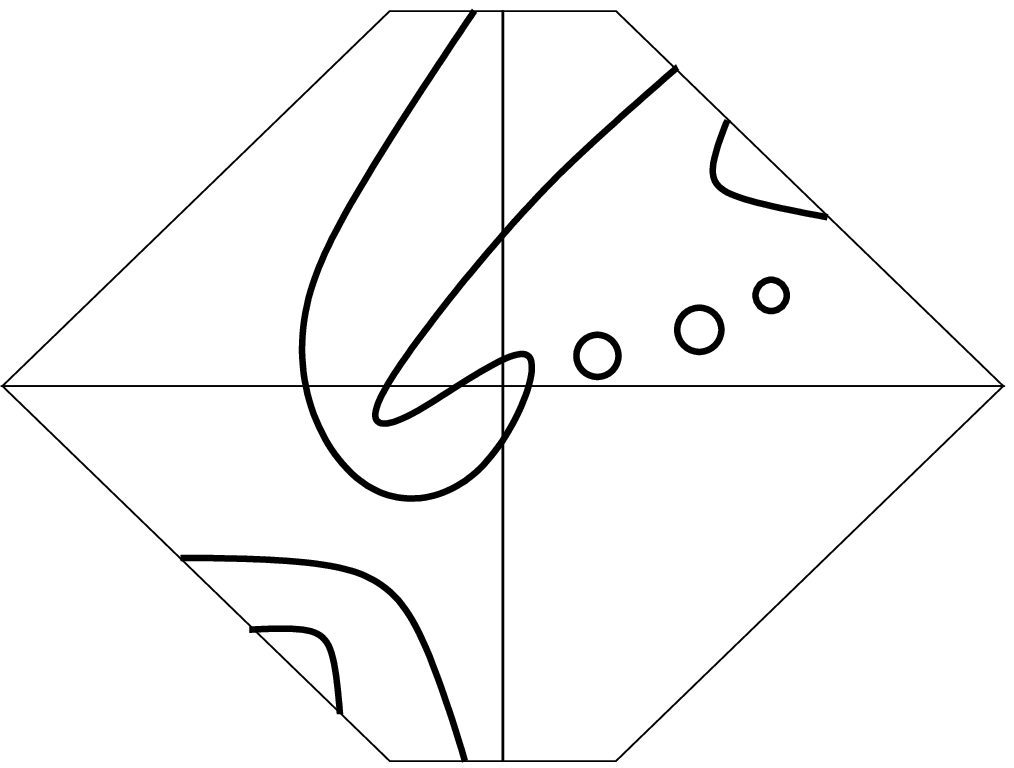} &
\includegraphics[height=3cm, angle=0]{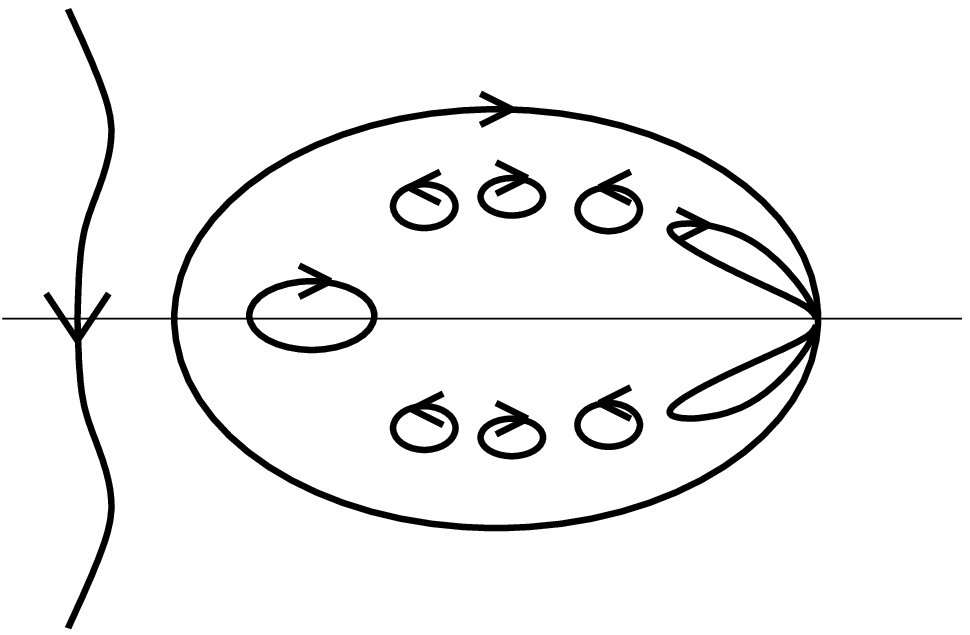}
\\ a)&b)&c)
 \end{tabular}
\caption{}
 \label{12_I}
\end{figure}

\begin{prop}\label{restr thm 4 7}
The complex scheme $\langle J\amalg  1_+\langle 6_+\amalg 6_-\rangle\rangle_I $ is realizable by nonsingular  symmetric real algebraic curves of degree $7$ on $\mathbb RP^2$.
\end{prop}
\textit{Proof. }Consider the curve of degree $4$ with a C-shaped oval constructed in \cite{K1} and the coordinate system shown in Figure \ref{12_I}a). In this
coordinate system, the Newton polygon of the curve is the trapeze with vertices $(0,0),(0,3),(1,3)$ and $(4,0)$, and its chart is depicted in Figure  \ref{12_I}b).
This curve can be seen as a singular real algebraic curve of bidegree $(3,1)$ on the surface $\Sigma_2$. The corresponding symmetric curve of degree $7$
 has a singular point $J_{10}$ at $[1:0:0]$ and is depicted in Figure \ref{12_I}c). 
This curve is maximal according
to Proposition \ref{lattice},
so of type I, and we can deduce its complex orientations from its
quotient curve.
 Finally, we symmetrically smooth the singular point
with the chart depicted in Figure \ref{degree 6}a) and obtain the desired curve.
~\findemo

Denote by $f_P$ the real birational transform of $\mathbb CP^2$ given by
$(x,y)\mapsto (x,y-P(x))$ in the affine coordinate $\{z=1\}$, where $P$ is a polynomial of degree 2.

\begin{figure}[h]
      \centering
 \begin{tabular}{cccc}
\includegraphics[height=2.5cm, angle=0]{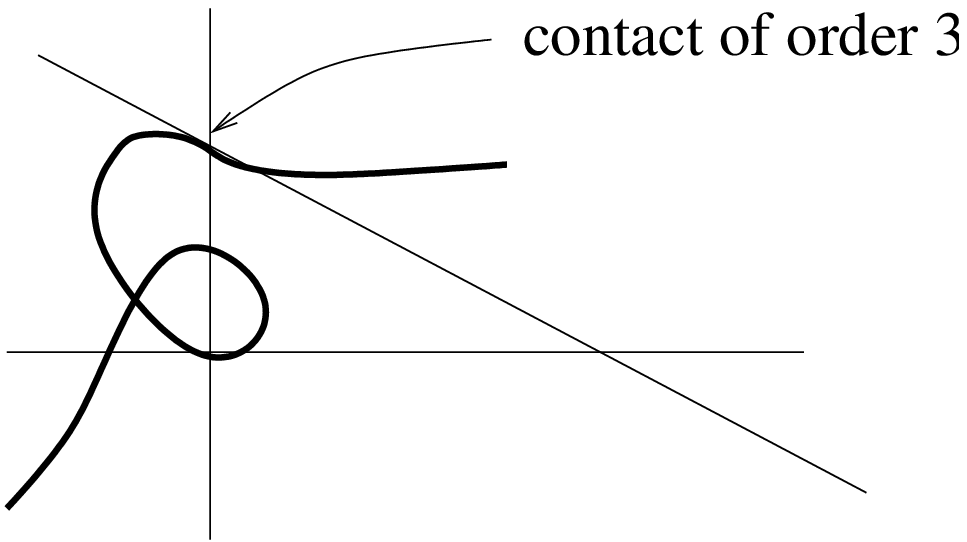}    &
 \includegraphics[height=2.5cm, angle=0]{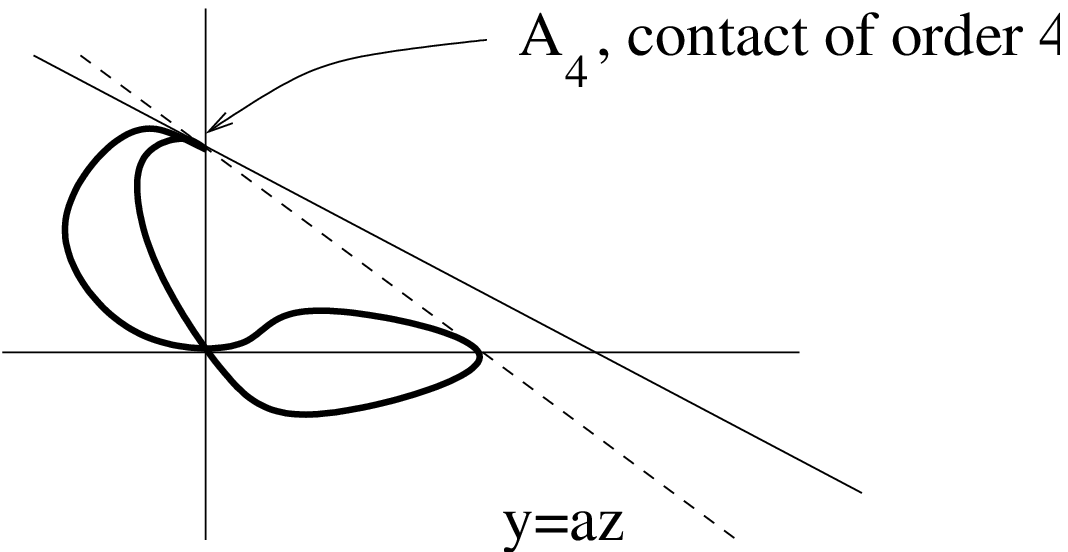} &
\includegraphics[height=2.5cm, angle=0]{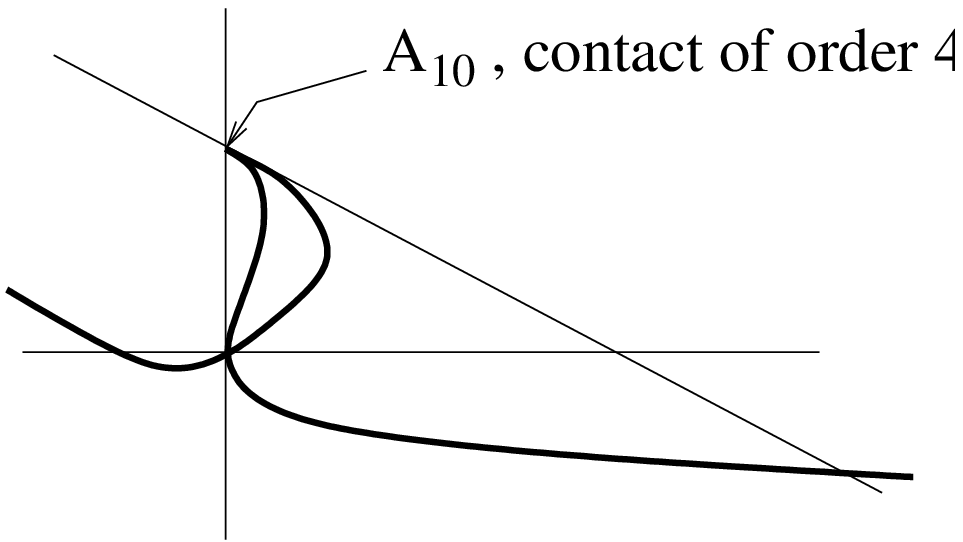}
\\ a)&b)&c)\\
 \end{tabular}
\caption{}
 \label{10_2_IA}
\end{figure}

\begin{figure}[h]
      \centering
 \begin{tabular}{cc}
\includegraphics[height=3cm, angle=0]{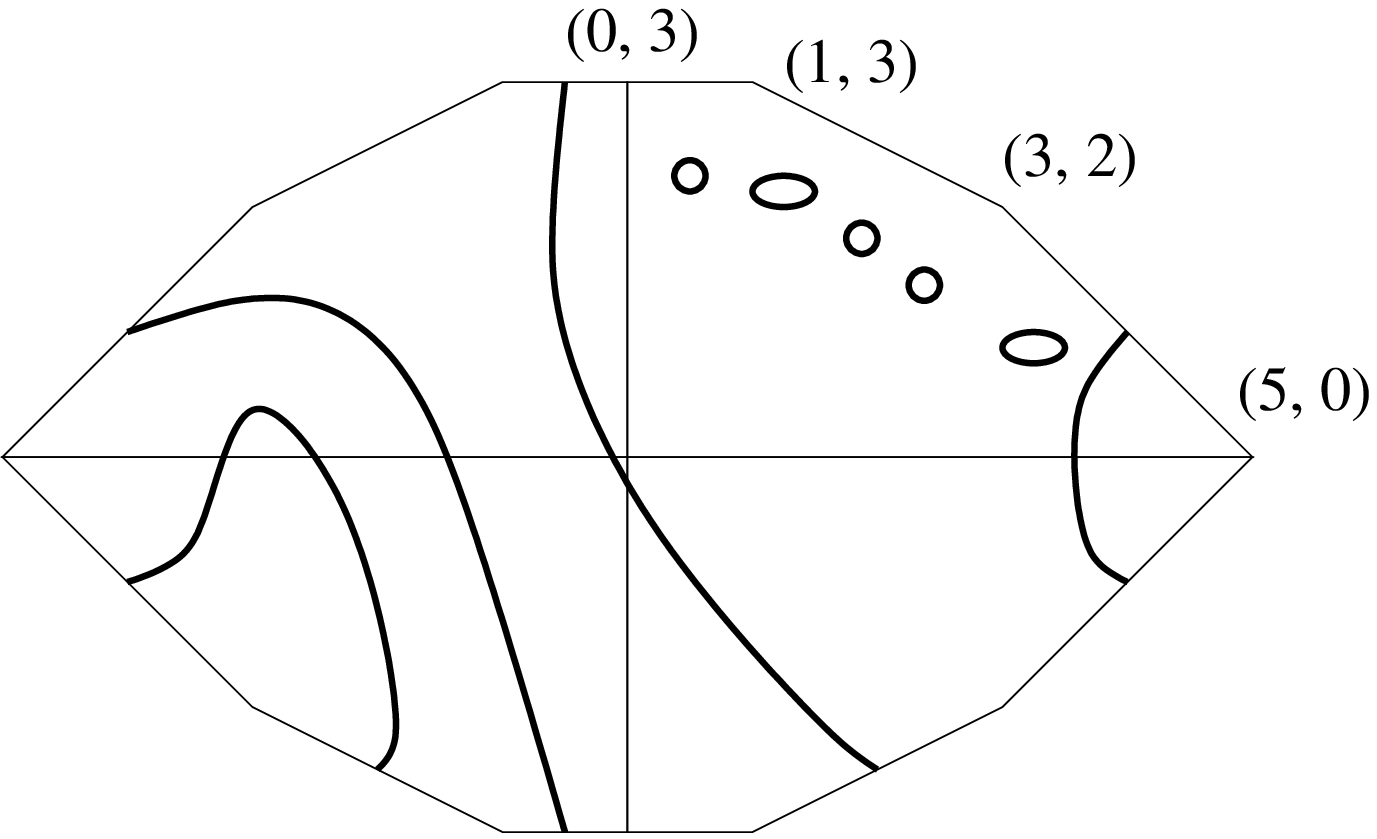}    &
 \includegraphics[height=3cm, angle=0]{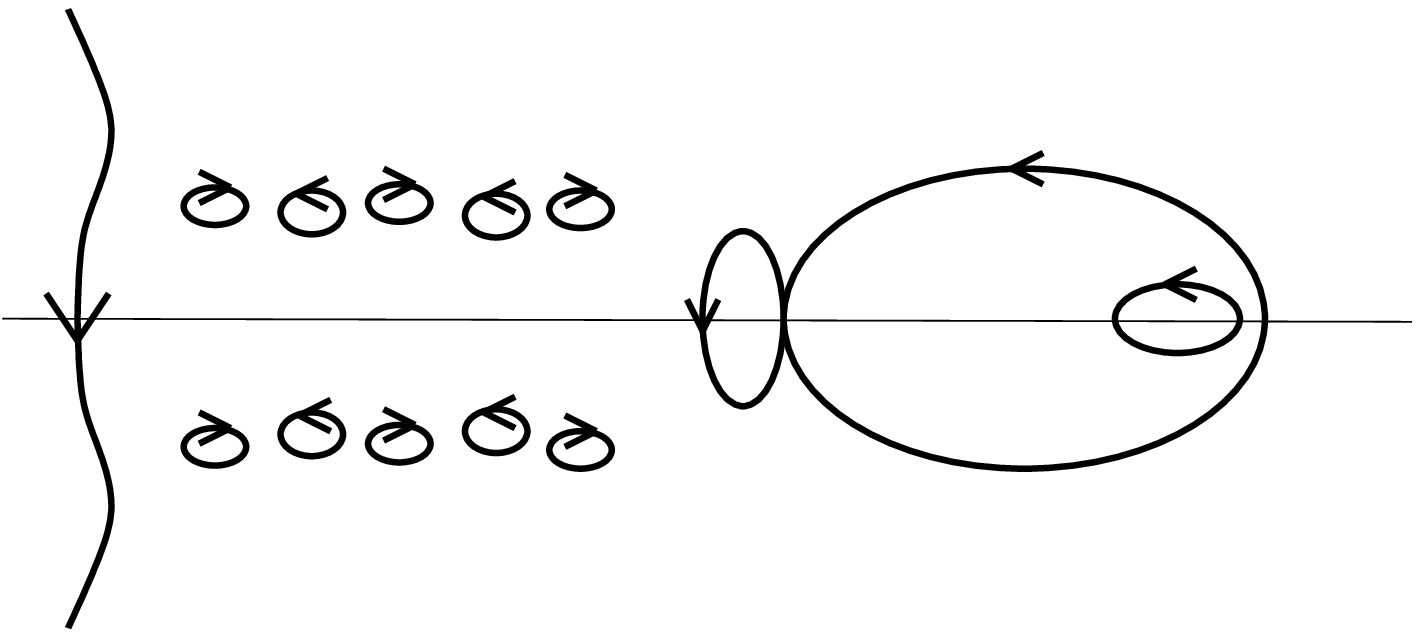}
\\ a)&b)
 \end{tabular}
\caption{}
 \label{10_2_IB}
\end{figure}

\begin{prop}\label{restr thm 4 8}
The complex scheme $\langle J\amalg 6_+\amalg 4_-\amalg  1_+\langle 1_+\amalg 1_-\rangle\rangle_I$ is realizable by nonsingular  symmetric real algebraic curves of degree $7$ on $\mathbb RP^2$.
\end{prop}
\textit{Proof. }Consider the nodal curve of degree $3$ depicted in Figure \ref{10_2_IA}a) with a contact of order $3$ at
the point $[0:1:0]$ with the line $\{z=0\}$.
 Then there exists a unique polynomial $P$ of degree 2 such that the image of the cubic under $f_P$ is the curve of
degree $4$ depicted in Figure \ref{10_2_IA}b), with a singular point of type $A_4$ at $[0:1:0]$ and a contact of order $4$
 at this point with the line $\{z=0\}$. Moreover, the line $\{y=0\}$ intersects the quartic in two points, one
of them is the node, and this line is tangent at one of the local branches at the node
and at the second intersection point, a line $\{y=az\}$ is tangent at the curve of degree $4$. Perform the change of coordinates
of $\mathbb CP^2$ $[x:y:z]\mapsto [y:x:y-az]$.
For this new coordinate system, there exists a polynomial $Q$ such that the image of the quartic under $f_Q$ is the curve of
degree $5$ depicted in Figure \ref{10_2_IA}c), with a singular point of type $A_{10}$ at $[0:1:0]$ and a contact of order $4$
 at this point with the line $\{z=0\}$. Applying the change of coordinates of $\mathbb CP^2$ $[x:y:z]\mapsto [y:z:x]$ and using Shustin's results on the independent perturbations of generalized semi-quasihomogeneous singular points of a curve keeping the same Newton polygon(see \cite{Sh1} and \cite{Sh2}), we can
smooth the singular point $A_{10}$ in order to obtain a curve
with the chart shown in Figure \ref{10_2_IB}a).

Hence, we can see this curve as a  singular curve of bidegree $(3,1)$ on the surface $\Sigma_2$. The corresponding symmetric
 curve of degree $7$ has a singular point $A_{3}$ at $[1:0:0]$ and is depicted in Figure \ref{10_2_IB}b). This curve is maximal according
to Proposition \ref{lattice},
so of type I, and we can deduce its complex orientations from its
 quotient curve. 
Finally, we symmetrically smooth the singular point
with the chart depicted in Figure \ref{degree 6}e) and obtain the desired curve.
~\findemo

\small
\def\rightmark{\em Bibliography}
\addcontentsline{toc}{section}{References}

\bibliographystyle{alpha}
\bibliography{../../Biblio.bib}

\hspace*{2 ex}
\textbf{\\Erwan Brugall\'e}\\
Université de Paris 6\\
Projet Analyse algébrique\\
175 rue du Chevaleret\\
75013 Paris\\
France\\
\\
E-mail : brugalle@math.jussieu.fr
\end{document}